\documentclass[10pt,leqno]{amsart}
\usepackage[T1]{fontenc}
\usepackage[utf8]{inputenc}
\usepackage{lmodern,amsmath,amssymb,amsthm,mathtools}
\usepackage[a4paper,margin=25mm]{geometry}
\usepackage{microtype,enumitem,booktabs,array}
\usepackage[colorlinks=true,linkcolor=blue,citecolor=blue,urlcolor=blue]{hyperref}
\usepackage{xcolor}
\hypersetup{pdftitle={Hodge--Riemann relations for m-positive and generalized m-positive forms},pdfauthor={Enhan Li}}
\numberwithin{equation}{section}
\newtheorem{theorem}{Theorem}[section]
\newtheorem{proposition}[theorem]{Proposition}
\newtheorem{lemma}[theorem]{Lemma}
\newtheorem{corollary}[theorem]{Corollary}
\newtheorem{conjecture}[theorem]{Conjecture}

\newtheorem{mainthm}{Theorem}

\theoremstyle{definition}
\newtheorem{definition}[theorem]{Definition}
\newtheorem{remark}[theorem]{Remark}
\newtheorem{example}[theorem]{Example}
\newcommand{\R}{\mathbb R}\newcommand{\C}{\mathbb C}
\newcommand{\wtG}{\widetilde\Gamma}\newcommand{\vol}{\operatorname{vol}}
\newcommand{\In}{\operatorname{In}}
\newcommand{\diag}{\operatorname{diag}}\newcommand{\tr}{\operatorname{tr}}

\setlist[enumerate]{label=(\alph*),leftmargin=2.2em,itemsep=3pt}
\allowdisplaybreaks[2]
\title[Hodge--Riemann relations and m-positivity]{Hodge--Riemann relations for $m$-positive and generalized $m$-positive forms}
\author{Enhan Li} 
\address{Department of Mathematics, University of Science and Technology of China}
\email{lablehustc0719@mail.ustc.edu.cn}
\date{}
\newcommand{\ii}{\sqrt{-1}}
\newcommand{\Ann}{\operatorname{Ann}}
\newcommand{\rank}{\operatorname{rank}}

\newcommand{\rad}{\operatorname{rad}}

\begin{document}
\begin{abstract}
Our goal in this paper is to show strong positivity and Hodge--Riemann relations for real $(1,1)$-forms with various generalizations of $m$-positivity. For a good prefix, which we define using the hyperbolicity cones of iterative derivatives of the determinant, we show that the wedge product is strictly strongly positive, and mixed Hodge--Riemann relations hold in all bidegrees. We also determine the hyperbolicity cone of every nonzero product of an almost-good prefix as the polarizing component of its natural exterior algebra quotient. This yields the mixed version of Hodge--Riemann relations for $m$-positive forms, without the assumption of semipositivity. For any weakly $(m,l)$-positive form, i.e. a form which is $(m-l)$-positive on any subspace of codimension $l$, we show that Hodge--Riemann relations hold in bidegrees $(p,q)$ with $p+q \leq m$ and $\max\{p,q\} \geq l$. For $l=1$ and $m \ge 2$, the $(m-1)$-th power map provides a diffeomorphism between the weakly $(m,1)$-positive cone and the positive $(d-1,d-1)$-cone. For smooth closed forms on a compact K\"ahler manifold satisfying the ordinary, weak $(m,l)$-positive, or good-prefix hypotheses pointwise, we use the method of Dinh--Nguy\^en to prove hard Lefschetz theorems and Hodge--Riemann relations. We also treat generalized Hessian equations and associated intersection inequalities, and prove necessary numerical conditions for sequential positivity.
\end{abstract}
\maketitle
\setcounter{tocdepth}{1}
\tableofcontents

\section{Introduction}

%What positivity conditions of a real $(1,1)$-form imply that its wedge products are strongly positive and satisfy Hodge--Riemann relations?
Let $V$ be a $d$-dimensional complex vector space with real positive $(1,1)$-forms $\omega,\omega_1,\dots,\omega_d$. By definition, the product $\Omega_k=\omega_1\wedge\dotsm\wedge\omega_k$ is strictly strongly positive for any $k\le d$. The mixed Hodge--Riemann relation by Timorin \cite{timorin} asserts that multiplication by $\Omega_{p+q}$ induces an isomorphism 
\[
 \Lambda^{p,q}V^*\xrightarrow{\ \wedge \Omega_{p+q}\ }
 \Lambda^{d-q,d-p}V^*,
\]
and that the hermitian form
\begin{gather*}
 Q_{\Omega_{p+q}}\colon \Lambda^{p,q}V^*\times\Lambda^{p,q}V^*\to \mathbb{C},\\    
 Q_{\Omega_{p+q}}(\alpha,\beta)
 = (\sqrt{-1})^{p-q}(-1)^{(p+q)(p+q-1)/2}
       *\bigl(\Omega\wedge\alpha\wedge\overline\beta\bigr)
\end{gather*}
is positive definite on the primitive space 
\[
 \ker\bigl(\Omega_{p+q}\wedge\omega\wedge-\colon
 \Lambda^{p,q}V^*\longrightarrow \Lambda^{d-q,d-p}V^*).
\]

It is natural to ask what under which positivity conditions on a sequence of real $(1,1)$ forms on $V$ does their product satisfy strongly positive and satisfy Hodge--Riemann relations.

We call a real $(1,1)$-form $\eta$ $m$-positive in the ordinary sense with respect to $\omega$ if
\[
 \eta^j\wedge\omega^{d-j}>0\qquad(1\le j\le m).
\]
When $m<d$, ordinary $m$-positive forms can have negative eigenvalues, but they must have at least $m$ positive ones. We extend the hard Lefschetz theorem and Hodge--Riemann bilinear relations to other notions of positivity.

Vanishing theorems on complex manifolds, involving some sort of partial positivity of real $(1,1)$-forms, have been obtained by Andreotti and Grauert \cite{andreotti-grauert}, Xiaokui Yang \cite{yang-ag} and Dinew and Popovici \cite{dinew-popovici}. Demailly \cite{demailly-positivity} and Qilin Yang \cite{yang-ks} also studied partial positivity of forms on Hermitian vector bundles. Many versions of partial positivity have been defined but all of them are based on bounding from below the number of positive eigenvalues of curvature forms. 

In this work we study generalizations of the $m$-positivity of real $(1,1)$-forms. The first generalization is based on the concept of hyperbolic polynomials, of which the $m$-th elementary symmetric polynomial $\sigma_m$ is an example. We call a sequence of real $(1,1)$-forms sequentially positive if each form in the sequence belongs to the hyperbolicity cone of the preceding ones. It is well known in the theory of hyperbolic polynomials \cite{garding} that consecutive cones contain the ones before them. 

Ordinary $m$-positivity also has generalizations related to subsolutions of elliptic PDEs. A real $(1,1)$-form $\eta$ is weakly $(m,l)$-positive with respect to $\omega$ if for any subspace $W\subseteq V$ of codimension $l$, $\eta|_W$ is $(m-l)$-positive with respect to $\omega|_W$. For $l=1$ we obtain weak positivity. Weak $m$-positivity of $\chi+\ii\partial\bar\partial\underline u$ is the $\mathcal C$-subsolution condition for $\underline u$ in the $m$-th Hessian equation on a compact Hermitian manifold with metric $\omega$ \cite{szekelyhidi,trudinger-hessian}. The cone of weak $(m,l)$-positive forms is convex and contains the ordinary $m$-positive cone. In both cases $\eta$ has at least $m$ positive eigenvalues. For sequential positivity we assume that $\eta$ belongs to the hyperbolicity cone of some nonzero almost-good product of length $d-m$, and for weak $(m,l)$-positivity we assume that $l<m$ (Propositions~\ref{seq:positive-index} and~\ref{ord:positive-index}).

\subsection{Sequential positivity and polarization}
For a prefix $(\eta_1,\ldots,\eta_r)$, put
$\Theta_0=1$ and $\Theta_j=\eta_1\wedge\cdots\wedge\eta_j$.
For $1\le j\le r$, let $P_{\Theta_{j-1}}(A)=*\bigl(\Theta_{j-1}\wedge A^{d-j+1}\bigr)$, where $*$ means we extract the coefficient of a top-degree form relative to $\omega^d/d!$. Let $\Gamma_{\Theta_{j-1}}(\omega)$ be the connected component of $\{P_{\Theta_{j-1}}\ne 0\}$ containing $\omega$. We say the prefix is sequentially positive, or good, if for all $j$, $\eta_j\in\Gamma_{\Theta_{j-1}}(\omega)$. We say the prefix is sequentially almost positive, or almost good, if for all $j$, $\eta_j\in\overline{\Gamma_{\Theta_{j-1}}(\omega)}$ and $\Theta_j\ne0$.
The first open cone, $\Gamma_1(\omega)$, is the positive definite cone. Taking derivatives of the determinant repeatedly yields hyperbolic polynomials $P_{\Theta_{j-1}}$. For $j\ge 2$, the cone corresponding to $P_{\Theta_{j-1}}$ contains the cone of $P_{\Theta_{j-2}}$. If $\Theta_{j-1}=\omega^{j-1}$, then the open-cone condition on $\eta_j$ is ordinary $(d-j+1)$-positivity. For any nonzero preceding product $\Theta_{j-1}$ of an almost-good prefix, $\eta_j$ belongs to the open cone (resp. its closure) if and only if
\[
 \Theta_{j-1}\wedge\eta_j^a\wedge\omega^{d-j+1-a}>0
 \qquad(1\le a\le d-j+1)
\]
(resp.\ $\Theta_{j-1}\wedge\eta_j^a\wedge\omega^{d-j+1-a}\ge0$ for the same range of $a$).
Whether a sequence is good or almost good does not depend on the choice of the positive form $\omega$.

By hyperbolicity (Proposition~\ref{seq:restriction}), the product $\Theta_r$ of an almost good prefix is weakly positive. If the prefix is good, $\Theta_r$ is strictly weakly positive. Our main theorem states that $\Theta_r$ is in fact strongly positive.

\begin{mainthm}[Strong positivity]\label{main:strong}
The product of every almost-good prefix is strongly positive.
The product of every good prefix is strictly strongly positive.
In particular, if $(\eta_1,\ldots,\eta_r)$ is good, then for some
$\varepsilon>0$ the form
\[
 \eta_1\wedge\cdots\wedge\eta_r-\varepsilon\omega^r
\]
is a finite sum of products of $r$ semipositive rank-one $(1,1)$-forms.
\end{mainthm}

The proof reduces the rank of the first factor while varying the
later factors at the same time. Each partial product is affine in
the parameter, and the highest positivity condition at every step
is preserved. A rank-one first factor reduces the dimension by
restriction to a hyperplane. This proves strong positivity by
induction. Averaging small sign perturbations of the factors of a good prefix proves strict strong positivity.

The classical hard Lefschetz theorem and Hodge--Riemann bilinear relations are formulated for a K\"ahler class on a compact K\"ahler manifold \cite{voisin,huybrechts}. Their mixed versions replace a power of a K\"ahler class by a product of several K\"ahler classes. Timorin proved linear algebra versions of mixed relations in \cite{timorin}. On compact K\"ahler manifolds, mixed relations were proved by Dinh and Nguy\^en in \cite{dn} using $L^2$-methods. Cattani \cite{cattani} also gave a proof using polarized Hodge--Lefschetz modules and a theorem of Cattani, Kaplan and Schmid \cite{cks}.

In recent years there have been extensions of both the hard Lefschetz theorem and Hodge--Riemann relations for forms of other types of positivity. Xiao proved the mixed Hodge index theorem for $m$-positive forms \cite[Theorem~A]{xiao-index} and mixed Hodge--Riemann relations in higher degrees assuming the forms are semipositive \cite[Theorems~3.1, A]{xiao-mixed}. Hu and Xiao also proved mixed Hodge--Riemann relations for semipositive forms satisfying a Hall--Rado condition \cite[Proposition~3.2]{hu-xiao-hl}. Zhang proved a mixed Hodge index theorem for sequentially positive forms and high degree relations for semipositive forms \cite[Theorems 1.7, 1.10]{zhang}.

Sommese studied $(d-m)$-ample vector bundles \cite{sommese}. De Cataldo--Migliorini proved the hard Lefschetz theorem and Hodge--Riemann relations for lef line bundles \cite[Theorem~2.3.1]{dCM}, and Hu--Shang--Xiao proved mixed relations for free line bundles under $m$-lefness hypotheses \cite[Theorem~A]{hu-shang-xiao}. Li--Zheng and Hu proved related partial Lefschetz results in \cite{li-zheng,hu-kernels}. The $L^2$ method of Dinh and Nguy\^en \cite{dn,dn-abstract} passes from the Hodge--Riemann relations on the exterior algebra to those on the cohomology of compact K\"ahler manifolds.

For a real $(d-p-q,d-p-q)$-form $\Omega$, write
\[
 Q_\Omega(\alpha,\beta)
 = (\sqrt{-1})^{p-q}(-1)^{(p+q)(p+q-1)/2}
       *\bigl(\Omega\wedge\alpha\wedge\overline\beta\bigr),
\]
and, for a real $(1,1)$-form $\lambda$, set
\[
 P^{p,q}(\Omega,\lambda)=\ker(\Omega\wedge\lambda\wedge-).
\]
We use integration to define the cohomological pairing. We analyze the hard Lefschetz isomorphism, positivity on the primitive subspace and the orthogonal Lefschetz decomposition.

Let $\Theta$ be the nonzero product of an almost-good prefix
of length $r\le d$. Let $n=d-r$. For $r=d$, the product of forms is a positive volume form, and we define the hyperbolicity cone of the resulting constant polynomial to be the space of real $(1,1)$-forms. The quotient
\[
 M_\Theta=\Lambda^{*,*}V^*/\operatorname{Ann}(\Theta)
\]
has the Poincar\'e pairing induced by
$(\alpha,\beta)\mapsto *\bigl(\Theta\wedge\alpha\wedge\beta\bigr)$.
Write $\operatorname{Pol}_\Theta(\omega)$ for the connected component,
containing $\omega$, of the real $(1,1)$-forms that polarize this
Hodge--Lefschetz module, i.e., that satisfy the hard Lefschetz property and the Hodge--Riemann relations on $M_{\Theta}$ in all bidegrees. From the definition of the hyperbolicity cone we have $\operatorname{Pol}_{\Theta}(\omega)\subseteq\Gamma_{\Theta}(\omega)$. Our next main result is that these two sets are equal.

\begin{mainthm}[Hyperbolicity and polarizing cones]\label{main:blocks}
For every nonzero almost-good-prefix product $\Theta$, the module
$M_\Theta$ is polarized by the positive definite cone, and
\[
 \Gamma_\Theta(\omega)=\operatorname{Pol}_\Theta(\omega).
\]
More precisely, let $p,q\ge0$, $p+q\le n$, and
$\beta_1,\ldots,\beta_{n-p-q},\lambda\in\Gamma_\Theta(\omega)$.
Set $C=\beta_1\wedge\cdots\wedge\beta_{n-p-q}$. Then
multiplication by $C$ induces an isomorphism
\[
 M_\Theta^{p,q}\xrightarrow{\ \wedge C\ }
 M_\Theta^{n-q,n-p},
\]
and $Q_{\Theta\wedge C}$ is positive definite on
\[
 \ker\bigl(C\wedge\lambda\wedge-\colon
 M_\Theta^{p,q}\longrightarrow M_\Theta^{n-q+1,n-p+1}\bigr).
\]
The corresponding orthogonal Lefschetz decomposition holds.
For a good prefix, multiplication by $\Theta\wedge C$ gives an
isomorphism $\Lambda^{p,q}V^*\to\Lambda^{d-q,d-p}V^*$, with the
Hodge--Riemann relations on its $\lambda$-primitive subspace.
\end{mainthm}

With a good prefix we obtain positivity in bidegree $(0,0)$. Theorem~\ref{main:strong} implies Hermitian positivity and equivalently the Hodge--Riemann relations in bidegree $(0,q)$. The Hodge--Riemann relations on $(1,1)$ are the Hodge index theorem and follow from Zhang \cite{zhang} (see also Proposition~\ref{seq:low}). Proposition~\ref{seq:moore} gives another proof in bidegree $(0,2)$ using the Moore determinant. 

With good prefixes Theorem~\ref{main:blocks} recovers some known mixed Hodge--Riemann relations for forms. Theorem~\ref{glob:theorem} gives the cohomological version. The classical mixed K\"ahler relation is the case of the empty prefix. Taking $\Theta = \omega^{d-m}$ gives Xiao’s Hodge index and semipositivity results as $\Gamma_{\Theta} = \Gamma_m (\omega)$. Successively taking blocks of factors in the cone of the preceding product gives Zhang’s mixed relations \cite[Theorems~1.7 and~1.10]{zhang}.

For almost good prefixes Theorem~\ref{main:blocks} strengthens, for exterior-algebra quotients, the descent lemma of Cattani, Kaplan and Schmid \cite[Lemma~1.16]{cks}. If $\eta \in \overline{\Gamma_\Theta}$ and $\Theta \wedge \eta \neq 0$, descent says $\Gamma_\Theta$ polarizes $M_{\Theta \wedge \eta}$ \cite[Theorem~3.2 and Remark~3.3]{cattani}. After descent Theorem~\ref{main:blocks} gives the polarizing component.
\[
 \Gamma_\Theta(\omega)\subseteq\Gamma_{\Theta\wedge\eta}(\omega)
 =\operatorname{Pol}_{\Theta\wedge\eta}(\omega).
\]

The initial step in the proof is to consider tuples $(\eta_1,\dots,\eta_r)$ with $\eta_1\ne0$ semipositive and $\eta_j\in\Gamma_{\Theta_{j-1}}$ for $2\le j\le r$. Using rank deformation, we can adjust the signature. For kernel determination we restrict to hyperplanes, varying the last term of the tuple. For more general boundary terms, we can find successive kernels by using central quotients of a polarized Hodge--Lefschetz module \cite{ahl}. To determine the signs in the perturbation, we can lift the kernels. Both arguments are contained in Section~\ref{seq:proof-section}, which also contains a proof of strong positivity.

The cone comparison concerns the quotient of the exterior algebra.
Hyperbolicity of a cohomological volume polynomial alone does not
give higher Lefschetz isomorphisms; Appendix~\ref{app:volume}
gives an example. In Section~\ref{seq:weak-section} we compare positivity and weak positivity with respect to the same prefix.

\subsection{Ordinary \texorpdfstring{$m$}{m}-positivity}
Taking $\Theta = \omega^{d-m}$ in Theorem~\ref{main:blocks} gives the following corollary which confirms a prediction of Xiao \cite[Subsection 3.4]{xiao-index}.

\begin{corollary}\label{main:ordinary}
Let $1 \le p+q \le m \le d$ and let $\eta_1, \ldots, \eta_{m-p-q+1}$ be real $(1,1)$ forms, $m$-positive with respect to $\omega$. Let
\[
 \Omega=\omega^{d-m}\wedge\eta_1\wedge\cdots\wedge\eta_{m-p-q},\qquad
 \lambda=\eta_{m-p-q+1}.
\]
Then multiplication by $\Omega$ is an isomorphism
\[
 \Lambda^{p,q}V^*\longrightarrow\Lambda^{d-q,d-p}V^*,
\]
and $Q_\Omega$ is positive definite on $P^{p,q}(\Omega,\lambda)$. There is also an orthogonal decomposition
\[
 \Lambda^{p,q}V^*=P^{p,q}(\Omega,\lambda)
                   \oplus\lambda\wedge\Lambda^{p-1,q-1}V^*.
\]
The same is true in the cohomology of a compact K\"ahler manifold for classes represented by smooth closed forms that are pointwise $m$-positive with respect to a common K\"ahler form $\omega$.
\end{corollary}

In the case where all the forms are equal to an $m$-positive form $\eta$, we first prove the hard Lefschetz theorem. When $m<d$, Lemma~\ref{lem:goodhyperplanes} gives an orthonormal family of hyperplanes $\{H_i\}_{i=1}^d$ in $V$ such that each restriction $\eta|_{H_i}$ is $m$-positive. From the primitive relations on $H_i$, we obtain injectivity. Using homotopy we obtain the Hodge--Riemann relations.

Equivalently, multiplication by $\eta$ polarizes the Hodge--Lefschetz module
\[
 \Lambda^{*,*}V^*/\operatorname{Ann}(\omega^{d-m}).
\]
Since the cone of $m$-positive forms is convex, we can use Cattani's mixed theorem on this module. This gives a direct proof of Corollary~\ref{main:ordinary}, without the semipositivity assumption in Xiao's theorem \cite[Theorem~A]{xiao-mixed}.

\subsection{Weak \texorpdfstring{$(m,l)$}{(m,l)}-positivity}
Let $0 \leq l \leq m \leq d$. We call a real $(1,1)$-form $\eta$ weakly $(m,l)$-positive with respect to $\omega$ if its restriction to any complex subspace of codimension $l$ is $(m-l)$-positive with respect to the restriction of $\omega$. Positivity of order 0 is always satisfied. For eigenvalue coordinates, the cone is
\[
 C_{m,l}^d
 =\bigcap_{\substack{I\subseteq\{1,\ldots,d\}\\|I|=l}}
   \{a\in\R^d:a_{\widehat I}\in\Gamma_{m-l}^{d-l}\}.
\]
For $m \geq 1$, we have $C_{m,0}^d=\Gamma_m^d$, which is the $m$-positive cone, and $C_{m,1}^d=\widetilde\Gamma_m^d$, which is the weak $m$-positive cone. For $m \geq 2$, the weak $m$-positive cone is also described by the $(m-1)$-positive forms satisfying the inequality
\[
 \omega^{d-m}\wedge\eta^{m-1}>0
\]
as a $(d-1,d-1)$-form.

For $l=1$, weak $m$-positivity is the $\mathcal C$-subsolution condition for the complex $m$-Hessian equation. Let $\chi$ be a real $(1,1)$-form and $h$ be a positive function. Consider the equation
\[
 \sigma_m\bigl(\lambda_\omega(\chi+\ii\partial\bar\partial u)\bigr)=h
 \qquad\text{on }\Gamma_m.
\]
If we set $a=\lambda_\omega(\chi+\ii\partial\bar\partial\underline u)$ for the eigenvalues at a point, we can apply Sz\'ekelyhidi’s criterion \cite[Section~2]{szekelyhidi} to each ray $a+t e_i$. Such a ray enters $\Gamma_m$ if and only if $a_{\widehat i}\in\Gamma_{m-1}^{d-1}$. Then $\sigma_m(a+t e_i)=\sigma_m(a)+t\sigma_{m-1}(a_{\widehat i})$ tends to infinity as $t\to+\infty$. Thus $\underline u$ is a $\mathcal C$-subsolution if and only if $a\in C_{m,1}^d$ at all points.

The cones $C_{m,l}^d$ are open and convex. For $m\ge 1$ they are the higher derivative cones $\mathcal C_{\sigma_m}^{(l)}$ defined by Fang and Ma in \emph{G\aa rding polynomials} \cite[Definition 8.4]{fang-ma-garding}.
\[
 \partial_I\sigma_m(a)=\sigma_{m-l}(a_{\widehat I})\qquad(|I|=l).
\]

\begin{mainthm}[Hodge--Riemann relations for weak $(m,l)$-positive forms]\label{main:weak}
Let $\eta$ be weakly $(m,l)$-positive with respect to $\omega$, where $0\le l\le m\le d$. For $p,q\ge0$ satisfying
\[
 p+q\le m,\qquad \max\{p,q\}\ge l,
\]
put $\Omega_{p,q}=\omega^{d-m}\wedge\eta^{m-p-q}$. Then multiplication by $\Omega_{p,q}$ is an isomorphism
\[
 \Lambda^{p,q}V^*\longrightarrow\Lambda^{d-q,d-p}V^*.
\]
The Hodge--Riemann form $Q_{\Omega_{p,q}}$ is positive definite on $P^{p,q}(\Omega_{p,q},\omega)$, and
\[
 \Lambda^{p,q}V^*
 =P^{p,q}(\Omega_{p,q},\omega)
   \oplus\omega\wedge\Lambda^{p-1,q-1}V^*
\]
is an orthogonal decomposition. These claims hold true on compact K\"ahler cohomology for smooth closed $\eta$ satisfying the condition pointwise with respect to a K\"ahler form $\omega$.
\end{mainthm}

Thus, if $\eta$ is weakly $m$-positive, then the hard Lefschetz theorem, the Hodge--Riemann relations and the primitive decomposition hold in all positive total degrees $1 \le p+q \le m$, with primitive subspaces defined by $\omega$. In general, $\eta$ cannot be used in place of $\omega$ in the primitive condition as weakly $m$-positive forms are not necessarily positive with respect to the Hodge index pairing.

The proof proceeds by induction on $l$ and on the dimension. After decomposition of the exterior algebra by the diagonal torus we reduce the problem to certain higher Hessian matrices. This is then reduced to the case of equal positive eigenvalues and equal negative eigenvalues. The remaining determinant is given by a rectangular Schur polynomial and a Jacobi integral. We evaluate this integral to be a positive finite sum using Gaussian quadrature with the signs controlled by the cone condition. We give the proof in Section~\ref{sec:lorentz}.

Theorem~\ref{main:weak} with $l=1$ and $(p,q)=(1,1)$ asserts the signature in the statement below. We show global invertibility in Section~\ref{sec:gradient}.

\begin{corollary}\label{main:gradient}
Let $2\le m\le d$. For each $\eta\in\widetilde\Gamma_m(\omega)$, the real quadratic form
\[
 \alpha\longmapsto *\bigl(\omega^{d-m}\wedge\eta^{m-2}\wedge\alpha^2\bigr)
\]
has one positive and $d^2-1$ negative eigenvalues. The power map
\[
 \eta\longmapsto\omega^{d-m}\wedge\eta^{m-1}
\]
is a real-analytic diffeomorphism from $\widetilde\Gamma_m(\omega)$ to the cone of strictly positive $(d-1,d-1)$-forms.
\end{corollary}

We use Maclaurin’s inequalities to prove properness of the power map. This will establish that the local diffeomorphism is a global covering of one sheet, as the domain is connected and the target, the cone of strictly positive $(d-1,d-1)$ forms, is convex. This pointwise question was motivated by Fu and Xiao’s study of the map from the K\"ahler cone to the balanced cone \cite{fu-xiao}.

Weak positivity also has a mixed version defined by subsolutions. For a good prefix $(\eta_1,\ldots,\eta_{d-m})$, with $m\ge2$, let $\Theta=\eta_1\wedge\cdots\wedge\eta_{d-m}$, and we require that
\[
 \Theta\wedge\omega^{m-1-j}\wedge\eta^j>0
 \qquad(1\le j\le m-1)
\]
as $(d-1,d-1)$-forms. By Proposition~\ref{seq:weak-hyperplanes}, this is the same as $\eta|_H\in\Gamma_{\Theta|_H}(\omega|_H)$ for every complex hyperplane $H$. In general, in the mixed setting the quadratic form $\alpha\mapsto*(\Theta\wedge\eta^{m-2}\wedge\alpha^2)$ may not be Lorentzian. One example of this failure is presented in Example~\ref{seq:weak-failure} on $\C^5$ and $(\mathbb P^1)^5$.

\subsection{Globalization and numerical characterization}
We follow Dinh and Nguy\^en \cite{dn,dn-abstract} to establish cohomology relations from relations on points. In bidegree $(p,q)$, we assume that the pairs $(\Omega\wedge\lambda^{2j},\lambda)$ satisfy the pointwise Lefschetz and Hodge--Riemann properties in bidegree $(p-j,q-j)$ for $0\le j\le\min(p,q)$ for some smooth closed forms $\Omega$ and $\lambda$. Any globally good prefix also has the relations in each bidegree. Here, a globally good prefix is a sequence of smooth closed real $(1,1)$-forms that forms a good prefix at every point. Compact K\"ahler manifolds are assumed to be connected. To give a case where the global Hodge--Riemann relations apply, we provide the following criterion for the existence of an $m$-positive form in any given real $(1,1)$ class. The proof uses the numerical criterion of Fu--Li--Zhang--Zhang \cite[Theorem~10.1]{flzz}, as shown in Proposition~\ref{glob:complete-intersection}.

\begin{proposition}\label{intro:complete-intersection}
Let $X$ be a smooth projective manifold of dimension $d$, let $H$ be a
very ample line bundle, and let $1\le m\le d$. Assume that a smooth
closed real $(1,1)$-form $\eta$ restricts to a K\"ahler form on the smooth
complete intersection of $d-m$ generic hypersurfaces in $|H|$. Then $[\eta]$ contains
an $m$-positive representative with respect to every K\"ahler form in
$c_1(H)$.
\end{proposition}

Subsection~\ref{hess:section} discusses the equation
\begin{equation}\label{eq:generalizedhessian} \Theta\wedge\bigl(\theta+\sqrt{-1}\partial\bar\partial u\bigr)^m=e^F\omega^d,
 \qquad m=d-\deg\Theta,
\end{equation}
for a fixed good-prefix product $\Theta$. By hyperbolicity, the operator 
\[u\longmapsto \frac{\Theta\wedge(\theta+\sqrt{-1}\partial\bar\partial u)^m}{\omega^d}\]
is elliptic and its $m$-th root is concave on the admissible cone. By the maximum principle, any two smooth admissible solutions of \eqref{eq:generalizedhessian} differ by a constant. At present, existence of a smooth admissible solution of \eqref{eq:generalizedhessian} is unknown in general. For admissible classes and normalized right-hand sides, existence is known for the ordinary $m$-Hessian equation \cite{dk,szekelyhidi} and in the linear case $m=1$.

Chen, Nie and Xu established a numerical criterion for ordinary $m$-positivity on projective manifolds with respect to uniform lower bounds of normalized intersection numbers \cite[Corollary 1.5]{cnx}. On the other hand Fu, Li, Zhang and Zhang established criteria using connected components and paths in $H^{1,1}(X,\R)$ on compact K\"ahler manifolds $X$ \cite{flzz}. Their theorem extends the connected-component description of the K\"ahler cone of Demailly and P\u aun \cite{demailly-paun}. It then follows that the cone of $m$-positive classes is determined by the class of the reference K\"ahler form as in Corollary~\ref{num:reference-class}. Here we discuss conjectures concerning sequential $m$-positivity and existence of solutions to the generalized $m$-Hessian equation. 

Assume $X$ is a compact K\"ahler manifold of dimension $d$. Choose $1\le m\le d$ and let $\eta_1,\dots,\eta_{d-m}$ be smooth closed real $(1,1)$-forms on $X$ which form a good prefix at every point. Let $\Theta$ be their wedge product. Let $\mathcal A_\Theta(X)$ be the cone of real $(1,1)$-classes represented by closed smooth forms $\theta$ such that $(\eta_1,\dots,\eta_{d-m},\theta)$ is pointwise good.

\begin{conjecture}[Numerical characterization of sequential $m$-positivity]\label{intro:numerical}
The cone $\mathcal A_\Theta(X)$ is the connected component containing the K\"ahler cone of the set of classes $\alpha\in H^{1,1}(X,\R)$ with
\[
 \int_Y[\Theta]\cdot\alpha^{p-d+m}>0,
 \qquad d-m+1\le p=\dim Y\le d,
\]
for every irreducible analytic subvariety $Y\subseteq X$ of the indicated dimensions.
\end{conjecture}

We will prove in Section~\ref{num:section} that every admissible class is contained in the component mentioned above. We restate Conjecture~\ref{intro:numerical} as Conjecture~\ref{num:component} and formulate a uniform conjecture as done by Chen, Nie and Xu. Another problem of interest is the solvability of the generalized $m$-Hessian equation.

\begin{conjecture}[Solvability of the generalized $m$-Hessian equation]\label{intro:solvability}
Let $\theta$ be a closed smooth real $(1,1)$-form with $[\theta]\in\mathcal A_\Theta(X)$. For $F\in C^\infty(X,\R)$ with
\[
 \int_X e^F\omega^d=\int_X\Theta\wedge\theta^m,
\]
there is a unique smooth real function $u$ (up to addition of a constant) such that
\[
 \Theta\wedge\theta_u^m=e^F\omega^d,\qquad
 \theta_u=\theta+\sqrt{-1}\partial\bar\partial u,\qquad
 \theta_u(x)\in\Gamma_{\Theta_x}\quad(x\in X).
\]
\end{conjecture}

We formulate Conjecture~\ref{intro:solvability} as Conjecture~\ref{num:existence}. Conjectures~\ref{intro:numerical} and~\ref{intro:solvability} together would give a cohomological criterion for the existence of a solution of the equation.

\subsection{The K\"ahler package and reverse inequalities}
The mixed Khovanskii--Teissier inequalities, log-concavity, and Lorentzian nature of intersection polynomials follow from the bidegree $(1,1)$ Hodge--Riemann relations. Hu and Xiao established ‘reverse Khovanskii--Teissier inequalities’ in \cite{hu-xiao-rkt,hu-xiao-shadow} using these results. We will apply their argument to the cone defined by the product $\Theta$ of a good prefix of length $d-m$ for some $m\ge 2$. For $1\le k\le m-1$ and $a,b_1,\ldots,b_k,c_1,\ldots,c_{m-k}\in\Gamma_{\Theta}(\omega)$, let $B=b_1\wedge\cdots\wedge b_k$ and $C=c_1\wedge\cdots\wedge c_{m-k}$. Then
\[
 \bigl(*\bigl(\Theta\wedge a^m\bigr)\bigr)\bigl(*\bigl(\Theta\wedge B\wedge C\bigr)\bigr)
 \le\binom{m}{k}\bigl(*\bigl(\Theta\wedge a^{m-k}\wedge B\bigr)\bigr)\bigl(*\bigl(\Theta\wedge a^k\wedge C\bigr)\bigr).
\]
This pointwise version is achieved by applying Kummer, Plaumann, and Vinzant's hyperbolic Rayleigh inequality \cite[Theorem 3.1]{kpv}. The global version of the inequality follows with the constant $2^{k(m-k)}$ by applying the Lorentzian argument of Hu and Xiao to the global sequential Hodge index theorem, as in \cite[Theorem 2.5 and Corollary 3.10]{hu-xiao-rkt} (Remark~\ref{rkt:global}). Assuming Conjecture~\ref{intro:solvability}, the sharp inequality holds with the constant $\binom{m}{k}$.

\begin{proposition}[Sharp reverse Khovanskii--Teissier inequality under solvability]\label{intro:rkt-sharp}
Suppose $\Theta$ is the product of a globally good prefix of length $d-m$ on a compact $d$-dimensional K\"ahler manifold $X$ with $2\le m\le d$. Suppose Conjecture~\ref{intro:solvability} holds with respect to $\Theta$. Choose $1\le k\le m-1$ and classes
\[
 \alpha,\beta_1,\ldots,\beta_k,\gamma_1,\ldots,\gamma_{m-k}\in\mathcal A_\Theta(X).
\]
Put $\beta=\beta_1\cdot\cdots\cdot\beta_k$ and $\gamma=\gamma_1\cdot\cdots\cdot\gamma_{m-k}$. Then
\[
\begin{aligned}
 &\left(\int_X[\Theta]\cdot\alpha^m\right)
  \left(\int_X[\Theta]\cdot\beta\cdot\gamma\right)\\
 &\qquad\le\binom{m}{k}
  \left(\int_X[\Theta]\cdot\alpha^{m-k}\cdot\beta\right)
  \left(\int_X[\Theta]\cdot\alpha^k\cdot\gamma\right).
\end{aligned}
\]
\end{proposition}

We employ the normalization technique of Lehmann and Xiao \cite[Theorem~5.7]{lehmann-xiao} as discussed in Proposition~\ref{rkt:conditional}. For ordinary $m$-positive classes, the sharp inequality holds unconditionally (Corollary~\ref{rkt:ordinary}). Huh, Micha\l ek, Wang and Xiao \cite[Theorem~1.2 and Corollary~1.3]{huh-michalek-wang-xiao} include the reverse inequalities in a family of correlation inequalities for volume polynomials. See Remark~\ref{rkt:correlation}.

\subsection*{Organization}
In Section~\ref{sec:definitions} we define hyperbolicity cones, positivity, polarized Hodge--Lefschetz modules and descent. The statement of the cone comparison theorem for sequential positivity is given in Section~\ref{seq:section} along with some boundary cases. From these results we obtain the mixed Hodge--Riemann relations for ordinary $m$-positive forms, and prove results using hyperplanes and Moore determinants. Weak $(m,l)$-positivity, power maps, and weak positivity with respect to a given prefix are considered in Section~\ref{sec:weak-power}. Theorems about sequential and weak positivity are proved in Sections~\ref{seq:proof-section} and~\ref{sec:lorentz} respectively. In Section~\ref{glob:section} we present and prove the complete intersection criterion on cohomology level. Some sequential numerical criteria as well as the generalized Hessian equation are presented in Section~\ref{num:section}. Some inequalities on intersection numbers of forms are obtained, including a sharp reverse inequality in the presence of solutions of the generalized Hessian equation. In Appendix~\ref{app:volume} we explain the difference between hyperbolicity of volume polynomials on cohomology and higher Lefschetz properties.

Weak and sequential $m$-positivity were defined by the author in April and May 2026, respectively. In October 2026, the author was informed that Jian Xiao had also formulated the Hodge--Riemann question for good sequences of the present article (special case of Theorem~\ref{main:blocks}) and Conjecture~\ref{intro:solvability}.

\subsection*{Acknowledgments}
The author is grateful to Jian Xiao for valuable remarks. He would also like to thank Gao Chen, Hao Fang, Jiajun Hu, Biao Ma, Ruiheng Zhang and Weizhe Zheng for several discussions.  

\subsection*{AI disclosure}
AI tools were used for some of the proof verification and literature searches. AI was also used for part of the numerical calculations and while experimenting to prove Theorem~\ref{main:strong} and Theorem~\ref{main:blocks}. The author independently checked all the proofs discovered by AI and takes full responsibility for the correctness of the paper.

\section{Preliminaries on \texorpdfstring{$m$}{m}-positivity}\label{sec:definitions}
Let $V$ be a $d$-dimensional complex vector space with a real $(1,1)$-form $\omega > 0$. Let $\eta$ be a real $(1,1)$-form on $V$. We choose an $\omega$-unitary coframe such that $\eta$ is diagonal and let
\[
 e_i=\sqrt{-1}\,dz_i\wedge d\bar z_i,\qquad
 \omega=\sum_i e_i,\qquad \eta=\sum_i a_i e_i.
\]
We set $\Lambda^{p,q}V^*=0$ if $p$ or $q$ is not in $\{0,\dots,d\}$. Let $\gamma$ be a form of the top degree. Define $*\gamma$ by $\gamma=(*\gamma)\vol_\omega$ with $\vol_\omega = \omega^d/d!$. Let $\wedge$ be the exterior product on forms and $\cdot$ be the cup product of cohomology classes. Let $e_I = \bigwedge_{i\in I}e_i$, $a_I = (a_i)_{i\in I}$ and $a_{\widehat I} = a_{I^c}$ for $I \subseteq [d] = \{1,\dots,d\}$. Let $\sigma_j$ denote the elementary symmetric polynomial of degree $j$ on $n$ variables, with $\sigma_0 = 1$ and $\sigma_j = 0$ for $j < 0$ or $j > n$.

\subsection{Hyperbolicity, ordinary positivity, and weak positivity}\label{pre:hyperbolicity}
Let $E$ be a finite dimensional real vector space. A homogeneous polynomial $P\colon E\to\mathbb R$ of degree $n \ge 1$ is said to be hyperbolic with respect to $e \in E$ if $P(e) > 0$ and, for all $x \in E$, the polynomial $t \mapsto P(x + t e)$ has only real roots. The negatives of these roots are the hyperbolic eigenvalues of $x$ with respect to $e$, denoted by $\lambda_\nu^e(x)$, and so
\[
 P(x+te)=P(e)\prod_{\nu=1}^n\bigl(t+\lambda_\nu^e(x)\bigr).
\]
The hyperbolicity cone of $P$ with respect to $e$, denoted by $\Gamma(P,e)$, is the connected component of $\{P>0\}$ containing $e$. Equivalently, it is the set of all $x$ such that $\lambda_\nu^e(x) > 0$ for all $\nu$. Following the work of G\aa rding on hyperbolic polynomials \cite{garding,harvey-lawson}, we note three properties.
\begin{enumerate}[label=(\alph*)]
\item $\Gamma(P,e)$ is open and convex, and each $a \in \Gamma(P,e)$ is a hyperbolicity direction of $P$ with $\Gamma(P,a) = \Gamma(P,e)$.
\item The function $P^{1/n}$ is concave on $\Gamma(P,e)$. With the symmetric polarization $\mathsf P$ of $P$ satisfying $\mathsf P(a,\dots,a) = P(a)$, we have
\[
 \mathsf P(a_1,\ldots,a_n)
 \ge\prod_{j=1}^n P(a_j)^{1/n}>0
 \qquad(a_1,\ldots,a_n\in\Gamma(P,e)).
\]
\item For $n\ge 2$ and $a\in\Gamma(P,e)$, the polynomial $D_aP$ is hyperbolic with respect to $e$ and
\[
 \Gamma(P,e)\subseteq\Gamma(D_aP,e).
\]
For $a\in\overline{\Gamma(P,e)}$, if $D_aP\not\equiv 0$ then $D_aP$ is hyperbolic with respect to $e$ and positive on $\Gamma(P,e)$.
\end{enumerate}
Thus, it is possible to differentiate iteratively with respect to directions in $\Gamma(P,e)$.  To establish sequential positivity we can choose the next direction in the cone of the last derivative.

A hyperbolic polynomial is \emph{complete} if $\overline{\Gamma(P,e)}\cap(-\overline{\Gamma(P,e)})=\{0\}$.
Equivalently, for one, and hence every, $a\in\Gamma(P,e)$, the factorization
\[
 P(a+tz)=P(a)\prod_{\nu=1}^{n}(1+t\lambda_\nu),\qquad\lambda_\nu\in\R,
\]
has $\lambda_1=\cdots=\lambda_n=0$ only for $z=0$.

\begin{lemma}\label{lem:complete}
Let $P$ be complete, hyperbolic of degree $\ge 3$. Then for $a\in\Gamma(P,e)$, the derivative $D_aP$ is also complete and hyperbolic.
\end{lemma}
\begin{proof}
Let $r=\deg P$ and $q=D_aP$. Assume that the hyperbolic eigenvalues of $z$ with respect to $q$ and $a$ are zero.  Then $q(a+tz)=q(a)$, so $D_zq(a)=D_z^2q(a)=0.$  By Euler’s theorem,
\[
 D_zq(a)=(r-1)D_zP(a),\qquad
 D_z^2q(a)=(r-2)D_z^2P(a).
\]
Thus the coefficients of $t$ and $t^2$ in $P(a+tz)/P(a)=\prod_\nu(1+t\lambda_\nu)$ must be 0. Thus, $\sum_\nu\lambda_\nu=0$ and $\sum_{\mu<\nu}\lambda_\mu\lambda_\nu=0$. Therefore, $\sum_\nu\lambda_\nu^2=0$. This implies each of the $\lambda_\nu = 0$, as the $\lambda_\nu$ are real numbers. From the completeness of $P$ it then follows that $z = 0$.
\end{proof}

\begin{definition}
For $0\le r\le n$, let $\Gamma_r^n=\{a\in\R^n:\sigma_j(a)>0\text{ for }1\le j\le r\}$. In particular, $\Gamma_0^n=\R^n$.
An element $\eta$ is said to be ordinarily $(\omega,m)$-positive if $a\in\Gamma_m^d$, i.e., $\omega^{d-j}\wedge\eta^j>0$ for $1\le j\le m$.
For $0\le l\le m\le d$, define
\begin{equation}\label{codim:cone}
 C_{m,l}^d=\{a\in\R^d:a_{\widehat I}\in\Gamma_{m-l}^{d-l}
       \text{ for every }I\subseteq[d],\ |I|=l\}.
\end{equation}
We say $\eta$ is weakly $(m,l)$-positive with respect to $\omega$ if $a\in C_{m,l}^d$. By the eigenvalue interlacing theorem this is equivalent to ordinary $(\omega|_W, m-l)$-positivity of $\eta|_W$ on any complex subspace $W$ of codimension $l$. See Lemma~\ref{codim:cones}.
For $2\le m\le d$, the weak $m$-positive cone is
\begin{equation}\label{eq:weak}
 \wtG_m^d=\{a\in\Gamma_{m-1}^d:
              \sigma_{m-1}(a_{\widehat i})>0\text{ for every }i\}.
\end{equation}
By Proposition~\ref{prop:deletion}, $\wtG_m^d=C_{m,1}^d$. We set $\wtG_1^d=C_{1,1}^d=\R^d$.
\end{definition}
Let $\Gamma_m(\omega)$ and $\widetilde\Gamma_m(\omega)$ denote the respective cones of forms. For $l=0$ we have $C_{m,0}^d=\Gamma_m^d$, and for $l=m$ we have $C_{m,m}^d=\R^d$. The G\aa rding component of the set $\{\sigma_r>0\}$ for $r\ge 1$ is defined by the inequality constraints in the definition of $\Gamma_r^n$ for $1\le j\le r$.

For $m\ge 1$, the cones in \eqref{codim:cone} are components of the domain filtration of Fang and Ma \cite[Definition 5.4, Proposition 5.5, Definition 8.4]{fang-ma-garding}. They define
\[
 \mathcal C_f^{(l)}=\bigcap_{|\alpha|=l}\mathcal C_{\partial^\alpha f},
\]
where $\mathcal C_{\partial^\alpha f}=\R^d$ for nonnegative constant derivatives. Since for $|I|=l$, $\partial_I\sigma_m(a)=\sigma_{m-l}(a_{\widehat I})$, and derivatives involving a repeated coordinate are zero, and since $\mathcal C_{\sigma_r}=\Gamma_r$ for the other coordinates, we get
\begin{equation}\label{codim:fang-ma}
 C_{m,l}^d=\mathcal C_{\sigma_m}^{(l)}.
\end{equation}
In particular,
$C_{m,0}^d\subseteq C_{m,1}^d\subseteq\cdots\subseteq C_{m,m}^d$.

To study real $(k,k)$-forms we use the positivity cones of Harvey and Knapp \cite{harvey-knapp}. A form $\Phi$ is weakly positive if its restriction to any complex $k$-plane is nonnegative with respect to the complex orientation.  A form $\Phi$ is strongly positive if it can be written as a sum of products of $k$ semipositive $(1,1)$-forms of rank 1. A form $\Phi$ is Hermitian positive if its associated Hermitian matrix on $\Lambda^k V$ is positive semidefinite. In all these cases we can define strictly positive forms $\Phi$ if $\Phi-\varepsilon\omega^k$ is in the corresponding cone for some positive number $\varepsilon$. If $k=1, d-1,$ or $d,$ then these three cones are the same.

For $m\ge 2$ we have strict inequalities as in \eqref{eq:weak}. If $m=l=1$ then there is nothing to prove and the Lefschetz coefficient for $p+q=1$ is $\omega^{d-1}$.

\begin{lemma}[Hyperplane criterion for strong positivity]\label{seq:strong-hyperplane}
Let $1\le k\le d$. Let $\Psi$ be a real $(k-1, k-1)$-form on $V$. Let $(\phi_i)_{i=1}^d$ be an $\omega$-orthonormal coframe on $V$. Suppose $\Psi|_{W_i}$ is strongly positive on $W_i=\ker\phi_i$ for each $i$. Then $\omega\wedge\Psi$ is strongly positive. If each $\Psi|_{W_i}$ is strictly strongly positive, then $\omega\wedge\Psi$ is strictly strongly positive.
\end{lemma}
\begin{proof}
Set $e_i=\sqrt{-1}\phi_i\wedge\bar\phi_i$. Extend any of the covectors of a strongly positive decomposition of $\Psi|_{W_i}$ on $W_i$ to $V$. If we take the wedge product of $e_i$ with any of these extended covectors, the terms in the wedge product involving $\phi_i$ and $\bar\phi_i$ vanish. Thus, any extension defines a strongly positive representation of $e_i\wedge\Psi$. Thus, $\omega\wedge\Psi$ is strongly positive. If each $\Psi|_{W_i}$ is strictly strongly positive, choose $\delta>0$ small enough such that $(\Psi-\delta\omega^{k-1})|_{W_i}$ is strongly positive for each $i$. Now, use the first claim with $\Psi-\delta\omega^{k-1}$ in place of $\Psi$.
\end{proof}

\begin{lemma}[Simultaneous equal diagonals]\label{seq:two-equal-diagonals}\label{lem:equaldiagonal}
Two Hermitian forms on a $d$-dimensional complex Hermitian vector space admit a common orthonormal basis in which their diagonal entries are constant, equal to their respective traces divided by $d$.
\end{lemma}
\begin{proof}
Let $A,B$ be the Hermitian matrices in question and let $C = A+iB$. Then the Toeplitz--Hausdorff theorem states that the numerical range of $C$ is convex \cite{toeplitz,hausdorff}. For any orthonormal basis, the diagonal elements of $C$ lie in its numerical range, and hence so does its average $\tr(C)/d$. This implies there exists a unit vector $v$ such that $v^*Cv = \tr(C)/d$. Then $v^*Av = \tr(A)/d$ and $v^*Bv = \tr(B)/d$. Let $d>1$. Then the restrictions of $A$ and $B$ to $v^\perp$ also have average eigenvalue $\tr(A)/d$ and $\tr(B)/d$, respectively. We may apply induction to determine the desired basis for the whole space.
\end{proof}

\begin{proposition}[Deletion characterization]\label{prop:deletion}
Let the eigenvalues of $\eta$ with respect to $\omega$ be $a_1,\dots,a_d$. Let $1 \le m \le d$. The following are equivalent.
\begin{enumerate}[label=(\alph*)]
\item $a\in\wtG_m^d$;
\item $a_{\widehat i}\in\Gamma_{m-1}^{d-1}$ for every $i$;
\item $\eta|_H$ is $(\omega|_H,m-1)$-positive on every complex hyperplane $H$.
\end{enumerate}
Equivalently, $\eta$ is ordinarily $(m-1)$-positive and
$\omega^{d-m}\wedge\eta^{m-1}$ is strictly positive as a $(d-1,d-1)$-form. Moreover,
\begin{equation}\label{eq:inclusions}
 \Gamma_m^d\subseteq\wtG_m^d\subseteq\Gamma_{m-1}^d,
\end{equation}
and $\wtG_m^d$ is open and convex.
\end{proposition}
\begin{proof}
By ordinary deletion, if $a\in \Gamma_r^n$ and $1\le r\le n$, then $a_{\widehat{i}} \in \Gamma_{r-1}^{n-1}$ for all $i$. This is established by differentiating $\sigma_r$ with respect to the $i$th coordinate direction, as in \cite{garding,harvey-lawson}. Under condition (a), deletion results in inequalities of degree 1 up to $m-2$ in (b), and the inequality of degree $m-1$ in \eqref{eq:weak} holds when $m \ge 2$. Under condition (b), we have
\[
 (d-j)\sigma_j(a)=\sum_i\sigma_j(a_{\widehat i})>0,
 \qquad1\le j\le m-1,
\]
Thus $a\in \Gamma_{m-1}^d$ and condition (a) is met. Let $v=(v_1,\dots,v_d)$ be a unit normal in the eigenbasis of $\eta$. By the Cauchy--Binet formula
\begin{equation}\label{eq:compression}
 \sigma_j(\eta|_{v^\perp})=\sum_i|v_i|^2\sigma_j(a_{\widehat i}).
\end{equation}
Thus (b) implies (c). Restricting to coordinate hyperplanes gives the converse. An equivalent formulation in terms of $\omega^{d-m}\wedge\eta^{m-1}$ is given by the expansion
\[
 \omega^{d-m}\wedge\eta^{m-1}
 =(d-m)!(m-1)!\sum_i\sigma_{m-1}(a_{\widehat i})e_{[d]\setminus\{i\}}.
\]
The inclusions in \eqref{eq:inclusions} hold true by ordinary deletion and \eqref{eq:weak}. By (b), $\wtG_m^d$ is the intersection of the preimages of the open convex cone $\Gamma_{m-1}^{d-1}$ under the linear maps $a\mapsto a_{\widehat{i}}$. Since restriction to each fixed hyperplane is linear and the ordinary cone on that hyperplane is convex, (c) also implies the cone of forms is convex.
\end{proof}

\begin{corollary}\label{cor:repeated}
If $a\in C_{m,l}^d$ and $l\le |S|\le m$, then
$a_{S^c}\in\Gamma_{m-|S|}^{d-|S|}$. Therefore, for $l\le i\le m$,
\begin{equation}\label{eq:strong}
 \omega^{d-m}\wedge\eta^{m-i}
 =(d-m)!(m-i)!\sum_{|K|=d-i}\sigma_{m-i}(a_K)e_K
\end{equation}
has positive coefficients on the basis $(e_K)_{|K|=d-i}$ of diagonal forms. Thus, $\omega^{d-m}\wedge\eta^{m-i}$ is strictly strongly positive with positive definite Hermitian matrix of coefficients. This holds for $1\le i\le m$ if $a\in\wtG_m^d$, and for $i=0$ if $a\in\Gamma_m^d$.
\end{corollary}
\begin{proof}
Delete $l$ coordinates indexed by $S$ using \eqref{codim:cone}. Now delete the remaining $|S|-l$ coordinates of $S$ using ordinary deletion. From the eigenbasis expansion of the wedge product of $\eta$ we get \eqref{eq:strong}. By the first assertion with $S=K^c$, we have $\sigma_{m-i}(a_K) > 0$. Thus, after removing a small positive multiple of $\omega^{d-i}$, the coefficients will remain positive.
\end{proof}
We prove the higher Lefschetz isomorphisms for these forms in Section~\ref{sec:higher}.

\begin{proposition}[Number of positive eigenvalues]\label{ord:positive-index}
Let $0\le l < m\le d$. If the eigenvalues $\lambda_1\le\dots\le\lambda_d$ of $\eta$ lie in $C_{m,l}^d$, then
\[
 \lambda_1+\cdots+\lambda_{d-m+1}>0.
\]
In particular, $\eta$ has at least $m$ positive eigenvalues.
\end{proposition}
\begin{proof}
Apply Corollary~\ref{cor:repeated} with $S$ indexing the $m-1$ largest eigenvalues. Since $(\lambda_1,\ldots,\lambda_{d-m+1})\in\Gamma_1^{d-m+1}$,
\[
 0<\lambda_1+\cdots+\lambda_{d-m+1}
   \le(d-m+1)\lambda_{d-m+1},
\]
so the last $m$ eigenvalues are positive.
\end{proof}

\begin{remark}[Comparison of positivity notions]\label{pre:curvature-comparison}
A scalar curvature form is said to be $q$-positive in the sense of Andreotti and Grauert if it has at least $d-q$ positive eigenvalues \cite{andreotti-grauert}. Sommese, in his theory, as stated in \cite[Definition 2.3]{yang-ks}, also requires the form to be semipositive. Uniform $q$-positivity of $\eta$ with respect to $\omega$ means that the sum of any $q+1$ eigenvalues of $\eta$ is positive \cite[Definition 2.1]{yang-ag}. In the terminology of Dinew and Popovici, $s$-positivity means strict strong positivity of $\eta\wedge\omega^{s-1}$, which is the same as the sum of any $s$ eigenvalues of the form being positive \cite[Definition~1.1 and Lemma~1.2]{dinew-popovici}. Thus, by Proposition~\ref{ord:positive-index}, for $l<m$, any weakly $(m,l)$-positive form is uniformly $(d-m)$-positive with respect to $\omega$, or equivalently, $(d-m+1)$-positive in the sense of Dinew and Popovici.

For vector bundles, Demailly defines $s$-positivity of the curvature Hermitian form by positivity of the form on all nonzero tensors of rank at most $s$ \cite[Definitions~2.1--2.2]{demailly-positivity}. Testing on rank 1 tensors yields the fact that the scalar curvature form in each nonzero fiber direction is positive definite. For a line bundle, this means the line bundle has positive curvature.
\end{remark}

\subsection{Positivity relative to a prefix and sequential positivity}\label{pre:relative}\label{pre:sequential}

Fix $0\le k<d$ and a tuple $\boldsymbol\omega=(\omega_1,\ldots,\omega_k)$ of real $(1,1)$-forms. Put
\[
 \Theta=\omega_1\wedge\cdots\wedge\omega_k,\qquad n=d-k.
\]
The empty prefix has $\Theta = 1$. Consider prefixes for which the homogeneous polynomial
\begin{equation}\label{pre:prefix-polynomial}
 P_\Theta(A)=*\bigl(\Theta\wedge A^n\bigr)
 \qquad\bigl(A\in\Lambda^{1,1}_{\R}V^*\bigr)
\end{equation}
is hyperbolic with respect to $\omega$. Then $P_{\Theta}(\omega) > 0$. For such prefixes, the cone
\[
 \Gamma_\Theta(\omega)=\Gamma(P_\Theta,\omega)
\]
is the hyperbolicity cone with respect to the prefix $\boldsymbol{\omega}$ which is completely determined by $\Theta$ and the component containing $\omega$. We denote it as $\Gamma_{\Theta}$ for brevity.

For $1\le m\le n$, a real $(1,1)$-form $\eta$ is \emph{$(\Theta,\omega,m)$-positive} if
\begin{equation}\label{pre:relative-m-positive}
 *\bigl(\Theta\wedge\eta^j\wedge\omega^{n-j}\bigr)>0
 \qquad(1\le j\le m).
\end{equation}
Let $\Gamma_{\Theta,m}(\omega)$ be the cone defined by \eqref{pre:relative-m-positive}. We can write this cone as a hyperbolicity cone by considering the polynomial
\[
 P_{\Theta,m}(A)=*\bigl(\Theta\wedge\omega^{n-m}\wedge A^m\bigr)
 =\frac{m!}{n!}D_\omega^{n-m}P_\Theta(A).
\]
This polynomial is hyperbolic with respect to $\omega$ by repeated differentiation, and
\[
 P_{\Theta,m}(\eta+t\omega)
 =\sum_{j=0}^{m}\binom mj
   *\bigl(\Theta\wedge\eta^j\wedge\omega^{n-j}\bigr)t^{m-j}.
\]
Since the first coefficient of this polynomial is positive and all of its roots are real, the requirements in \eqref{pre:relative-m-positive} say that all of these roots are negative. Therefore,
\begin{equation}\label{pre:relative-cones}
 \Gamma_{\Theta,m}(\omega)=\Gamma(P_{\Theta,m},\omega),\qquad
 \Gamma_\Theta(\omega)=\Gamma_{\Theta,n}(\omega),
\end{equation}
and
\[
 \Gamma_{\Theta,m}(\omega)\subseteq\Gamma_{\Theta,m-1}(\omega)
 \qquad(2\le m\le n).
\]
In particular, these cones are open and convex, and G\aa rding's inequality holds for them.

By \eqref{pre:relative-m-positive} and \eqref{pre:prefix-polynomial}, ordinary $m$-positivity corresponds to $\Gamma_m(\omega)=\Gamma_{1,m}(\omega)=\Gamma_{\omega^{d-m}}(\omega)$.

Let $\Theta$ be fixed. Every choice of reference direction in $\Gamma_{\Theta}$ defines the same cone. For $m < n$, the cone $\Gamma_{\Theta,m}(\omega)$ can be dependent on $\omega^{n-m}$.

We use relative positivity as in Zhang’s Hodge index theorem \cite{zhang}.

\begin{definition}\label{seq:definition}
Let $(\eta_1,\dots,\eta_k)$ be real $(1,1)$-forms, where $0\le k\le d$. Let $\Theta_0=1$ and $\Theta_j = \eta_1\wedge\dots\wedge\eta_j$. We call $(\eta_1,\dots,\eta_k)$ a good prefix if
\begin{equation}\label{seq:inequalities}
 *\bigl(\Theta_{j-1}\wedge\eta_j^a\wedge\omega^{d-j+1-a}\bigr)>0,
 \qquad 1\le j\le k,\quad 1\le a\le d-j+1.
\end{equation}
We then say $(\eta_1,\dots,\eta_k)$ is sequentially positive. If $k=d$, we call it a good sequence.
\end{definition}
The order of factors is relevant. At stage $j$, we assume that $\eta_j$ is positive with respect to $(\Theta_{j-1},\omega,d-j+1)$. This means that $\eta_1$ is positive definite. Other factors can be indefinite.

Using $\omega$, identify real $(1,1)$-forms with Hermitian matrices. For $0\le j\le\min\{k,d-1\}$, set
\begin{equation}\label{seq:descendant}
 F_j(A)=D_{\eta_1}\cdots D_{\eta_j}\det(A)
       =\frac{P_{\Theta_j}(A)}{(d-j)!}=\frac{*\bigl(\Theta_j\wedge A^{d-j}\bigr)}{(d-j)!}.
\end{equation}

\begin{proposition}\label{seq:hyperbolicity}
Let $(\eta_1,\ldots,\eta_k)$ be a good prefix. For $0\le j\le\min\{k,d-1\}$, $F_j$ is hyperbolic with respect to $\omega$. The prefix $(\eta_1,\ldots,\eta_j,A)$ is good if and only if $A\in\Gamma_{\Theta_j}=\Gamma(F_j,\omega)$. These are open, convex cones which contain all positive definite forms and satisfy
\[
 \Gamma(F_j,\omega)\subseteq\Gamma(F_{j+1},\omega)
 \qquad(0\le j<\min\{k,d-1\}).
\]
We can use any positive definite form for $\omega$. Good prefixes stay good upon positive scaling of their individual terms and repetition of their last term, as long as the total length is at most $d$.
\end{proposition}
\begin{proof}
We proceed by induction on $j$. For $F_0=\det$, the hyperbolicity cone is the cone of positive definite forms. Assume $F_{j-1}$ is hyperbolic. By \eqref{pre:relative-cones} and \eqref{seq:inequalities}, $(\eta_1,\ldots,\eta_{j-1},A)$ is good if and only if $A\in\Gamma(F_{j-1},\omega)$. Hence, $\eta_j$ is in this cone. By G\aa rding's theorem for $F_j=D_{\eta_j}F_{j-1}$, $F_j$ is hyperbolic and $\Gamma(F_{j-1},\omega)\subseteq\Gamma(F_j,\omega)$ \cite{garding,harvey-lawson}. By \eqref{pre:relative-cones}, possible factors in the next step are elements of $\Gamma(F_j,\omega)$. Thus, all cones contain the positive definite cone and any positive definite reference form selects the same component. Rescaling prefix polynomials by a positive constant maintains their hyperbolicity. Also repetition of terms is allowed since $\eta_j\in\Gamma(F_{j-1},\omega)\subseteq\Gamma(F_j,\omega)$ for $j<d$.
\end{proof}
Let $\eta_1,\ldots,\eta_k$ be fixed and $\omega$ vary. For the ordinary cone $\Gamma_m(\omega)=\Gamma_{\omega^{d-m}}(\omega)$, various values of $\omega$ yield various prefixes if $m<d$.

\begin{corollary}\label{seq:complete}
For a good prefix of length $k$, $F_j$ is complete hyperbolic for $0\le j\le\min\{k,d-2\}$.
\end{corollary}
\begin{proof}
The determinant is complete hyperbolic since a Hermitian matrix with all eigenvalues zero is the zero matrix. Applying Lemma~\ref{lem:complete} for $\eta_1,\ldots,\eta_j$ completes the proof.
\end{proof}

\begin{definition}\label{seq:almost-definition}\label{seq:sp-relative}
Let $H$ be a nonzero semipositive form of type $(1,1)$. Let $\Theta_0=H$ and $\Theta_j=H\wedge B_1\wedge\dots\wedge B_j$. We say $(H;B_1,\dots,B_s)$ is relatively good if $B_j\in \Gamma_{\Theta_{j-1}}$ for $1\le j\le s$. It is almost good if $B_j\in\overline{\Gamma_{\Theta_{j-1}}}$ and $\Theta_j\ne 0$ for $1\le j\le s$. A nonempty prefix is good if and only if it is relatively good and $H>0$. We allow the empty prefix.
\end{definition}
The directional derivative property shows each $P_{\Theta_j}$ is non-zero and hyperbolic. Thus $\Gamma_{\Theta_j}$ contains the positive cone. When $P_{\Theta}$ is constant, we define the cone to be $\Lambda^{1,1}_{\R}V^*$. For weak positivity relative to a prefix, see Subsection~\ref{seq:weak-section}.

\subsection{Hodge--Lefschetz modules}\label{pre:modules}
Let $p,q \geq 0$ be such that $p+q \leq d$. For a real $(d-p-q,d-p-q)$-form $\Omega$, a real $(1,1)$-form $\lambda$, and $\alpha,\beta \in \Lambda^{p,q}V^*$ we set
\begin{align}
 Q_\Omega(\alpha,\beta)
 &=(\sqrt{-1})^{p-q}(-1)^{(p+q)(p+q-1)/2}
                         *\bigl(\Omega\wedge\alpha\wedge\bar\beta\bigr),\label{eq:Q}\\
 P^{p,q}(\Omega,\lambda)
 &=\ker(\Omega\wedge\lambda\wedge-\colon
       \Lambda^{p,q}V^*\to\Lambda^{d-q+1,d-p+1}V^*).\label{eq:primitive}
\end{align}
We say $\Omega$ has the hard Lefschetz property in bidegree $(p,q)$ if the map $\Omega\wedge-\colon \Lambda^{p,q}V^* \rightarrow \Lambda^{d-q,d-p}V^*$ is an isomorphism. We call $(\Omega, \lambda)$ a Hodge--Riemann pair in bidegree $(p,q)$ if $Q_{\Omega}$ is positive definite on $P^{p,q}(\Omega,\lambda)$. HR implies HL. Indeed, if $\Omega \wedge \alpha = 0$, then $\alpha \in P^{p,q}(\Omega,\lambda)$ and $Q_{\Omega}(\alpha,\alpha) = 0$. Hence $\alpha = 0$. The dimensions of $\Lambda^{p,q}V^*$ and $\Lambda^{d-q,d-p}V^*$ are the same.

Li and Zheng \cite[Definition~1.1]{li-zheng} assume both HL and primitive positivity for $(\Omega \wedge \lambda^{2j}, \lambda)$ in bidegree $(p-j, q-j)$, with $0 \leq j \leq \min(p,q)$. Our definition is limited to $(p,q)$. We state lower-bidegree hypotheses explicitly in Theorem~\ref{glob:theorem}. Nondegeneracy of $Q_{\Omega \wedge \lambda^2}$ on $\Lambda^{p-1,q-1}V^*$ implies nondegeneracy of $Q_{\Omega}$ on $\lambda \wedge \Lambda^{p-1,q-1}V^*$, by (\ref{gp:eq:shift-sign}). The $Q_\Omega$-orthogonal complement is $P^{p,q}(\Omega, \lambda)$ as a result of exterior duality. Hence,
\begin{equation}\label{eq:LD}
 \Lambda^{p,q}V^*=P^{p,q}(\Omega,\lambda)
                         \oplus\lambda\wedge\Lambda^{p-1,q-1}V^*.
\end{equation}
We use identical definitions on the cohomology of a compact Kähler manifold by replacing $*$ with integration.

Let $\Theta$ be a nonzero real $(r,r)$-form. Let $n = d - r$. Its exterior quotient is defined by
\[
 M_\Theta=\Lambda^{*,*}V^*/\Ann(\Theta),\qquad
 \Ann(\Theta)=\{u:\Theta\wedge u=0\}.
\]
Conjugation and the Hodge bigrading descend to $M_\Theta$. Exterior duality induces the pairing $(u,v)\mapsto *(\Theta\wedge u\wedge v)$ which is a perfect pairing of $M_\Theta^{p,q}$ and $M_\Theta^{n-p,n-q}$. Since $\Ann(\Theta)$ is an ideal, each real $(1,1)$-form operates on $M_\Theta$.

A module $M$ of weight $n$ is a finite-dimensional real vector space with a Hodge bigrading, and with a perfect graded-symmetric real pairing $B$ of type $(n,n)$ on complementary degrees. Complex conjugation maps $M^{p,q}$ to $M^{q,p}$. On $M$ we consider a family of commuting real operators $a$ of bidegree $(1,1)$ satisfying $B(au,v)=B(u,av)$. Such an operator $a$ polarizes $M$ if for all $p+q\le n$ the map
\[
 a^{n-p-q}\colon M^{p,q}\longrightarrow M^{n-q,n-p}
\]
is an isomorphism, and if the pairing on $\ker(a^{n-p-q+1}|_{M^{p,q}})$ induced by $a^{n-p-q}$ with the sign given by \eqref{eq:Q} is positive definite. For $M_\Theta$ this pairing is $Q_{\Theta\wedge a^{n-p-q}}$.

According to \cite[Proposition~1.2.2]{kashiwara-kawai}, each connected component of the set of polarizing operators is an open convex cone in the real vector space spanned by the commuting operators. Also see \cite[Proposition~3.11]{li-zheng}. Let $\omega$ polarize $M_\Theta$, and let $\operatorname{Pol}_\Theta(\omega)$ be the inverse image in $\Lambda^{1,1}_\R V^*$ of the component containing the map of multiplication by $\omega$. For any $a\in\operatorname{Pol}_\Theta(\omega)$, $P_\Theta(a)>0$ because of positivity in bidegree $(0,0)$. If $P_\Theta$ is hyperbolic with respect to $\omega$, then connectedness implies
\begin{equation}\label{seq:cone-inclusion}
 \operatorname{Pol}_\Theta(\omega)\subseteq\Gamma_\Theta(\omega).
\end{equation}
For an almost-good prefix, Theorem~\ref{main:blocks} yields the reverse inclusion.

To shift to the standard grading of a Hodge--Lefschetz module, we introduce the following definitions: $M_\ell=M^{n-\ell}$ and $M_{\rm mod}^{u,v}=M^{n-v,n-u}$. Let $M=M_\Theta$, and define $\tau(u\wedge v)=*(\Theta\wedge u\wedge v)$ and
\[
 S(u,v)=(-1)^{k(k-1)/2}\tau(u\wedge v),
 \qquad u\in M^k,\quad v\in M^{2n-k}.
\]
Here the operators are of bidegree $(-1,-1)$ and degree $-2$.

\begin{lemma}[Mixed relations on a polarizing cone]\label{pre:mixed-module}
Let $\mathcal K$ be a polarizing cone of a Hodge--Lefschetz module $M$ of weight $n$. Fix $p+q\le n$, $a_1,\ldots,a_{n-p-q},\lambda\in\mathcal K$, and let $C=a_1\cdots a_{n-p-q}$. Then $C\colon M^{p,q}\to M^{n-q,n-p}$ is an isomorphism. The pairing on $P=\ker(C\lambda|_{M^{p,q}})$ is positive definite. Also, there exists an orthogonal splitting $M^{p,q} = P\oplus\lambda M^{p-1,q-1}$.
\end{lemma}
From the descent lemma of Cattani, Kaplan and Schmid \cite[Lemma 1.16]{cks}, the mixed hard Lefschetz theorem holds for $C$. The mixed primitive relations also hold, for $p+q \geq 2$ \cite[Theorems~3.2 and~4.2--4.4]{cattani}. For $p+q \leq 1$, all elements of $M^{p,q}$ are primitive. Move $a_i$ continuously in $\mathcal K$ to $\lambda$. By mixed hard Lefschetz, the associated pairing on $M^{p,q}$ is still nondegenerate and so has the same signature as the positive pairing on $\lambda^{n-p-q}$.

\begin{lemma}[Boundary descent]\label{gp:lem:descent}
Let $M$ be a polarized Hodge--Lefschetz module of weight $w$ with
polarizing cone $\mathcal K$, and let $T\in\overline{\mathcal K}$.
The module $M/\ker T^a$, with pairing induced by $T^a$, has weight $w-a$ and polarization $\mathcal K$ for $1 \le a \le w$. This is also the case for the zero module.
\end{lemma}
\begin{proof}
Let $L \in \mathcal K$. If $T=0$, $M/\ker T^a=0$. Let $\delta > 0$ such that $L-\delta T \in \mathcal K$. Then $L$ lies in the relative interior of the cone spanned by $T$ and $L-\delta T$. Since $T \in \overline{\mathcal K}$ and $L-\delta T \in \mathcal K$, any positive linear combination of $T$ and $L-\delta T$ is in $\mathcal K$. Applying the descent lemma of Cattani, Kaplan and Schmid \cite[Lemma 1.16]{cks} to the split polarized mixed Hodge structure associated to the given grading we obtain a polarization of the regraded image $TM$ by $L$ with the pairing $\widetilde S(Tu,Tv)=S(u,Tv)$. By iterating this and using the isomorphism $[u] \mapsto T^a u$ from $M/\ker T^a$ to $T^a M$ we obtain the result. For representations as quotients and as Hodge--Lefschetz modules, see \cite[Lemma 4.5]{ahl} and \cite[Theorem 3.2, Remark 3.3]{cattani}.
\end{proof}

\begin{lemma}[The central quotient]\label{gp:lem:central}
Under the hypotheses of Lemma~\ref{gp:lem:descent}, the quotient $\mathcal C_T(M) = \ker T / (\ker T \cap TM)$ has a polarization of weight $w$ induced by $\mathcal K$.
\end{lemma}
\begin{proof}
In the representation of $M$ as a sum of $T$-strings, $\ker T$ corresponds to the endpoints of strings. The ends of $T$-strings of length $>1$ belong to $TM$. Therefore, $\mathcal C_T(M)$ is the zero-index piece of the perverse decomposition of $M$ (see Section 4.3 of \cite{ahl}). The commuting operators descend to the quotient, and its pairing is induced by restriction (see Section 4.4 of \cite{ahl}). From \cite[Corollary 1.11, Theorem 3.13]{ahl}, $\mathcal K$ polarizes the quotient with the induced Hodge bigrading.
\end{proof}

We write $Q_\Xi^{p,q}$ for the pairing in bidegree $(p,q)$. From \eqref{eq:Q}, for $u,v \in \Lambda^{p-1,q-1}V^*$ we have
\begin{equation}\label{gp:eq:shift-sign}
 Q_\Xi^{p,q}(\lambda\wedge u,\lambda\wedge v)
 =-Q_{\Xi\wedge\lambda^2}^{p-1,q-1}(u,v).
\end{equation}
For a Hermitian form $Q$, we let $n_+(Q)$ and $n_-(Q)$ be the positive and negative indices of $Q$, respectively. Let $M$ be a polarized module. Let $h_M^{a,b}=\dim M^{a,b}$. The pairing is positive definite on the primitive subspace. In view of \eqref{gp:eq:shift-sign}, each polarizer multiplication switches the sign of the pairing in the Lefschetz decomposition. Consequently, the negative index in bidegree $(p,q)$ with $p+q\le n$ is
\begin{equation}\label{gp:eq:expected-index}
 N_M(p,q)=\sum_{j\ge1}(-1)^{j-1}h_M^{p-j,q-j}.
\end{equation}
Here $h_M^{a,b}=0$ if $a<0$ or $b<0$. For $p+q=n$ and a polarizer $a$, there is a maximal negative subspace in $aM^{p-1,q-1}$.

\begin{lemma}\label{gp:lem:hyperbolic-facts}
Let $p$ be a homogeneous hyperbolic polynomial with cone $\Gamma$ and $v$ be in the closure of $\Gamma$. Then $D_vp=0$ or $D_vp$ is hyperbolic and its cone contains $\Gamma$. If $v\in\Gamma$ and $\deg p>0$, then $D_vp>0$ on $\Gamma$.

An almost good prefix of length at most $d-1$ restricts to an almost good prefix on a dense open set of complex hyperplanes $W$. For such $W$, restriction maps each $\Gamma_{\Theta_j}$ into $\Gamma_{\Theta_j|_W}$. If $(H;B_1,\ldots,B_s)$ is relatively good of length at most $d-1$, then $(H|_W;B_1|_W,\ldots,B_s|_W)$ is relatively good for any complex hyperplane $W$ for which $H|_W\ne0$.
\end{lemma}
\begin{proof}
For any $x\in\Gamma$ and $v\in\overline\Gamma$, the eigenvalues of $v$ with respect to $x$ are non-negative.
\[
 D_vp(x)=p(x)\sum_\nu\lambda_\nu^x(v).
\]
Suppose $D_vp(x)=0$. Then all of these eigenvalues are 0, which means $v\in\overline\Gamma$ and $-v\in\overline\Gamma$. Using the same logic on any other $y\in\Gamma$ yields $D_vp(y)\ge0$ and $D_{-v}p(y)\ge0$. Hence, $D_vp\equiv 0$. On the other hand, if $D_vp(x)>0$ then $D_vp(y)>0$ for all $y\in\Gamma$. To show that $D_vp$ is hyperbolic, we use the fact that $D_{v+\varepsilon x}p(x)>0$ for small positive $\varepsilon$ and apply G\aa rding’s theorem \cite{garding}. The limit is real rooted with respect to $x$ because $D_vp(x)>0$. It follows that $\Gamma$ is inside the hyperbolicity cone of $D_vp$. In particular, if $v\in\Gamma$ and $\deg p>0$, the sum above is strictly positive.

For any non-zero covector $u$, let $W=\ker u$ and $\rho_u=\ii u\wedge\bar u$. Let $\Theta_j$ be any partial product of length $j\le d-1$. Then the polynomial $D_{\rho_u}P_{\Theta_j}(A)$ is a positive multiple of $P_{\Theta_j|_W}(A|_W)$. This is because multiplying by $\rho_u$ with wedge kills terms containing $u$ and $\bar u$ in the polynomial. The function $u\mapsto D_{\rho_u}P_{\Theta_j}(\omega)$ takes non-negative values because $\rho_u$ lies in the closed cone of $P_{\Theta_j}$. It is a non-zero quadratic form in $u$, since its sum over an $\omega$-unitary coframe equals $D_\omega P_{\Theta_j}(\omega)=(d-j)P_{\Theta_j}(\omega)>0$. Hence, on a dense open set of $u$, these values are positive for all $j$, and the cone inclusion of $D_{\rho_u}P_{\Theta_j}$ shows that restriction to $W$ takes $\Gamma_{\Theta_j}$ to $\Gamma_{\Theta_j|_W}$. For $u$ in this set, restriction also takes $\overline{\Gamma_{\Theta_j}}$ to $\overline{\Gamma_{\Theta_j|_W}}$, because of continuity.

Assume $(H;B_1,\ldots,B_s)$ is relatively good, and take a hyperplane $W=\ker u$ such that $H|_W\ne0$. Because $H|_W$ is semipositive, $D_{\rho_u}P_H(\omega)>0$. For each $j=1,\dots,s$, let $P_j=P_{H\wedge B_1\wedge\cdots\wedge B_j}$ with $P_0=P_H$. If $D_{\rho_u}P_{j-1}\ne0$, then $B_j$ lies in the cone of $D_{\rho_u}P_{j-1}$. Since the degree of $D_{\rho_u}P_{j-1}$ is positive for $j\le s$, the directional derivative of $D_{\rho_u}P_{j-1}$ with respect to $B_j$ is nonzero. These directional derivatives commute, thus $D_{\rho_u}P_j\ne0$. Hence, each $B_j|_W$ is in the open cone of the product of the previous restrictions.
\end{proof}

\section{Sequential positivity, its boundary case, and ordinary \texorpdfstring{$m$}{m}-positivity}\label{seq:section}\label{sec:ordinary}
Almost-good prefixes have product $\Theta$ which is strongly positive. If $\Theta\ne0$, then $\Gamma_\Theta$ polarizes $M_\Theta$. The proofs can be found in Section~\ref{seq:proof-section}. Moreover, bidegree $(0,2)$ relations are given by the Moore determinant, and ordinary mixed relations by restricting to hyperplanes.

\subsection{Hyperbolicity cones and polarizing cones}\label{seq:polarizing-section}
The descent lemma preserves the polarizing cone when a boundary direction is appended to the prefix. The hyperbolicity cone may grow larger here. Theorem~\ref{gp:thm:boundary} identifies the polarizing cone of the quotient with this larger cone.

\begin{theorem}[The hyperbolicity cone polarizes the exterior quotient]\label{gp:thm:boundary}\label{seq:cone-comparison}\label{seq:general-hr}
Let $\Omega$ be the product of an almost-good prefix of length $r$ on $V$ such that $\Omega\ne0$. Let $n=d-r$ and fix $p,q\ge0$ such that $m=p+q\le n$. For any $\lambda_1,\dots,\lambda_{n-m},\lambda$ in $\Gamma_\Omega$, we can define $C=\lambda_1\wedge\dots\wedge\lambda_{n-m}$. Then
\begin{equation}\label{gp:eq:boundary-radical}
 \rad Q_{\Omega\wedge C}^{p,q}
 =\Ann(\Omega)\cap\Lambda^{p,q}V^*.
\end{equation}
Multiplication by $C$ induces an isomorphism
\[
 M_\Omega^{p,q}\longrightarrow M_\Omega^{n-q,n-p},\qquad
 [\alpha]\longmapsto[C\wedge\alpha],
\]
and $Q_{\Omega\wedge C}^{p,q}$ is positive definite on
\[
 \{[\alpha]\in M_\Omega^{p,q}:\Omega\wedge C\wedge\lambda\wedge\alpha=0\}.
\]
In particular, $\Gamma_\Omega$ polarizes $M_\Omega$, and
\begin{equation}\label{seq:cone-equality}
 \operatorname{Pol}_\Omega(\omega)=\Gamma_\Omega(\omega).
\end{equation}
\end{theorem}
\begin{proof}
For \eqref{gp:eq:boundary-radical} and the mixed relations of $M_\Omega$ see Subsection~\ref{gp:sec:boundary-proof}. If we let $\lambda_i=\lambda$ then every $\lambda\in\Gamma_\Omega$ polarizes $M_\Omega$. Since $\Gamma_\Omega$ is connected and contains $\omega$ it is contained in $\operatorname{Pol}_\Omega(\omega)$. By \eqref{seq:cone-inclusion} $\operatorname{Pol}_\Omega(\omega)\subseteq\Gamma_\Omega$.
\end{proof}

\begin{corollary}[Strengthened descent]\label{seq:strengthened-descent}
Let $\Theta$ be an almost good product. Let $a \in \overline{\Gamma_\Theta}$ be such that $\Theta \wedge a \neq 0$. Let $[\alpha]_{\Theta} \mapsto [\alpha]_{\Theta \wedge a}$ be the projection. This induces an isomorphism
\[
 M_\Theta/\ker(a\colon M_\Theta\to M_\Theta)\simeq M_{\Theta\wedge a}
\]
with the pairing induced by $a$. Moreover,
\[
 \Gamma_\Theta\subseteq\Gamma_{\Theta\wedge a}
 =\operatorname{Pol}_{\Theta\wedge a}(\omega).
\]
The same is true after any number of boundary steps, provided all the partial products are non-zero.
\end{corollary}
\begin{proof}
The kernel of the projection is $\Ann(\Theta\wedge a)/\Ann(\Theta)=\ker(a\colon M_\Theta\to M_\Theta)$. According to Lemma~\ref{gp:lem:hyperbolic-facts}, we have $\Gamma_\Theta\subseteq\Gamma_{\Theta\wedge a}$. Applying Theorem~\ref{gp:thm:boundary} to the extended prefix, we find that $\Gamma_{\Theta\wedge a}=\operatorname{Pol}_{\Theta\wedge a}(\omega)$.
\end{proof}

\begin{theorem}[Strong positivity]\label{seq:strong-general}
Every almost-good product is strongly positive. Let $\Theta_k$ be the product of a good prefix of length $k$. Then there exists $\varepsilon>0$ such that $\Theta_k-\varepsilon\omega^k$ is strongly positive. This implies the existence of positive coefficients $c_\alpha$ and covectors $u_{\alpha,j}\in V^*$ such that
\begin{equation}\label{seq:sp-finite}
 \Theta_k=\sum_{\alpha=1}^N c_\alpha
 \bigwedge_{j=1}^k\bigl(\sqrt{-1}\,u_{\alpha,j}\wedge\overline{u_{\alpha,j}}\bigr).
\end{equation}
\end{theorem}
See Subsection~\ref{seq:strong-section} for the proof.

\begin{corollary}[Good prefixes]\label{seq:ordinary-good}
Let $(\eta_1,\ldots,\eta_r)$ be a good prefix and let $p+q\le d-r$.  Choose $\lambda_1,\ldots,\lambda_{d-r-p-q},\lambda\in\Gamma_{\Theta_r}$ and set $\Omega=\Theta_r\wedge\lambda_1\wedge\cdots\wedge\lambda_{d-r-p-q}$.  Then $\Omega$ satisfies hard Lefschetz on $\Lambda^{p,q}V^*$ and $(\Omega,\lambda)$ satisfies the Hodge--Riemann relations. This is also true for any strict continuation of the prefix in the relevant bidegrees.
\end{corollary}
\begin{proof}
For $r\ge1$, Theorem~\ref{gp:thm:relative} equates $\rad Q_\Omega^{p,q}$ and $\Ann(\eta_1)^{p,q}$. Since $\eta_1>0$ and $p+q\le d-1$, multiplying by $\eta_1$ defines an injective operator on $\Lambda^{p,q}V^*$ by hard Lefschetz. Therefore, $Q_\Omega^{p,q}$ is non-degenerate and, using Theorem~\ref{gp:thm:relative} once more, positive definite on $P^{p,q}(\Omega,\lambda)$. For $r=0$, this statement is the mixed theorem.
\end{proof}

\begin{corollary}\label{seq:blocks}\label{seq:blocks-section}
To any semipositive good prefix with product $\Theta$, any number of forms $\eta\in\Gamma_\Theta$ may be appended up to the remaining length of the prefix. Then the mixed relations of Corollary~\ref{seq:ordinary-good} hold for forms in the cone of the extended prefix. For any non-zero product $\Psi$ of semipositive forms, the mixed relations hold for all forms in $\Gamma_\Psi$ on the space $M_\Psi$.
\end{corollary}
\begin{proof}
By Proposition~\ref{seq:hyperbolicity}, each successive cone contains $\Gamma_\Theta$, so the extended prefix is good.  To prove the second part of the statement, note that every semipositive form belongs to the closure of any successive cone. Positive definite forms are contained in these cones. No partial product of $\Psi$ can be zero.  Thus Theorem~\ref{gp:thm:boundary} applies.
\end{proof}

\subsection{Restriction and the relations in low bidegrees}
\begin{proposition}[Hereditary restriction]\label{seq:restriction}
If $(\eta_1,\ldots,\eta_k)$ is good and $W\subseteq V$ has dimension $n$, then
\[
 (\eta_1|_W,\ldots,\eta_{\min(k,n)}|_W)
\]
is good on $W$. In particular, $\Theta_k$ is strictly weakly positive because, restricted to each complex $k$-plane, it is a positive volume form.
\end{proposition}
\begin{proof}
Let $W=\ker\phi$ for some nonzero linear form $\phi$, and set $\rho=\sqrt{-1}\phi\wedge\bar\phi$. Let $0\le j<\min(k,d-1)$. Since $\rho\ge 0$, $\rho\in\overline{\Gamma(F_j,\omega)}$. By Corollary~\ref{seq:complete}, $F_j$ is complete. Hence $D_\rho F_j\not\equiv 0$. By the boundary derivative theorem, the hyperbolicity cone of $D_\rho F_j$ contains $\Gamma(F_j,\omega)$ (see \cite{harvey-lawson}). A positive multiple of $D_\rho F_j$ is given by the polynomial
\[
 A\longmapsto *_{W}\bigl((\Theta_j|_W)\wedge(A|_W)^{d-j-1}\bigr).
\]
Thus $\eta_{j+1}|_W$ belongs to the hyperbolicity cone of $\Theta_j|_W$ with respect to $\omega|_W$. Therefore goodness holds on $W$. Repeating the argument for other hyperplanes gives the result for arbitrary $W$.
\end{proof}

\begin{proposition}[Positive eigenvalues relative to an almost-good prefix]\label{seq:positive-index}
Let $1\le m\le d$ and let $\Theta$ be the nonzero product of an almost good prefix of length $d-m$. Any form $\eta\in\Gamma_\Theta(\omega)$ has at least $m$ positive eigenvalues.
\end{proposition}
\begin{proof}
Let $\eta=\eta_+-\eta_-$ be the spectral decomposition of $\eta$ with respect to $\omega$. Because $\eta_-\ge 0$, and forms that are positive definite are in $\Gamma_\Theta$, we have $\eta_-\in\overline{\Gamma_\Theta}$. Since $\Gamma_\Theta$ is stable under addition by elements in its closure, we have $\eta_+=\eta+\eta_-\in\Gamma_\Theta$.
\[
 0<P_\Theta(\eta_+)=*\bigl(\Theta\wedge\eta_+^m\bigr).
\]
Hence $\eta_+^m$ is nonzero and so the rank of $\eta_+$ is at least $m$.
\end{proof}

Let $\Theta$ be the product of an almost-good prefix of length $k$, and let $1\le m\le d-k$. Adding $d-k-m$ copies of $\omega$ gives an almost-good prefix of product $\Psi = \Theta \wedge \omega^{d-k-m}$. By \eqref{pre:relative-cones}, we have $\Gamma_{\Theta,m}(\omega)=\Gamma_\Psi$. Consequently, any $\eta \in \Gamma_{\Theta,m}(\omega)$ has at least $m$ positive eigenvalues. If $\eta_j\in\Gamma_{\Theta_{j-1}}$, then $\eta_j$ has at least $d-j+1$ positive eigenvalues.

\begin{proposition}\label{seq:low}
Let $p, q \ge 0$, $p + q \le d$, and $r = d - p - q$.  If $(\eta_1,\dots,\eta_r)$ is good, then $\Omega = \Theta_r$ satisfies the hard Lefschetz property if
\[
 \min\{p,q\}=0\qquad\text{or}\qquad(p,q)=(1,1).
\]
If $p + q \ge 1$ and $(\eta_1,\dots,\eta_r,\lambda)$ is good, then $(\Omega,\lambda)$ satisfies the Hodge--Riemann relations as well as the orthogonal Lefschetz decomposition in the bidegrees listed above.
\end{proposition}
\begin{proof}
For $p=q=0$, $\Theta_d$ is a positive volume form. Let $q=0$ and $p>0$. Then, by Theorem~\ref{seq:strong-general}, for some $\varepsilon >0$ the form $\Omega - \varepsilon \omega^{d-p}$ is strongly positive. We set $\Omega=1$ if $r=0$. Pairing of this form with signed square of any $(p,0)$-form is nonnegative. Hence
\[
 Q_\Omega(\alpha,\alpha)
 \ge\varepsilon Q_{\omega^{d-p}}(\alpha,\alpha)>0
 \qquad(0\ne\alpha\in\Lambda^{p,0}V^*).
\]
All $(p,0)$-forms are primitive. Hard Lefschetz is true as $Q_\Omega$ is non-degenerate and using exterior duality. The case $p=0$ can be shown similarly.

As shown in Corollary~\ref{seq:complete}, in the case of bidegree $(1,1)$, the quadratic polynomial $F_{d-2}$ is complete and hyperbolic. Its polarization on real $(1,1)$-forms has signature $(1,d^2-1)$. For $\lambda\in\Gamma(F_{d-2},\omega)$ the line $\R\lambda$ is positive so its orthogonal complement is negative definite. Complexification of this negative definite complement gives $P^{1,1}(\Omega,\lambda)$, so the sign in \eqref{eq:Q} makes $Q_\Omega$ positive definite on this space. Thus, positivity and nondegeneracy of $Q_\Omega$ give the other claims.
\end{proof}

\begin{proposition}[The Moore determinant in bidegree $(0,2)$]\label{seq:moore}
Let $d\ge2$ and $(\eta_1,\dots,\eta_{d-2})$ be a good prefix. Then $Q_{\Theta_{d-2}}$ is positive definite on $\Lambda^{0,2}V^*$. This implies $\Theta_{d-2}$ is hard Lefschetz in bidegree $(0,2)$.
\end{proposition}
\begin{proof}
Let $P$ be the Moore determinant of hyperhermitian $d\times d$ matrices. It has been proved using the quaternionic spectral theorem that $P$ is complete, hyperbolic and its hyperbolicity cone is given by the positive definite hyperhermitian matrices \cite{alesker}. Let a hyperhermitian matrix be expressed as $A+B\mathbf j$, with $A^*=A$ and $B^T=-B$. In terms of these components we have
\[
 \operatorname{Herm}_{\mathbb H}(d)
   =\operatorname{Herm}_{\mathbb C}(d)\oplus\Lambda^2(\mathbb C^d)^*,
 \qquad P(A,0)=\det A.
\]
Since $P(A,0)=\det A$, the polynomials $F_j$ of \eqref{seq:descendant} are successive directional derivatives restricted to $B=0$ in the directions $(\eta_j,0)$. Therefore $(\eta_j,0)$ is contained in the derivative cone of the prior derivative with respect to $(\omega,0)$. By Lemma~\ref{lem:complete}, the 2nd degree derivative is complete and is Lorentzian.

In an $\omega$-unitary coframe put $\beta_B=\sum_{i<j}B_{ij}\,dz_i\wedge dz_j$. We claim that
\begin{equation}\label{seq:moore-identity}
 D_{(\eta_1,0)}\cdots D_{(\eta_{d-2},0)}P(A,B)
 =F_{d-2}(A)-*\bigl(\Theta_{d-2}\wedge\beta_B\wedge\bar\beta_B\bigr).
\end{equation}
Indeed, the Pfaffian formula for the Moore determinant reads
\[
 P(A,B)=(-1)^{d(d-1)/2}
 \operatorname{Pf}\begin{pmatrix}B&A\\-\bar A&\bar B\end{pmatrix}.
\]
Using the correspondence between $A$ and its real $(1,1)$-form, let $\gamma=\beta_B+\bar\beta_B-\sqrt{-1}\,A$. $P$ can then be represented as $(\sqrt{-1})^d*\gamma^d/d!$. After differentiating $d-2$ times, only the $(2,2)$ component of $\gamma^2/2$ remains when wedging with $\Theta_{d-2}$. This is $-A^2/2+\beta_B\wedge\bar\beta_B$. Its scalar factor is $(\sqrt{-1})^d(-\sqrt{-1})^{d-2}=-1$. Thus $*(\Theta_{d-2}\wedge A^2)/2$ is equal to $F_{d-2}(A)$ (see \eqref{seq:descendant}).

From \eqref{seq:moore-identity}, the subspace $A=0$ is orthogonal to the positive vector $(\omega,0)$ of this Lorentzian form, hence on this subspace the form is negative definite. This means that $*(\Theta_{d-2}\wedge\beta_B\wedge\bar\beta_B)>0$ for all non-zero $B$. Thus positivity on $\Lambda^{2,0}V^*$ and, by taking the conjugate, on $\Lambda^{0,2}V^*$ follows. Nondegeneracy of the form concludes the proof.
\end{proof}

\subsection{Ordinary \texorpdfstring{$m$}{m}-positivity}\label{seq:ordinary-consequences}\label{ord:mixed-subsection}
For $\Theta = \omega^{d-m}$, the mixed relations for $\Gamma_m(\omega)$ are stated in Corollary~\ref{seq:ordinary-good}. The next argument utilizes restriction to a hyperplane.

\begin{lemma}[Orthonormal hyperplanes preserving $m$-positivity]\label{lem:goodhyperplanes}
Let $m < d$ and let $\eta \in \Gamma_m(\omega)$. There is an $\omega$-orthonormal coframe $\phi_1, \dots, \phi_d$ such that $\eta|_{W_i} \in \Gamma_m(\omega|_{W_i})$ where $W_i = \text{ker } \phi_i$.
\end{lemma}
\begin{proof}
Using the discrete Fourier transform on the $\omega$-unitary eigencoframe of $\eta$ allows us to write $\phi_i$ so that each of its squared modulus coordinates is $1/d$. Let $a$ be the eigenvalue vector of $\eta$. Then, by \eqref{eq:compression},
\[
 \sigma_j(\eta|_{W_i})=\frac1d\sum_{k=1}^d\sigma_j(a_{\widehat k})
                      =\frac{d-j}{d}\sigma_j(a)>0
 \qquad(1\le j\le m),
\]
since each degree-$j$ monomial of $\sigma_j(a)$ occurs in exactly $d-j$ summands.
\end{proof}

\begin{theorem}[Mixed ordinary $m$-positive Hodge--Riemann relations]\label{ord:mixed}\label{thm:ordinary-mixed}
Let $p+q\le m$, $1\le m\le d$, and $r=m-p-q$. For $\eta_1,\ldots,\eta_r\in\Gamma_m(\omega)$ put
\[
 \Omega=\omega^{d-m}\wedge\eta_1\wedge\cdots\wedge\eta_r.
\]
The hard Lefschetz theorem in bidegree $(p, q)$ holds for $\Omega$. For $p + q \ge 1$ and $\lambda \in \Gamma_m(\omega)$ we have the Hodge--Riemann relations for $(\Omega, \lambda)$ and also the orthogonal decomposition \eqref{eq:LD}.
\end{theorem}
\begin{proof}
To start, let $\eta \in \Gamma_m(\omega)$ be fixed and let $\eta_i = \eta$ for all $i$. Let $s = p + q$. We will show HL and HR for $\Omega = \omega^{d-m} \wedge \eta^{m-s}$ and $\lambda \in \{\omega, \eta\}$ using induction on the dimension $d$. Within each dimension, we will first show HL for each $m$ and then prove HR. For $s = 0$, $*(\omega^{d-m} \wedge \eta^m) > 0$ by the $m$-positivity assumption. For $m = d$, both $\omega$ and $\eta$ are positive definite and so the mixed theorem applies. This also holds for dimension $d = 1$.

Assume $1\le s\le m<d$ and $\Omega\wedge\alpha=0$. Using $W_i$ from Lemma~\ref{lem:goodhyperplanes}, let $\alpha_i=\alpha|_{W_i}$ and $\Omega_i=(\omega^{d-m-1}\wedge\eta^{m-s})|_{W_i}$. Since $\Omega_i\wedge\omega|_{W_i}=\Omega|_{W_i}$ and $\Omega\wedge\alpha=0$, we have that $\alpha_i$ is primitive with respect to $(\Omega_i,\omega|_{W_i})$. By the Hodge--Riemann relations in dimension $d-1$, we have that $Q_{\Omega_i}(\alpha_i,\alpha_i)\ge0$ and equals zero only when $\alpha_i=0$. Since $\sum_i\sqrt{-1}\,\phi_i\wedge\bar\phi_i=\omega$, with respect to the induced unitary volumes we have
\begin{equation}\label{eq:hyperplanesum}
 \sum_i Q_{\Omega_i}(\alpha_i,\alpha_i)=Q_\Omega(\alpha,\alpha)=0.
\end{equation}
So $\alpha_i=0$ for all $i$. Each exterior basis monomial of $\alpha$ contains at most $s<d$ indices and will survive on some $W_i$. Hence, $\alpha=0$. Since $\Lambda^{p,q}V^*$ and $\Lambda^{d-q,d-p}V^*$ have the same dimension, multiplication by $\Omega$ is an isomorphism.

Let $p,q\ge1$. Let $\lambda\in\{\omega,\eta\}$. Multiplication by $\lambda$ on $\Lambda^{p-1,q-1}V^*$ pulls back $Q_\Omega$ to $-Q_{\Omega\wedge\lambda^2}$. Thus we have
\begin{equation}\label{eq:ordinarylower}
 \Omega\wedge\omega^2=\omega^{d-(m-2)}\wedge\eta^{(m-2)-(s-2)},\qquad
 \Omega\wedge\eta^2=\omega^{d-m}\wedge\eta^{m-(s-2)}.
\end{equation}
Using HL of orders $m-2$ and $m$ on the lower pairings in \eqref{eq:ordinarylower} respectively, we see that they are both non degenerate (the coefficient for order zero is $\omega^d$). Thus, multiplication by $\lambda$ is injective and the image of this operation is non degenerate with its orthogonal complement being $P^{p,q}(\Omega,\lambda)$. When $pq = 0$, then the whole space is primitive. Let us consider the deformation of $\eta$ to $\omega$ via $\eta_t=(1-t)\eta+t\omega\in\Gamma_m(\omega)$. For $\lambda = \omega$ and $\lambda = \eta_t$, the primitive pairing is non degenerate along the path, so its signature remains unchanged. For $t=1$, we get positive definiteness by classical HR. Therefore, the primitive pairing is positive definite for $t=0$ as well.

For the mixed case we look at the exterior quotient of weight $m$, $A = M_{\omega^{d-m}}$.  Using classical HL for $\omega$ and taking the exterior dual yields
\begin{equation}\label{ord:quotient-identification}
 A^{p,q}=\Lambda^{p,q}V^*,\qquad
 A^{m-q,m-p}\xrightarrow[\ \omega^{d-m}\wedge-\ ]{\ \simeq\ }
      \Lambda^{d-q,d-p}V^*\qquad(p+q\le m).
\end{equation}
Letting $\lambda=\eta$ in the repeated form case, it follows that each $\eta\in\Gamma_m(\omega)$ polarizes $A$. Because $\Gamma_m(\omega)$ contains $\omega$ and is convex, each of these forms is in the same polarization cone. Using Lemma~\ref{pre:mixed-module} on $\eta_1, \dots, \eta_r, \lambda$ with the identification \eqref{ord:quotient-identification} gives $Q_\Omega$ as the mixed pairing and $\ker(\Omega\wedge\lambda\wedge-)$ as the primitive space. This proves the mixed relations as well as \eqref{eq:LD}.
\end{proof}

For $p=q=1$, Theorem~\ref{ord:mixed} becomes the mixed Hodge index theorem. This theorem was proved in cohomology by Xiao \cite[Theorem~A]{xiao-index}. When $\eta_1,\ldots,\eta_r$ and $\lambda$ are also semipositive, Theorem~\ref{ord:mixed} provides the pointwise mixed Hodge--Riemann relations of Xiao \cite[Theorems~3.1 and~A]{xiao-mixed}.

\begin{remark}[Omitting the first factor on a hyperplane]\label{seq:four-hyperplanes}
In Lemma~\ref{lem:goodhyperplanes}, one of the $\omega$ factors can be dropped when restricting to $W_i$. In general, omitting the first factor of a good prefix after restriction to a hyperplane need not leave a good prefix. For example, $(I,\diag(1,1,-1/4))$ is good on $\C^3$ since the second factor has $\sigma_1=7/4$ and $\sigma_2=1/2$. If we omit $I$ after restriction to $\ker dz_2$, the remaining form is $(1, -1/4)$, so it is not a good prefix on that hyperplane.

A good quadruple $(\omega,A,B,C)$ on $\C^4$ can be constructed so that for any orthonormal coframe $v_1,\ldots,v_4$, there is some $j$ for which $(A,B,C)|_{\ker v_j}$ is not good. Let $\omega=I$, $\varepsilon=1/100$, and $A=\diag(-1+\varepsilon,3,3,3)$. Choose
\[
 B=\varepsilon I+
 \left(\begin{pmatrix}1/6&1/2\\1/2&-1/2\end{pmatrix}\oplus I_2\right),\qquad
 C=\varepsilon I+
 \left(\begin{pmatrix}0&\ii/2\\-\ii/2&0\end{pmatrix}\oplus 0_2\right).
\]
Indeed, $\sigma_1(A)=8.01$, $\sigma_2(A)=18.09$, and $\sigma_3(A)=0.27$. The densities of $\omega^2\wedge A\wedge B$, $\omega\wedge A\wedge B^2$ and $\omega\wedge A\wedge B\wedge C$ with respect to $\omega^4/4!$ are $18.5106,0.365406$ and $0.185106$ respectively.

Let $v$ be a unit covector. Let $p=|v_1|^2$, $z=\bar v_1v_2$ and $\rho_v=\ii v\wedge\bar v$. If $(A,B,C)|_{\ker v}$ is good, then the densities below are positive. Calculating gives
\begin{align*}
 *(\rho_v\wedge A^3)&=54(3.99p-0.99),\\
 *(\rho_v\wedge A\wedge B^2)&=v^*D v-12.12\Re z,\\
 *(\rho_v\wedge A\wedge B\wedge C)&=\varepsilon v^*U v+6.06\Im z,
\end{align*}
with
\[
 D=\diag(0.1818,0.121402,0.031102,0.031102),\qquad
 U=\begin{pmatrix}9.18&-3\\-3&5.1202\end{pmatrix}\oplus 2.1052 I_2 .
\]
Each identity follows by extracting the coefficient of $t_1t_2t_3$ in
$v^*\operatorname{adj}(t_1M_1+t_2M_2+t_3M_3)v$ for the three indicated factors. Since $D<19\varepsilon I$ and $\|U\|_{\mathrm{op}}\le12.18$, we have $p>1/5$, $\Re z<2\varepsilon$ and $\Im z>-3\varepsilon$.

Now let $v_j=(v_{ij})_i$ be such a unitary coframe. Let $p_j=|v_{1j}|^2$ and $z_j=\bar v_{1j}v_{2j}$. Because the first two rows are orthogonal we must have $\sum z_j=0$. The remaining three columns provide estimates on the values of $z_j$, with $\Re z_j>-6\varepsilon$ and $\Im z_j<9\varepsilon$ for all $j$, contradicting the fact that the second row is a unit vector, since
\[
 1=\sum_j|v_{2j}|^2=\sum_j |z_j|^2/p_j
 <4\cdot5(6^2+9^2)\varepsilon^2=0.234 .
\]
Restriction to $\ker dz_1$, however, is good. It yields the terms $3I_3$, $\diag(-1/2+\varepsilon,1+\varepsilon,1+\varepsilon)$, and $\varepsilon I_3$. The second form has $\sigma_1=3/2+3\varepsilon>0$ and $\sigma_2=3\varepsilon+3\varepsilon^2>0$.
\end{remark}

\section{Weak positivity and the power map}\label{sec:weak-power}\label{sec:gradient}

In Section~\ref{sec:lorentz}, we show that weakly $(m,l)$-positive forms satisfy Hodge--Riemann relations, which determine the signature of the Hessian of $\sigma_m$. We use this signature to show that the power map is a diffeomorphism.

\subsection{Weak positivity and the Hodge--Riemann relations}

\begin{lemma}\label{codim:cones}
Let $0\le l\le m\le d$ and $\lambda\in\R^d$.
\begin{enumerate}
\item[(a)] Let $J$ be the set of indices of the $l$ largest coordinates of $\lambda$. Then $\lambda \in C_{m,l}^d$ if and only if $\lambda_{\widehat J} \in \Gamma_{m-l}^{d-l}$. This is equivalent to saying that a Hermitian form with eigenvalues $\lambda$, restricted to any subspace of codimension $l$, is $(m-l)$-positive.
\item[(b)] The set $C_{m,l}^d$ is open, convex, and closed under positive coordinate increments. For $l<m$, we have $C_{m,l}^d\subseteq C_{m,l+1}^d$. For $l<m$, any element in $C_{m,l}^d$ has at least $m$ positive coordinates.
\item[(c)] Let $\lambda \in C_{m,l}^d$ and let $t = |I|$. For $t\le l$, we have $\lambda_{\widehat I} \in C_{m-t,l-t}^{d-t}$. For $l < t \le m$, we have $\lambda_{\widehat I} \in \Gamma_{m-t}^{d-t}$.
\item[(d)] Membership in $\Gamma_r$ is preserved when a nonpositive coordinate of $\lambda$ is deleted, as long as $r$ is at most the number of coordinates that remain. For $m>l$, if $\lambda\in C_{m,l}^d$ has a zero coordinate, then removing it gives a tuple in $C_{m,l}^{d-1}$.
\item[(e)] For $l\ge2$, $C_{m,l}^d\subseteq C_{m-2,l-2}^d$. For $l=1$ and $m\ge2$, $C_{m,1}^d\subseteq\Gamma_{m-2}^d$.
\end{enumerate}
\end{lemma}
\begin{proof}
Let the coordinates of $\lambda$ be listed in increasing order. By eigenvalue interlacing, any restriction of codimension $l$ has eigenvalues that are, when listed in order, at least as large as the smallest $d-l$ coordinates of $\lambda$. On the span of the corresponding eigenvectors they are equal. Since increasing coordinates of an element in $\Gamma_{m-l}^{d-l}$ results in an element in the same set, this proves (a), first for coordinate subspaces and then for any subspace.

From \eqref{codim:cone}, we see that $C_{m,l}^d$ is the intersection of the inverse images of $\Gamma_{m-l}^{d-l}$ under the coordinate projections $\lambda\mapsto\lambda_{\widehat I}$ for $|I|=l$. Thus, it is convex and open, and increasing any coordinate preserves membership in the cone. Deleting coordinates gives the inclusion in (b). The lower bound on the number of positive coordinates in (b) is given in Proposition~\ref{ord:positive-index}. For (c), use \eqref{codim:cone} if $t\le l$, and use Corollary~\ref{cor:repeated} if $t>l$.

If $(z,-a)\in\Gamma_r$ with $a\ge0$, the recursion
\[
 \sigma_j(z)=\sigma_j(z,-a)+a\sigma_{j-1}(z)
\]
means $\sigma_j(z)>0$ successively for $1\le j\le r$. For the second assertion of (d), we delete any $l$ of the other coordinates first, then use the recurrence relation with $a=0$.

For (e) each squarefree monomial of degree $j$ on the left hand side has coefficient $\binom{d-j}{t}$ in the following identity.
\begin{equation}\label{codim:eq-restore}
 \sum_{|I|=t}\sigma_j(\lambda_{\widehat I})
 =\binom{d-j}{t}\sigma_j(\lambda).
\end{equation}
For $l \ge 2$, delete $l-2$ coordinates of a tuple and denote the remaining coordinates by $z$. By deleting another pair of coordinates from $z$ we obtain a tuple in $\Gamma_{m-l}$. Applying \eqref{codim:eq-restore} to $z$ with $t=2$ we see $z\in\Gamma_{m-l}$. If $l=1$, applying \eqref{codim:eq-restore} with $t=1$ to $\lambda$ yields $\lambda\in\Gamma_{m-1}\subset\Gamma_{m-2}$.
\end{proof}

\begin{theorem}\label{codim:theorem}
Let $0\le l\le m\le d$ and let $\eta$ be weakly $(m,l)$-positive with respect to $\omega$. For $p,q\ge0$ satisfying
\begin{equation}\label{codim:eq-range}
 p+q\le m,\qquad \max\{p,q\}\ge l,
\end{equation}
put $\Omega_{p,q}=\omega^{d-m}\wedge\eta^{m-p-q}$. Then
\[
 \Omega_{p,q}\wedge-\colon\Lambda^{p,q}V^*\longrightarrow\Lambda^{d-q,d-p}V^*
\]
is an isomorphism. On $P^{p,q}(\Omega_{p,q},\omega)$, the form $Q_{\Omega_{p,q}}$ is positive definite, and
\[
 \Lambda^{p,q}V^*
 =P^{p,q}(\Omega_{p,q},\omega)
   \oplus\omega\wedge\Lambda^{p-1,q-1}V^*
\]
is an orthogonal decomposition with respect to $Q_{\Omega_{p,q}}$. We use the convention that $\Lambda^{a,b}V^*=0$ if either $a$ or $b$ is negative.
\end{theorem}

For $l=0$ and $m\ge1$ this is Theorem~\ref{thm:ordinary-mixed} with coefficient factors $\eta$ and with primitive space defined by $\omega$. For $l=1$, the range in \eqref{codim:eq-range} is $1\le p+q\le m$.

\begin{corollary}[The Lorentzian Hessian theorem]\label{thm:lorentz}
Let $2\le m\le d$ and let $a\in\wtG_m^d$ be the eigenvalue tuple of $\eta$ with respect to $\omega$. Then
\[
 \In\nabla^2\sigma_m(a)=(1,d-1).
\]
Letting $\Omega=\omega^{d-m}\wedge\eta^{m-2}$, the quadratic form on the space of real $(1,1)$-forms, given by $\alpha\mapsto *(\Omega\wedge\alpha^2)$, has inertia $(1,d^2-1)$. It follows that $(\Omega,\omega)$ satisfies the hard Lefschetz and Hodge--Riemann properties in degree $(1,1)$.
\end{corollary}
\begin{proof}
By Theorem~\ref{codim:theorem} with $l=p=q=1$, the form $\alpha\mapsto *(\Omega\wedge\alpha^2)$ is negative definite on the real primitive hyperplane. It is positive on $\R\omega$ because $*(\Omega\wedge\omega^2)>0$, since $a\in\Gamma_{m-1}^d$. Being orthogonal, the two subspaces give inertia $(1,d^2-1)$.

Write $H_m=\nabla^2\sigma_m(a)$. Its diagonal entries vanish, and Lemma~\ref{codim:cones}(c) implies
\begin{equation}\label{eq:Hentry}
 (H_m)_{ij}=\sigma_{m-2}(a_{\widehat{i,j}})>0\qquad(i\ne j).
\end{equation}
Let $(\phi_i)$ be an $\omega$-unitary basis of eigenvectors of $\eta$, with dual coframe $(dz_i)$. For a Hermitian matrix $X$, let $x_i=X_{i\bar i}$ and let $\alpha_X=\sqrt{-1}\sum X_{i\bar j}dz_i\wedge d\bar z_j$. Then
\begin{equation}\label{eq:full}
 *(\Omega\wedge\alpha_X^2)=(d-m)!(m-2)!
 \left[x^TH_mx-2\sum_{i<j}\sigma_{m-2}(a_{\widehat{i,j}})|X_{i\bar j}|^2\right].
\end{equation}
By \eqref{eq:Hentry}, the quadratic form is negative definite on the $d(d-1)$-dimensional real subspace of off-diagonal Hermitian matrices. Thus, the restriction of the form to diagonal matrices has signature $(1, d-1)$.
\end{proof}

\begin{remark}[Polynomial terminology and the primitive reference]
\leavevmode\par\nobreak
The multiaffine quadratic polynomial $\sum_{i<j}\sigma_{m-2}(a_{\widehat{i,j}})x_ix_j$ is Lorentzian as defined by Br\"and\'en and Huh \cite{branden-huh}. The quadratic form on Hermitian matrices \eqref{eq:full} also has Lorentzian signature since off-diagonal squares have negative coefficients.

The choice of $\omega$ is crucial in the definition of $P^{1,1}(\Omega,\omega)$. In the case $d=3, m=2$, the eigenvalues $(1,1,-3/4)$ of $\eta$ yield a weakly 2-positive form with $\sigma_2(1,1,-3/4)=-1/2$. Thus $*(\omega\wedge\eta^2)<0$. The signature of the form $\alpha\mapsto *(\omega\wedge\alpha^2)$ is $(1,8)$, so its restriction to the orthogonal complement of $\eta$ is indefinite.
\end{remark}

\subsection{The power map}

Let $X$ be a compact K\"ahler $d$-fold with $d\ge2$. Write $\mathcal K_X$ for its K\"ahler cone and let $\mathcal B_X$ be the cone of classes represented by closed, strictly positive $(d-1,d-1)$-forms. Fu and Xiao show that the map
\begin{equation}\label{eq:FuXiao}
 b\colon\mathcal K_X\longrightarrow\mathcal B_X,\qquad
 [\alpha]\longmapsto[\alpha^{d-1}]
\end{equation}
is injective and give counterexamples to its surjectivity \cite[Proposition~1.1 and Section~3]{fu-xiao}. They also show injectivity in the case of nef and big classes \cite[Theorem~1.2]{fu-xiao}.

We consider the pointwise map
\begin{equation}\label{eq:rootmap}
 \Phi_{m,\omega}\colon\eta\longmapsto\omega^{d-m}\wedge\eta^{m-1}.
\end{equation}
For $m=d$, Michelsohn showed that every strictly positive $(d-1,d-1)$-form has a unique positive $(d-1)$-th root \cite{michelsohn}. We show that for $2\le m\le d$ each such form has a unique preimage in $\widetilde\Gamma_m(\omega)$. Thus the surjectivity of $b$ does not follow from this as the inverse of \eqref{eq:rootmap} does not necessarily map closed forms to closed forms.

\begin{theorem}\label{thm:gradient}
For $2\le m\le d$, the gradient
\[
 \nabla\sigma_m\colon\wtG_m^d\longrightarrow\R_{>0}^d
\]
is a real-analytic diffeomorphism.
\end{theorem}
\begin{proof}
By Corollary~\ref{thm:lorentz}, the Hessian $\nabla^2\sigma_m$ is nonsingular on $\wtG_m^d$. To show properness of $\nabla\sigma_m$, let $K\subseteq\R_{>0}^d$ be any compact set. Let $\varepsilon,M>0$ be such that $\varepsilon\le y_i\le M$ for all $y\in K$ and all $i$. Let $a\in\wtG_m^d$ and $y=\nabla\sigma_m(a)\in K$. Rename coordinates such that $a=(b,t)$ and $t=\max_i a_i$. Since $a\in\Gamma_{m-1}^d$, we have $\sigma_1(a)>0$ and so $t>0$. From Lemma~\ref{codim:cones}(c) we also have $b\in\Gamma_{m-1}^{d-1}$. Counting the monomials in $y_j=\sigma_{m-1}(a_{\widehat j})$, we have
\begin{equation}\label{eq:gradidentity}
 y_d=\sigma_{m-1}(b),\qquad
 \sum_{j<d}y_j=(d-m)y_d+(d-m+1)t\sigma_{m-2}(b).
\end{equation}
Maclaurin's inequalities on $\Gamma_{m-1}^{d-1}$ imply
\[
 (d-m+1)\sigma_{m-2}(b)
 \ge(m-1)\binom{d-1}{m-1}^{1/(m-1)}
             y_d^{(m-2)/(m-1)}.
\]
Substituting this bound into \eqref{eq:gradidentity} and using $\sum_{j<d}y_j\le(d-1)M$, we find
\begin{equation}\label{eq:gradbound}
 0<t\le
 \frac{(d-1)M}{(m-1)\binom{d-1}{m-1}^{1/(m-1)}
                       \varepsilon^{(m-2)/(m-1)}}.
\end{equation}
The same formulas apply for $m=2$, with $\sigma_0=1$. Since $\sigma_1(a)>0$, each coordinate satisfies $-(d-1)t<a_i\le t$. Hence $(\nabla\sigma_m)^{-1}(K)$ is bounded.

For any $i$, $a_{\widehat i}\in\Gamma_{m-1}^{d-1}$ and $\sigma_{m-1}(a_{\widehat i})=y_i\ge\varepsilon$. According to Maclaurin’s inequalities, for $1\le j\le m-1$, we have
\begin{equation}\label{eq:gradboundary}
 \sigma_j(a_{\widehat i})\ge
 \binom{d-1}{j}\left(\frac{\varepsilon}{\binom{d-1}{m-1}}\right)^{j/(m-1)}>0.
\end{equation}
Therefore, any limit point of $(\nabla\sigma_m)^{-1}(K)$ also belongs to $\wtG_m^d=C_{m,1}^d$. Since $K$ is closed, $(\nabla\sigma_m)^{-1}(K)$ is also closed.

Being a proper local diffeomorphism, $\nabla\sigma_m$ has open and closed image in the connected set $\R_{>0}^d$. Thus, $\nabla\sigma_m$ is a covering of $\R_{>0}^d$. Since $\R_{>0}^d$ is simply connected and $\wtG_m^d$ is convex, the covering has a single sheet. Its inverse is real analytic by the inverse function theorem.
\end{proof}

\begin{theorem}[Hermitian Newton map]\label{thm:newton}
Let $2\le m\le d$. With respect to the trace inner product, the gradient
\[
 T_{m-1}(A)=\nabla_A\sigma_m(A)
 =\sum_{j=0}^{m-1}(-1)^j\sigma_{m-1-j}(A)A^j
\]
is a real-analytic diffeomorphism of the cone of weakly $m$-positive Hermitian matrices onto the cone of positive definite Hermitian matrices. In terms of Hermitian forms, $\Phi_{m,\omega}$ is a real-analytic diffeomorphism of $\widetilde\Gamma_m(\omega)$ onto the cone of strictly positive $(d-1,d-1)$-forms.
\end{theorem}
\begin{proof}
Let $B$ be a positive definite matrix. Write $B=\operatorname{diag}(b)$ with respect to a unitary eigenbasis of $B$. By Theorem~\ref{thm:gradient}, let $a=(\nabla\sigma_m)^{-1}(b)$. Then $A=\operatorname{diag}(a)$ is such that $T_{m-1}(A)=B$ since $T_{m-1}(A)_{ii}=\sigma_{m-1}(a_{\widehat i})$. The inverse gradient commutes with permutations, so if $b_i=b_j$, then $a_i=a_j$. Thus $A$ is a scalar on the eigenspaces of $B$ and does not depend on the choice of the unitary eigenbasis. Every weakly $m$-positive solution $A'$ of $T_{m-1}(A')=B$ commutes with $B$, since $T_{m-1}(A')$ is a polynomial in $A'$. In a common eigenbasis, the eigenvalues of $A'$ are determined by $(\nabla\sigma_m)^{-1}(b)=a$. This proves the uniqueness of $A$. By \eqref{eq:full}, the derivative of $T_{m-1}$ is nonsingular and hence its inverse is real analytic. Under the trace pairing, $\Phi_{m,\omega}$ corresponds to $(d-m)!(m-1)!\,T_{m-1}$.
\end{proof}

Let $C\subseteq\R^d$ be a closed cone. Write $C^*=\{y\in\R^d:y\cdot x\ge0\text{ for all }x\in C\}$.
\begin{proposition}\label{prop:dual}
For $2\le m\le d$,
\[
 \nabla\sigma_m(\Gamma_m^d)
 =\operatorname{int}\bigl((\overline{\Gamma_m^d})^*\bigr).
\]
This also holds for Hermitian matrices, where duality is defined by the trace pairing.
\end{proposition}
\begin{proof}
Let $a\in\Gamma_m^d$ and $b\in\overline{\Gamma_m^d}$. G\aa rding’s positivity theorem says that $\nabla\sigma_m(a)\cdot b=D_b\sigma_m(a)\ge0$. Thus $\nabla\sigma_m(\Gamma_m^d)$ is contained in the dual cone. It is open by Theorem~\ref{thm:gradient}.

Now, let $y\in\operatorname{int}((\overline{\Gamma_m^d})^*)$. Then $y\cdot x$ has a positive minimum on the unit sphere in $\overline{\Gamma_m^d}$. Thus, $S=\{x\in\overline{\Gamma_m^d}:y\cdot x=1\}$ is a compact set. Furthermore, $S$ meets $\Gamma_m^d$, and $\sigma_m$ is zero on the boundary of $\Gamma_m^d$. Thus, the maximum of $\sigma_m$ on $S$ is attained at some $x\in\Gamma_m^d$. By Lagrange’s multiplier method, $\nabla\sigma_m(x)$ is a positive multiple of $y$. Taking the scalar product with $x$ and using Euler’s identity we have $\nabla\sigma_m(x)=m\sigma_m(x)y$. Since $\sigma_m$ is homogeneous we get $\nabla\sigma_m((m\sigma_m(x))^{-1/(m-1)}x)=y$. The same result holds for Hermitian matrices, using Theorem~\ref{thm:newton}.
\end{proof}
For $d=3$ and $m=2$, for example, $(\nabla\sigma_2)^{-1}(1/4,1/4,2)=(1,1,-3/4)$ lies in $\wtG_2^3\setminus\Gamma_2^3$.

\subsection{Weak positivity relative to a fixed prefix}\label{seq:weak-section}

Choose a good prefix $(\omega_1,\ldots,\omega_k)$ and define $\Theta=\omega_1\wedge\cdots\wedge\omega_k$ and $m=d-k\ge2$. By \eqref{pre:relative-cones}, $\Gamma_{\Theta,m-1}(\omega)=\Gamma(D_\omega P_\Theta,\omega)$. The cone for the subsolution condition is
\[
 \widetilde\Gamma_\Theta(\omega)
 =\{\eta\in\Gamma_{\Theta,m-1}(\omega):
       \Theta\wedge\eta^{m-1}>0\text{ as a }(d-1,d-1)\text{-form}\}.
\]
For $\Theta=\omega^{d-m}$, this is the ordinary weak cone $\widetilde\Gamma_m(\omega)$.

\begin{proposition}[Hyperplane characterization of the relative weak cone]\label{seq:weak-hyperplanes}
Let $\Theta$ be the product of a good prefix of length $d-m$ for $m \geq 2$. Then
\[
 \eta\in\widetilde\Gamma_\Theta(\omega)
 \quad\Longleftrightarrow\quad
 \eta|_H\in\Gamma_{\Theta|_H}
 \quad\text{for every complex hyperplane }H.
\]
Every such $\eta$ has at least $m$ positive eigenvalues.
\end{proposition}
\begin{proof}
Assume first that $\eta\in\widetilde\Gamma_\Theta(\omega)$ and let $H$ be a hyperplane. For $m\ge3$, $\eta\in\Gamma_{\Theta,m-1}(\omega)$ implies that $(\omega_1,\ldots,\omega_{d-m},\omega,\eta)$ is a good prefix. By Proposition~\ref{seq:restriction} with this prefix, we have
\[
 (\Theta\wedge\eta^j\wedge\omega^{m-1-j})|_H>0
 \qquad(1\le j\le m-2).
\]
For $j=m-1$, the same inequality follows from $\Theta\wedge\eta^{m-1}>0$. These inequalities determine $\Gamma_{\Theta|_H}$. For $m=2$, we have only the inequality for $j=m-1$.

Assume $\eta|_H\in\Gamma_{\Theta|_H}$ for any hyperplane $H$. Then for each $1\le j\le m-1$ the $(d-1,d-1)$-form $\Theta\wedge\eta^j\wedge\omega^{m-1-j}$ is positive on each hyperplane and hence strictly positive. Multiplying by $\omega$ we obtain $\eta\in\Gamma_{\Theta,m-1}(\omega)$. Since $\Theta\wedge\eta^{m-1}>0$, we conclude that $\eta\in\widetilde\Gamma_\Theta(\omega)$.

Let $H$ be a hyperplane orthogonal to an $\omega$-unit eigenvector of $\eta$ with maximum eigenvalue. From Proposition~\ref{seq:positive-index}, $\eta|_H$ has at least $m-1$ positive eigenvalues. The omitted eigenvalue is at least as large as each eigenvalue on $H$, so $\eta$ has at least $m$ positive eigenvalues.
\end{proof}

For a general good prefix, relative subsolutions need not satisfy the Hodge--Riemann relations in bidegree (1,1).

\begin{example}\label{seq:weak-failure}
On $V=\C^5$, put $e_i=\sqrt{-1}\,dz_i\wedge d\bar z_i$ and $\omega=\sum_{i=1}^5e_i$. Set
\[
 \rho=2e_1-e_2+5e_3+41e_4+36e_5,\qquad
 \eta=5e_1+8e_2-8e_3+39e_4+45e_5.
\]
Then $(\omega, \rho)$ is good and $\eta \in \widetilde{\Gamma}_{\Theta}(\omega)$ with $\Theta = \omega \wedge \rho$. But $(\Theta \wedge \eta, \omega)$ does not satisfy the Hodge--Riemann relation in bidegree $(1,1)$.
\end{example}
\begin{proof}
The first four elementary symmetric functions of $(2,-1,5,41,36)$ are
\[
 (\sigma_1,\sigma_2,\sigma_3,\sigma_4)=(83,1941,9077,3658).
\]
Hence $\rho$ is 4-positive with respect to $\omega$. Let us write the volume form as $\omega^5/5!=e_1\wedge\dots\wedge e_5$. After expanding the products in $(e_i)$ we obtain
\[
 *(\Theta\wedge\omega^3)=1992,\qquad
 *(\Theta\wedge\eta\wedge\omega^2)=25236,\qquad
 *(\Theta\wedge\eta^2\wedge\omega)=110684.
\]
Since the last two coefficients are positive, we obtain $\eta\in\Gamma_{\Theta,2}(\omega)$. We write $e_{\widehat I}=\bigwedge_{j\notin I}e_j$, abbreviating $\widehat{\{i\}}$ to $\widehat i$ and $\widehat{\{i,j\}}$ to $\widehat{i,j}$. Therefore we have
\[
 \Theta\wedge\eta^2=2\sum_{i=1}^5 c_i e_{\widehat i},\qquad
 (c_1,c_2,c_3,c_4,c_5)=(6124,214,47996,868,140).
\]
Since all five coefficients are positive, $\eta\in\widetilde\Gamma_\Theta(\omega)$. Moreover, $\eta\notin\Gamma_\Theta(\omega)$ since $P_\Theta(\eta)=*(\Theta\wedge\eta^3)=-622968$.

Put $\Omega=\Theta\wedge\eta$. On the span of $e_1,\ldots,e_5$, the pairing $(\alpha,\beta)\mapsto*(\Omega\wedge\alpha\wedge\beta)$ has matrix
\[
 M=\begin{pmatrix}
 0&3053&3781&228&204\\
 3053&0&3802&216&159\\
 3781&3802&0&524&583\\
 228&216&524&0&68\\
 204&159&583&68&0
 \end{pmatrix}.
\]
Let $v = (2, 4, -1, -13, -21)^T$ and let $\mathbf 1$ be the all 1s vector in $\mathbb{R}^5$. Then $\mathbf 1^TMv = 0$ and $v^TMv = 378$. This means $\gamma=\sum_i v_ie_i\in P^{1,1}(\Omega,\omega)$ and $Q_\Omega(\gamma,\gamma)=-378$.

Next we show that for any $h>0$, the set $\{A\in\Gamma_\Theta(\omega):P_\Theta(A)=h,\ A-\eta\ge0\}$ is bounded. Let $h>0$, and choose $A\in\Gamma_\Theta(\omega)$ such that $P_\Theta(A)=h$ and $S=A-\eta\ge0$. The coefficients of $\Theta$ in the basis $e_i\wedge e_j$ for $i<j$ are sums of pairs of distinct entries of $(2,-1,5,41,36)$ and hence positive, so $\Theta$ is strongly positive. $\Omega$ is positive as well, since on $e_{\widehat{i,j}}$, its coefficients are $M_{ij}>0$ for $i<j$. Since $S$ is semipositive we have $*(\Omega\wedge S^2)\ge0$ and $*(\Theta\wedge S^3)\ge0$. Now we can expand $\Theta\wedge\eta^2$:
\[
 *(\Theta\wedge\eta^2\wedge S)=2\sum_i c_iS_{i\bar i}
 \ge280\,\operatorname{tr}_\omega S.
\]
Expanding $P_\Theta(\eta+S)$ and applying these bounds, we obtain
\begin{align*}
 h&=-622968+3*(\Theta\wedge\eta^2\wedge S)+3*(\Omega\wedge S^2)+*(\Theta\wedge S^3)\\
  &\ge-622968+840\,\operatorname{tr}_\omega S.
\end{align*}
From this we get $\operatorname{tr}_\omega S\le(h+622968)/840$. Since $S \geq 0$, the largest eigenvalue of $S$ must also be at most this upper bound. This shows $\{A\in\Gamma_\Theta(\omega):P_\Theta(A)=h,\ A-\eta\ge0\}$ is bounded.
\end{proof}

The same construction works for $X = (\mathbb{P}^1)^5$. Let $e_i$ be the pullback of the Fubini--Study form of integral 1 on the $i$-th factor. We use the same formulas for $\omega, \rho, \eta$ and $\gamma$. Using a unitary coframe of the product metric, we obtain identical coefficients as in Example~\ref{seq:weak-failure}. Thus, $(\omega, \rho)$ is pointwise good, and $\eta\in\widetilde{\Gamma}_\Theta(\omega)$ pointwise. Since these are closed forms and $\int_X e_1\wedge\cdots\wedge e_5=1$, $[\gamma]$ is primitive relative to $([\Omega], [\omega])$, and
\[
 \mathcal Q_\Omega([\gamma],[\gamma])
 =-\int_X\Omega\wedge\gamma^2=-378.
\]
Hence, the cohomological Hodge--Riemann relations in bidegree $(1,1)$ may not hold on a smooth projective variety even with the relative weak positivity condition.

\section{Proof of the strong positivity and cone comparison theorems}\label{seq:proof-section}

Let $\xi$ be a top degree form.  We set $[\xi]=*\xi$ using compatible positive volume forms on subspaces.  For $u\in V^*$, we set $\rho_u=\ii u\wedge\bar u$.

\begin{lemma}\label{gp:lem:radical-variation}
Let $Q_t=Q_0+tD$ be Hermitian forms on a fixed finite-dimensional
space. If $n_-(Q_t)=n_-(Q_0)$ for both signs of sufficiently small
$t$, then $D(x,x)=0$ for every $x\in\rad Q_0$.
\end{lemma}
\begin{proof}
Let $L$ be a maximal negative subspace for $Q_0$ and let $x$ be in the radical of $Q_0$ other than $0$. For small $t$, $Q_t$ remains negative definite on $L$. Let $x_t$ be $x$ minus its $Q_t$-orthogonal projection on $L$. Since $Q_t(x,L)=tD(x,L)$, we have $Q_t(x_t,x_t)=tD(x,x)+O(t^2)$. This is nonnegative for both signs of sufficiently small $t$, since $n_-(Q_t)=\dim L$. Therefore $D(x,x)=0$.
\end{proof}

\begin{lemma}\label{gp:lem:detection}
Let $\xi\in\Lambda^{a,b}V^*$ with $a,b<d$. If $\xi$ vanishes on a dense set of complex hyperplanes, then $\xi=0$.
\end{lemma}

\begin{proof}
For $u\ne 0$, $\xi$ vanishes on the restriction to $\ker u$ if and only if $\rho_u\wedge\xi = 0$. By hypothesis this is true for dense hyperplanes. Thus it is true for all $u$. Rank one Hermitian forms span $\Lambda^{1,1}V^*$ over $\C$, so $e_i\wedge\bar e_j\wedge\xi=0$ for all $i, j$ with respect to any coframe $(e_i)$. Let $\xi_{I, J}$ be any coefficient. Choose $i$ and $j$ such that $i\notin I$, $j\notin J$.  The coefficient of the product above at $(I\cup\{i\},J\cup\{j\})$ is $\pm\xi_{I,J}$, so $\xi_{I,J}=0$.
\end{proof}

\subsection{Strong positivity and the pure cases}\label{seq:strong-section}

\begin{lemma}\label{gp:lem:rank-deformation}
Suppose that $(H;B_1,\ldots,B_s)$ is relatively good and
$h=\rank H\ge2$. There are real $(1,1)$-forms $R,K$ and real
numbers $a_1,\ldots,a_s$ such that
\[
 K\ge0,\qquad \rank K=1,\qquad K^2=RK=0,
\]
and, for every $t\in[-1,1]$, the tuple
\[
 \bigl(H+tR;B_1+t a_1K,\ldots,B_s+t a_sK\bigr)
\]
is relatively good. $H+R$ and $H-R$ have rank $h-1$. Every prefix product is affine with respect to $t$, so in particular
\begin{align*}
 H\prod_{j=1}^s B_j
 ={}&\frac12(H+R)\prod_{j=1}^s(B_j+a_jK)\\
 &+\frac12(H-R)\prod_{j=1}^s(B_j-a_jK).
\end{align*}
\end{lemma}

\begin{proof}
Choose a coframe where $H$ is of the form $\rho_u+\rho_v+H'$ and $H'$ is the sum of the remaining $h-2$ coordinate forms. Now define
\[
 K=\rho_u,\qquad R=\ii(u\wedge\bar v+v\wedge\bar u).
\]
Then $K^2=RK=0$, $HK\ne0$, and the eigenvalues of the $u,v$ block of $H+tR$ are $1\pm t$. Let $\Omega_j(t)=(H+tR)\prod_{i=1}^j(B_i+t a_iK)$.  As before, we can choose any $a_i$. Then
\[
 \Omega_j(t)=\Omega_j(0)+t\Delta_j,\qquad
 \Delta_0=R,\qquad \Delta_j=\Delta_{j-1}B_j+a_j\Omega_{j-1}(0)K.
\]
Indeed, $RK=K^2=0$ implies inductively that $\Delta_{j-1}K=0$, so multiplying by $B_j+t a_jK$ produces no quadratic term in $t$. Let $m_j=d-j$ and let
\begin{equation}\label{gp:eq:rank-correction}
 a_j=-\frac{[\Delta_{j-1}B_j^{m_j}]}
 {m_j[\Omega_{j-1}(0)KB_j^{m_j-1}]}.
\end{equation}
Set $g_i=D_KP_{\Omega_i(0)}$. The denominator of \eqref{gp:eq:rank-correction} is $g_{j-1}(B_j)$. Since $H,K\ge0$ and $HK\ne0$, $g_0(\omega)>0$. By Lemma~\ref{gp:lem:hyperbolic-facts}, $B_j$ lies in the cone of $g_{j-1}$ whenever $g_{j-1}\ne0$. Thus $g_{j-1}(B_j)>0$. For $j<s$, we also have $g_j=D_{B_j}g_{j-1}/m_j\ne 0$.

Now we use induction on $j$ to show $B_j+t a_jK$ is admissible for each $t\in [-1, 1]$. Since $H+tR\ne 0$ is semipositive, assuming the previous factors are admissible we look at
\[
 z\longmapsto P_{\Omega_{j-1}(t)}(B_j+t a_jK+z\omega).
\]
By hyperbolicity of $P_{\Omega_{j-1}(t)}$, all roots are real. The leading coefficient is also positive. The constant term of this polynomial is independent of $t$ because $K^2=\Delta_{j-1}K=0$ and the first degree term in $t$ vanishes by \eqref{gp:eq:rank-correction}. Since $B_j\in\Gamma(P_{\Omega_{j-1}(0)})$, the roots of this polynomial at $t=0$ are negative. Since the constant term never vanishes, they remain negative on the entire interval $[-1,1]$. Since $\Omega_s(t)$ is affine, $\Omega_s(0)$ is the average of its endpoint values.
\end{proof}

\begin{proof}[Proof of Theorem~\ref{seq:strong-general}]
Let $(H;B_1,\ldots,B_s)$ be relatively good. We first induct on $d$, then on rank of $H$. When $s=0$ we use the spectral decomposition of $H$. For $\rank H=1$, let $H=c\rho_u$. Then the prefix polynomials are positive multiples of those of $B_1|_{\ker u},\ldots,B_s|_{\ker u}$. This forms a good prefix on $\ker u$. By induction on $d$, their product is strongly positive. Wedging a strongly positive decomposition with $c\rho_u$ proves the assertion. For $\rank H\ge2$, we apply the induction on rank to the end products of Lemma~\ref{gp:lem:rank-deformation}.

Let $F=1$ for a good tuple $(B_1,\ldots,B_s)$ and let $F=H$ for a relatively good tuple $(H;B_1,\ldots,B_s)$. By openness, we may choose a small enough positive number $\delta$ so that all tuples with $B_j$ replaced by $B_j\pm\delta\omega$ are still admissible. For $s \ge 1$ we can write
\[
 F\prod_{j=1}^s B_j-\delta^s F\omega^s
 =2^{1-s}\!\sum_{\substack{\sigma\in\{\pm1\}^s\\\prod_j\sigma_j=-1}}
            F\prod_{j=1}^s(B_j+\sigma_j\delta\omega).
\]
All terms in the sum are strongly positive. This shows strict strong positivity when $F=1$. For $F=1$ and $s=0$, the product is empty, and hence 1.

We will now approximate an almost good prefix with good ones. To do this, replace the first factor with $H+\delta\omega$. If $B$ is the last factor, then $B+\delta\omega$ is in the open cone of the prefix before it for any positive value of $\delta$. By induction, approximate the preceding prefix by a good prefix closely enough that the strict inequalities for $B+\delta\omega$ persist. This holds for any positive value of $\delta$. We can take the limit as $\delta\to 0$ to get the result since the strongly positive cone is closed.
\end{proof}

\begin{lemma}\label{gp:lem:pure-kernel}
Let $\Omega$ be a strongly positive $(k,k)$-form, and put $n=d-k$.
For $0\le q\le n$, the pairing
$Q_{\Omega\omega^{n-q}}^{0,q}$ is semipositive and has kernel
$\Ann(\Omega)\cap\Lambda^{0,q}V^*$.
If $\Omega$ is an almost-good product, the same conclusion holds
with $\omega^{n-q}$ replaced by any product of $n-q$ vectors in
$\Gamma(P_\Omega)$.

If $\Omega=HB_1\cdots B_s$ is relatively good and
$s+1=d-q$, then
\[
 \ker Q_\Omega^{0,q}
 =\Ann(H)\cap\Lambda^{0,q}V^*.
\]
All these statements also hold in the conjugate bidegrees.
\end{lemma}

\begin{proof}
We decompose $\Omega$ into positive generators. Let $G=c\rho_{u_1}\cdots\rho_{u_k}$ be a nonzero generator. Then $Q_{G\omega^{n-q}}^{0,q}(\beta,\beta)$ is some positive constant times the squared norm of $\beta$ on the $n$-dimensional plane $\bigcap_{i=1}^k\ker u_i$. Thus if $Q_{\Omega\omega^{n-q}}^{0,q}(\beta,\beta)=0$ we have $G\beta=0$ for all $G$, which implies $\Omega\beta=0$. The other direction is obvious.

Let $\lambda_1,\ldots,\lambda_{n-q}\in\Gamma(P_\Omega)$. Since $\delta>0$ can be made as small as we please, we may assume $\lambda_i-\delta\omega\in\Gamma(P_\Omega)$ for each $i$. Expand the product of $\lambda_i$ using $\lambda_i=(\lambda_i-\delta\omega)+\delta\omega$. Each of these product terms extends the prefix defining $\Omega$ since the derivative cones are nested. Thus the product with $\Omega$ is strongly positive by Theorem~\ref{seq:strong-general}. This yields $Q_{\Omega\prod_i\lambda_i}^{0,q}\ge\delta^{n-q}Q_{\Omega\omega^{n-q}}^{0,q}$. The first assertion shows the kernel of the operator on the right is $\Ann(\Omega)^{0,q}$. Since this is obviously contained in the kernel of the operator on the left, the two kernels are identical.

Let $\Omega=HB_1\cdots B_s$ with $s\ge 1$. From the proof of Theorem~\ref{seq:strong-general}, we average signs to find $Q_\Omega^{0,q}\ge\delta^sQ_{H\omega^s}^{0,q}$. This along with the first assertion for $H$ and the fact that $\Ann(H)\subseteq\Ann(\Omega)$ yields the desired result. The first assertion applies directly to the case $s=0$. Complex conjugation switches bidegrees.
\end{proof}

\subsection{Interior steps after a semipositive first factor}\label{gp:sec:ordinary}

Let $H$ be a semipositive form of rank $h$ on a $d$-dimensional space. We define
\[
 a_h(d;a,b)=\dim\bigl(\Lambda^{a,b}V^*/\Ann(H)^{a,b}\bigr),
 \qquad
 N_h(d;p,q)=\sum_{j\ge1}(-1)^{j-1}a_h(d;p-j,q-j).
\]

For given bidegrees, $a_h$ and $N_h$ depend on $d$ and $h$ because semipositive forms of the same rank differ only by a change of basis. Thus summands with bidegree outside $0\le a,b\le d$ are zero.

\begin{theorem}\label{gp:thm:relative}
Let $(H;B_1,\ldots,B_s)$ be relatively good on $V$, where $H\ge0$ is
nonzero. Put $h=\rank H$ and
\[
 \Omega=H\wedge B_1\wedge\cdots\wedge B_s,
 \qquad m=d-s-1.
\]
For every $p,q$ with $p+q=m$, the form $Q_\Omega^{p,q}$ has radical
$\Ann(H)^{p,q}$ and negative index $N_h(d;p,q)$.
For every $\lambda\in\Gamma(P_\Omega)$, it induces a positive definite
form on
\[
 \frac{\{\alpha\in\Lambda^{p,q}V^*:
              \Omega\wedge\lambda\wedge\alpha=0\}}
      {\Ann(H)^{p,q}}.
\]
\end{theorem}

For $H>0$, multiplication by $H$ on $\Lambda^{p,q}V^*$ is injective for $p+q\le d-1$ by the Lefschetz theorem. Thus, for any good prefix and its strict extensions, the pairing on $\Lambda^{p,q}V^*$ is non-degenerate.

\begin{proof}
We proceed by induction on $(m,d)$ in lexicographical order. The cases $p=0$ and $q=0$ are dealt with in Lemma~\ref{gp:lem:pure-kernel}. So let $p,q\ge 1$. We use induction on $h$ to prove:
\begin{equation}\label{gp:eq:ordinary-upper}
 n_-(Q_\Omega^{p,q})\le N_h(d;p,q).
\end{equation}
Let $F$ be a relatively good product of weight $w<m$ with first factor $H_F$. For $a+b\le w$, append $w-a-b$ copies of any element of $\Gamma(P_F)$ to $F$. By induction on $m$, the pairing with the resulting product has radical $\Ann(H_F)^{a,b}$. This radical contains $\Ann(F)^{a,b}$, and $\Ann(H_F)\subseteq\Ann(F)$. Thus,
\begin{equation}\label{gp:eq:ordinary-lower-ann}
 \Ann(F)^{a,b}=\Ann(H_F)^{a,b}.
\end{equation}
The same induction argument together with the exterior duality for the complementary degree implies that the elements of $\Gamma(P_F)$ polarize $M_F$.

If $h=1$ we write $H = c\rho_u$ and $W = \ker u$. The restriction of forms gives an isomorphism of $\Lambda^{*,*}V^*/\Ann(H)$ with $\Lambda^{*,*}W^*$. $B_1|_W, \dots, B_s|_W$ are a good prefix of length $s$. Thus, the negative index is $N_1(d;p,q)$ by induction on $d$ (or the classical Hodge--Riemann relations for $s=0$).

For $h\ge 2$ we use the family of forms $\Omega_t$ defined in Lemma~\ref{gp:lem:rank-deformation} where $H = \rho_u + \rho_v + H'$, $K = \rho_u$ and $R = \ii(u\wedge\bar v + v\wedge\bar u)$ for some $u$ and $v$. Let $W = \ker u$. We set
\[
 H_0=H|_W,\qquad \Psi=\prod_j B_j|_W,\qquad F=H_0\Psi=\Omega_t|_W .
\]
Since $\rank H_0=h-1>0$, the restricted tuple is relatively good by Lemma~\ref{gp:lem:hyperbolic-facts}.

We write $\Lambda^{p,q}V^*$ according to the presence or absence of each of $u$ and $\bar u$. The summand $J=K\Lambda^{p-1,q-1}W^*$ is isotropic and can be paired only with $\Lambda^{p,q}W^*$ by multiplying with $F$. According to \eqref{gp:eq:ordinary-lower-ann}, the rank of this pairing is $j = a_{h-1}(d-1; p-1, q-1)$. The kernel of this pairing on $J$ is contained in $\rad Q_{\Omega_t}$ for all $t$ because $\Omega_t K\beta = KF\beta$. Thus we may divide out this kernel. Splitting off $j$ hyperbolic planes yields the following matrix
\begin{equation}\label{gp:eq:ordinary-blocks}
 \begin{pmatrix}A_t&U_t&V_t\\ U_t^*&B&0\\ V_t^*&0&C\end{pmatrix},
 \qquad B=Q_F^{p-1,q},\quad C=Q_F^{p,q-1},
\end{equation}
In the space $A\oplus u\Lambda^{p-1,q}W^*\oplus\bar u\Lambda^{p,q-1}W^*$, $A=\Ann(F)^{p,q}$ on $W$. The Schur complement contributes no correction to \eqref{gp:eq:ordinary-blocks}, as follows. The block we split off is $\left(\begin{smallmatrix}D&P\\P^*&0\end{smallmatrix}\right)$ where $P$ is invertible. Its inverse has zero upper left block, and the remaining vectors are orthogonal to $J$. The signs of blocks $B$ and $C$ are correct because the normal factor change cancels the Hodge sign change. The block between $B$ and $C$ is zero, because we have the same normal covector twice. Blocks $A_t$, $U_t$ and $V_t$ depend on $t$ in an affine way.

The negative index of the sum of hyperbolic planes and the blocks $B$ and $C$ by induction on $m$ is
\begin{equation}\label{gp:eq:ordinary-fixed-budget}
 j+N_{h-1}(d-1;p-1,q)+N_{h-1}(d-1;p,q-1)=N_{h-1}(d;p,q).
\end{equation}
Adding a coordinate to $H_0$ with value 0 multiplies the generating polynomial of $a_{h-1}$ by $(1+x)(1+y)$. Equation~\eqref{gp:eq:ordinary-fixed-budget} then follows from $N_h(d;a,b)+N_h(d;a-1,b-1)=a_h(d;a-1,b-1)$. At $t=\pm1$, by induction on rank, $n_-(Q_{\Omega_t}^{p,q})$ is bounded by $N_{h-1}(d;p,q)$. If $A$ was non-trivially paired with either $\rad B$ or $\rad C$ at the endpoints of the interval, then an additional hyperbolic plane would increase the negative index. Since $U_t$ and $V_t$ are affine in $t$, the pairings vanish for all $t$. On the quotients by $\rad B$ and $\rad C$, the Schur complement of $B \oplus C$ is
\[
 S_t=A_t-U_tB^{-1}U_t^*-V_tC^{-1}V_t^*.
\]
The endpoint bounds imply $S_{\pm1}\ge0$. Expanding at $t=0$, we obtain
\begin{equation}\label{gp:eq:ordinary-average}
 S_0=\tfrac12(S_{-1}+S_1)+\dot UB^{-1}\dot U^*+\dot VC^{-1}\dot V^* .
\end{equation}

Here the coefficients of $t$ in $U_t$ and $V_t$ are denoted by dots. For any $\alpha\in A$ the functionals $\dot U^*\alpha$ and $\dot V^*\alpha$ are, up to factors of modulus one,
\[
 \beta\longmapsto[\Psi\bar v\alpha\bar\beta],\qquad
 \gamma\longmapsto[\Psi v\alpha\bar\gamma].
\]
Since $H_0\Psi\alpha=0$, they are both zero on multiples of $H_0$. By Lemma~\ref{gp:lem:joint}, the negative index of the sum of the two inverse forms in \eqref{gp:eq:ordinary-average} is at most $N_h(d;p,q)-N_{h-1}(d;p,q)$. Together with equations \eqref{gp:eq:ordinary-fixed-budget} and \eqref{gp:eq:ordinary-average} we have \eqref{gp:eq:ordinary-upper}.

Take $\lambda\in\Gamma(P_\Omega)$. By induction on $m$, $Q_{\Omega\lambda^2}^{p-1,q-1}$ is non-degenerate modulo $\Ann(H)^{p-1,q-1}$. If $H\lambda\beta=0$ for some $\beta$, then $\beta$ lies in $\Ann(H)^{p-1,q-1}$. Multiplying by $\lambda$ gives an inclusion of $M_H^{p-1,q-1}$ into $M_H^{p,q}$. Thus, from \eqref{gp:eq:shift-sign}, $Q_\Omega^{p,q}$ pulls back to $-Q_{\Omega\lambda^2}^{p-1,q-1}$. This implies that $Q_\Omega^{p,q}$ has at least
\[
 a_h(d;p-1,q-1)-N_h(d;p-1,q-1)=N_h(d;p,q)
\]
negative directions. This proves equality in \eqref{gp:eq:ordinary-upper}.

Let $\alpha\in\rad Q_\Omega^{p,q}$. If $s=0$, exterior duality yields $H\alpha=0$. For $s>0$, write $\Omega = \Omega'B_s$. For every real $(1,1)$-form $\theta$ and all sufficiently small $|t|$, the factor $B_s+t\theta$ is admissible. Therefore, the negative index of $Q_{\Omega'(B_s+t\theta)}^{p,q}$ is the same for all sufficiently small $|t|$. Hence, $Q_{\Omega'\theta}^{p,q}(\alpha,\alpha) = 0$ by Lemma~\ref{gp:lem:radical-variation}. Now take $\theta = \rho_a$ and $W_a = \ker a$. Then $\alpha|_{W_a}$ is a zero of $Q_{\Omega'|_{W_a}}^{p,q}$. Since $\Omega\alpha=0$, this restriction is primitive relative to $(\Omega'|_{W_a}, B_s|_{W_a})$. If $H|_{W_a} \neq 0$, then by Lemma~\ref{gp:lem:hyperbolic-facts}, the restriction is relatively good, and by induction on $d$ we obtain $(H\alpha)|_{W_a} = 0$. This also holds for $H|_{W_a}=0$. Since $p,q \ge 1$ and $p+q \le d-1$, $H\alpha$ has bidegree $(p+1,q+1)$ with $p+1,q+1<d$. From Lemma~\ref{gp:lem:detection}, we obtain $H\alpha = 0$. The reverse containment follows from $\Omega=HB_1\cdots B_s$.

In $M_H^{p,q}$, there is a nondegenerate subspace $\lambda M_H^{p-1,q-1}$ with negative index $N_h(d;p,q)$. By exterior duality, its orthogonal complement is the positive definite primitive subspace.
\end{proof}

\begin{remark}\label{gp:rem:ordinary-dependencies}
The proof of Lemma~\ref{gp:lem:joint} uses the weight $m-1$ polarization on $M_F$. Since we prove the lemma by induction on $m$, we assume the existence of the polarization. Induction on $h$ is used only for $t=\pm1$. For the rank one step and the hyperplane restrictions of $H\alpha$, we use induction on $d$.
\end{remark}

\subsection{The joint estimate in the rank induction}\label{gp:sec:joint}

In the rank-reducing family, put $m=p+q$, $W=\ker u$, $H_0=H|_W$,
$F=H_0\Psi$, $r=h-1$ and $s=d-h$. Split $W=S\oplus T$ so that
$H_0$ is positive definite on $S$ of dimension $r$, $T=\ker H_0$
has dimension $s$, and $v\in S^*$. The coefficient $F$ has weight
$m-1$. By induction, the fixed blocks $B=Q_F^{p-1,q}$ and
$C=Q_F^{p,q-1}$ are nondegenerate modulo $\Ann(H_0)$ in their
respective bidegrees. Their inverse forms act on the conjugate-linear
duals of these quotients.

\begin{lemma}\label{gp:lem:joint}
Let $\mathcal A = \Ann(F) \cap \Lambda^{p,q}W^*$, and let $\dot{U}$ and $\dot{V}$ be the coefficients of $t$ in $U_t$ and $V_t$ in \eqref{gp:eq:ordinary-blocks}. Then
\[
 n_-\bigl(\dot UB^{-1}\dot U^*+\dot VC^{-1}\dot V^*\bigr)
 \le N_h(d;p,q)-N_{h-1}(d;p,q).
\]
If $h\ge m$, then $\dot UB^{-1}\dot U^*+\dot VC^{-1}\dot V^*\ge0$.
\end{lemma}
\begin{proof}
For $m=2$, both $B$ and $C$ are positive definite on their quotients by their radicals by Lemma~\ref{gp:lem:pure-kernel}. For these blocks, $N_h(d;p,q) = N_{h-1}(d;p,q) = 1$.
Assume $m\ge 3$. Take any block $D = Q_F^{a,b}$ on $M_F^{a,b}$. Let $M = H_0 M_F^{a-1,b-1}$. Let us write
$a_k(a,b)=\rank(H_0^k|_{\Lambda^{a,b}W^*})$ and
$N_k(a,b)=\sum_{j\ge1}(-1)^{j-1}a_k(a-j,b-j)$, with out-of-range
terms zero. For a nondegenerate Hermitian form $D$ and a subspace $M$,
\begin{equation}\label{gp:eq:joint-inertia}
 n_-(D^{-1}|_{\Ann M})=n_-(D)-\dim M+n_+(D|_M).
\end{equation}
Indeed, if $D|_M$ has inertia $(a',b',c')$, its $b'$ negative directions
and the $c'$ hyperbolic pairs from its radical leave $n_-(D)-b'-c'$ negative directions on $M^\perp$. We map $M^\perp$ to $\Ann M$ by $x\mapsto D(x,\cdot)$.

Here $\dim M=a_2(a-1,b-1)$ and $D|_M$ pulls back to
$-Q_{FH_0^2}^{a-1,b-1}$. By induction on $m$, $\Gamma(P_F)$ polarizes $M_F$. Since $H_0$ is
in the closure of its polarizing cone, Lemma~\ref{gp:lem:descent} applied twice shows this cone also polarizes $M_{FH_0^2}$ of weight $m-3$.
Using \eqref{gp:eq:expected-index}, only degrees up to $m-5$ affect $N_{M_{FH_0^2}}(a-1,b-1)$. For each homogeneous element $\beta$ of degree at most $m-5$, the element $H_0^2\beta$ has degree at most $m-1$. From \eqref{gp:eq:ordinary-lower-ann}, $FH_0^2\beta=0$ is equivalent to $H_0^3\beta=0$.
Hence
$n_-(Q_{FH_0^2}^{a-1,b-1})=N_3(a-1,b-1)$, and
\begin{equation}\label{gp:eq:joint-exact-index}
 n_-(D^{-1}|_{\Ann M})
 =N_1(a,b)-a_2(a-1,b-1)+N_3(a-1,b-1).
\end{equation}

Let $A_k^r$ be the generating polynomial for the ranks of the map given by multiplication by $H_0^k$
on $\Lambda^{*,*}S^*$, and let $P_t^r=\sum_{i+j=t}\dim P_{H_0}^{i,j}(S)x^iy^j$
be the generating polynomial of the primitive bidegrees of total degree $t$ on $S$.
With $z=xy$, one has
\begin{align}
 A_1^r-(1+z)A_2^r+zA_3^r&=P_{r-1}^r,\label{gp:eq:joint-strings}\\
 A_1^{r+1}-(1+x)(1+y)A_1^r&=P_r^r.\label{gp:eq:joint-branching}
\end{align}
For \eqref{gp:eq:joint-strings}, sum the contributions of the Lefschetz strings to obtain
$A_k^r=\sum_{t=0}^{r-k}P_t^r(1+z+\cdots+z^{r-t-k})$.
For \eqref{gp:eq:joint-branching}, the rank in bidegree $(i,j)$ is the smaller of
$\binom ri\binom rj$ and $\binom r{i+1}\binom r{j+1}$.
By Pascal’s identity, the left hand side has non-zero coefficients only when $i+j=r$. In that case, the coefficient is
$\binom ri\binom rj-\binom r{i-1}\binom r{j-1}$ on that diagonal.

Let $L=z(1+x)^s(1+y)^s/(1+z)$ be expanded near $z=0$. Equating coefficients in \eqref{gp:eq:joint-strings} and \eqref{gp:eq:joint-exact-index} gives $n_-(D^{-1}|_{\Ann M})=[x^ay^b]LP_{r-1}^r$.
By \eqref{gp:eq:joint-branching}, we likewise have
$N_h(d;p,q)-N_{h-1}(d;p,q)=[x^py^q]LP_r^r$.
When $h\ge m$, these coefficients are zero by total degree. For any $\alpha\in\mathcal A$, we have $H_0\Psi\alpha=0$. Thus the functionals $\dot U^*\alpha$ and $\dot V^*\alpha$ vanish on multiples of $H_0$. By \eqref{gp:eq:joint-exact-index}, the inverse forms are semi-positive when evaluated at these functionals. This yields the upper bound when $h\ge m$.

Let $2\le h\le m-1$. Set $\ell=m-r$. The restriction $\Psi_T$ is a good prefix of degree $s-\ell$ on $T$. Also, $D_{H_0}^{r-1}P_{H_0}$ is some positive multiple of the determinant on $T$. Because each nonzero derivative cone contains the preceding cone, admissibility is preserved as we commute the derivatives of $\Psi$. The derivatives in these interior directions remain nonzero through the required order.

For $i+j=r-1$, multiplication by $H_0^2$ kills $P_{H_0}^{i,j}(S)$. Hence the image of
\[
 \bigoplus_{i+j=r-1}P_{H_0}^{i,j}(S)\otimes\Lambda^{a-i,b-j}T^*
\]
is in $M^\perp$. On this subspace, $D$ is the orthogonal sum of pairings of the form
\begin{equation}\label{gp:eq:joint-tensor}
 Q_{H_0,S}^{i,j}\otimes Q_{\Psi_T,T}^{a-i,b-j}.
\end{equation}
If the bidegrees on $S$ are different, at least one of the wedge degrees on $S$ is larger than $r$. If the bidegrees on $S$ are identical, the factors on $S$ have degree $(r, r)$. Thus, the only relevant part is $\Psi_T$. Moving the form on $T$ to the other side of the form on $S$ introduces the sign $(-1)^{(r-1)\ell}$, which is equal to the ratio of the Hodge signs. The first factor of \eqref{gp:eq:joint-tensor} is positive definite on primitive forms. By induction on $m$, since $\ell < m$, the second factor is nondegenerate as well (which is also the case if $\Psi_T$ is empty, because of the classical Hodge--Riemann relations).

For both blocks choose the same maximal negative subspaces $N_T^{c,e}$ for $Q_{\Psi_T,T}^{c,e}$, with dimensions
$N_T(s;c,e)=[x^cy^e]L$ by \eqref{gp:eq:expected-index}. Then
\[
 N_D=\bigoplus_{i+j=r-1}P_{H_0}^{i,j}(S)\otimes N_T^{a-i,b-j}
\]
is negative definite and has dimension $[x^ay^b]LP_{r-1}^r$.
Using \eqref{gp:eq:joint-exact-index}, $N_D$ is a maximal negative subspace of $M^\perp$. Then, split off $N_D$ and express the inverse pairing on $\Ann M$ as
\begin{equation}\label{gp:eq:joint-inverse-split}
 D^{-1}(\xi,\xi)
 =R_D(\xi)-G_D(\xi|_{N_D},\xi|_{N_D}),
 \qquad R_D\ge0,\quad G_D>0,
\end{equation}
also when $M^\perp$ has a radical.

For $\alpha\in\mathcal A$ and $\tau\in N_T^{c,e}$, let
$\gamma_\tau(\alpha)$ be the coefficient of a fixed top form on $T$
in $\Psi\alpha\wedge\bar\tau$. Then $\gamma_\tau(\alpha)$ has bidegree $(i,j)=(p-c,q-e)$ on $S$, where $i+j=r$. Since $H_0\Psi\alpha=0$, it is primitive. On $N_D$, the functionals $\dot U^*\alpha$ and $\dot V^*\alpha$ pair $\gamma_\tau(\alpha)$ with some elements of $v\wedge P_{H_0}^{i-1,j}(S)$ and $\bar v\wedge P_{H_0}^{i,j-1}(S)$, respectively. Multiples of $H_0$ pair with $\gamma_\tau(\alpha)$ to zero. This means that both restrictions factor through the same map
\[
 \Gamma\colon\mathcal A\longrightarrow
 \bigoplus_{i+j=r}P_{H_0}^{i,j}(S)
                \otimes\overline{N_T^{p-i,q-j}}^{\,*},
 \qquad \alpha\longmapsto(\tau\mapsto\gamma_\tau(\alpha)).
\]
The map $\tau\mapsto\gamma_\tau(\alpha)$ is conjugate-linear.
By \eqref{gp:eq:joint-inverse-split}, the sum of the two inverse
forms pulled back to $\mathcal A$ is $R-\Gamma^*G\Gamma$ with
$R,G\ge0$. Hence, its negative index is bounded by the rank of $\Gamma$ which is bounded by
\begin{equation}\label{gp:eq:joint-budget}
 \sum_{i+j=r}\dim P_{H_0}^{i,j}(S)\,N_T(s;p-i,q-j)
 =[x^py^q]LP_r^r=N_h(d;p,q)-N_{h-1}(d;p,q).
\end{equation}
\end{proof}

\begin{remark}
Factors of $\Psi$ do not need to commute as Hermitian endomorphisms. In \eqref{gp:eq:joint-tensor}, only $\Psi_T$ is involved, but other components of $\Psi$ with respect to $S \oplus T$ may affect $\gamma_\tau(\alpha)$.
\end{remark}

\subsection{Boundary factors in arbitrary bidegree}\label{gp:sec:boundary-proof}
\begin{proof}[Proof of Theorem~\ref{gp:thm:boundary}]
Let $b$ be the number of factors beyond the initial factor $H$ that are on the boundary of the hyperbolicity cone of the preceding products, with $H \neq 0$ and semipositive. We induct on the lexicographically ordered triple $(m, b, d)$ where $m = p+q$. For empty prefixes, we can use mixed Hodge--Riemann relations. Pure bidegrees may be reduced to Lemma~\ref{gp:lem:pure-kernel}.  When $b = 0$, we apply Theorem~\ref{gp:thm:relative} and add the factors in $C$ to the prefix. In bidegree $(p,q)$, we have
\[
 \Ann(H)\subseteq\Ann(\Omega)\subseteq\Ann(\Omega C)=\Ann(H).
\]
Assume henceforth that $p,q,b>0$.

Let $\eta$ be the last boundary factor and write $\Omega=\Theta\eta\Psi$.
Since openness holds, for any small enough $\varepsilon>0$, each of the factors of $\Psi$ stays in the open cone of the product of its preceding factors in
\begin{equation}\label{gp:eq:boundary-pencil}
 \Omega_\varepsilon=\Theta(\eta+\varepsilon\omega)\Psi
   =\Omega+\varepsilon\Omega_1,\qquad \Omega_1=\Theta\omega\Psi,
\end{equation}
Also, $\lambda_1, \ldots, \lambda_{n-m}, \lambda$ remain in $\Gamma_{\Omega_{\varepsilon}}$. So the prefix which defines $\Omega_{\varepsilon}$ has $b-1$ steps that are on the boundary. We can let $U_j = A^{p-j, q-j}$ with $A = M_{\Theta}$ and with
$a_{k,j}=\rank(\eta^k|_{U_j})$, and $b_{k,j}=a_{k-1,j}-a_{k,j}$
for $k,j\ge1$, with $a_{0,j}=\dim U_j$ and out-of-range bidegrees taken as zero.
In total degree $e<m$, multiplication by a strict factor on a prefix quotient of weight $w>e$ is injective, since its $(w-e)$th power is a Lefschetz isomorphism by induction on $m$. Removing the strict factors in degree $m-2j$, we get
\begin{equation}\label{gp:eq:lower-hodge-differences}
 \dim M_{\Omega_\varepsilon}^{p-j,q-j}=a_{0,j},\qquad
 \dim M_\Omega^{p-j,q-j}=a_{1,j}.
\end{equation}
By \eqref{gp:eq:expected-index}, the difference $N_{M_{\Omega_\varepsilon}}(p,q)-N_{M_\Omega}(p,q)$ is
\begin{equation}\label{gp:eq:delta-index}
 \Delta=\sum_{j\ge1}(-1)^{j-1}b_{1,j}.
\end{equation}

\emph{Successive kernels.}
For $1\le k\le\min(p,q)$ put
\begin{equation}\label{gp:eq:Bk-definition}
 B_k=\ker(\eta^k|_{U_k})/\ker(\eta^{k-1}|_{U_k}).
\end{equation}
For $k=1$, let $F_\varepsilon=\Omega_\varepsilon C\lambda^2$.
Using the induction hypothesis on $m$, we have $M_{F_\varepsilon}$ is polarized of weight $m-2$. When $a+b\leq m-2$, using strict factors we can identify $A^{a,b}$ with $M_{F_\varepsilon}^{a,b}$. Since $\eta B_1=0$ in $A$,
\begin{equation}\label{gp:eq:first-stage}
 Q_{F_\varepsilon}|_{B_1}=\varepsilon E_1,\qquad
 E_1=Q_{\Omega_1C\lambda^2}|_{B_1}.
\end{equation}
For $k\ge2$ take
\begin{equation}\label{gp:eq:higher-stage}
 G_k=\Theta\eta^{k-1}\Psi C\lambda^2\omega^k,\qquad
 E_k=Q_{G_k}^{p-k,q-k}|_{B_k}.
\end{equation}
If $\Theta\eta^{k-1}=0$, all subsequent $B_j$ for $j\ge k$ vanish.
In cases when $\Theta\eta^{k-1}\ne0$, additional insertions of $\eta$ leave subsequent factors in their open cones, since directional derivatives commute and each nonzero derivative cone contains the preceding cone (Lemma~\ref{gp:lem:hyperbolic-facts}). From the induction hypothesis on $m$, $M_{G_k}$ is polarized of weight $m-2k$. Since the multiplication by strict factors is injective, we get
\begin{equation}\label{gp:eq:stage-quotient}
 \Ann(G_k)=\Ann(\Theta\eta^{k-1})\qquad
 \text{in total degrees at most }m-2k.
\end{equation}
In particular, $E_k$ is well defined.

For $k\ge2$ write $M=M_{G_k}$, and for $k=1$ write $M=M_{F_\varepsilon}$.
The hyperbolicity cones of the products defining $M$ in both cases contain $\Gamma_\Theta$ and polarize $M$ by induction on $m$. Thus $\eta\in\overline{\Gamma_\Theta}$ is in the closure of its polarizing cone. Let $K=\ker(\eta\colon M\to M)$. In bidegree $(p-k, q-k)$ the map $A \to M$ embeds $B_k$ into $K$. In bidegree $(p-k-1, q-k-1)$ it identifies $K^{p-k-1, q-k-1}$ with
\begin{equation}\label{gp:eq:lower-stage-kernel}
 \ker(\eta^k|_{U_{k+1}})/\ker(\eta^{k-1}|_{U_{k+1}}).
\end{equation}
In both bidegrees the identifications with quotients of $A$ hold. Multiplication by $\omega$ maps \eqref{gp:eq:lower-stage-kernel} into $B_k$. Using perfect duality, and the fact that $\eta$ is self adjoint, the radical of the pairing on $K$ in bidegree $(p-k, q-k)$ is $(K \cap \eta M)^{p-k, q-k}$. By Lemma~\ref{gp:lem:central}, $K / (K \cap \eta M)$ is polarized by $\omega$ and in its Lefschetz decomposition, a maximal negative subspace lies in the image of multiplication by $\omega$ from bidegree $(p-k-1, q-k-1)$. Therefore, $B_k$ has a maximal negative subspace for the pairing on $K^{p-k, q-k}$.

If we split off this maximal negative subspace we get a positive semidefinite pairing on $K$. On restriction to $B_k$, the radical of this pairing is $B_k \cap \eta M$. Under the middle degree
identification above, a class in this intersection lifts in $A$ to
$\beta=\eta\gamma\bmod\ker\eta^{k-1}$, with $\gamma\in U_{k+1}$.
Since $\eta^k\beta=0$, we have $\eta^{k+1}\gamma=0$, and the
induced class is zero precisely when $\eta^k\gamma=0$. Conversely
every such $\eta\gamma$ pairs to zero with $B_k$. Hence
\begin{equation}\label{gp:eq:exact-stage-radical}
 \rad E_k=\eta B_{k+1}\simeq B_{k+1}.
\end{equation}
For $k=1$ the positive scalar in \eqref{gp:eq:first-stage} has no
effect on this calculation.

For $j\ge k+1$, we repeat the calculation for bidegree $(p-j, q-j)$ to find that $\dim(K/(K\cap\eta M))^{p-j,q-j}=b_{k,j}-b_{k+1,j+1}$. By \eqref{gp:eq:expected-index}, its middle negative index $\nu_k=n_-(E_k)$ is given by
\begin{equation}\label{gp:eq:nu-stage}
 \nu_k=\sum_{\ell\ge1}(-1)^{\ell-1}
       (b_{k,k+\ell}-b_{k+1,k+\ell+1}).
\end{equation}
Since $n_+(E_k)=b_{k,k}-b_{k+1,k+1}-\nu_k$, summing first in $k$
cancels successive terms in both differences, so
\begin{equation}\label{gp:eq:positive-stage-count}
 \sum_k n_+(E_k)=b_{1,1}-b_{1,2}+b_{1,3}-\cdots=\Delta .
\end{equation}

\emph{The negative index of the pencil.}
For each $k$, let $D_k$ be a subspace of $\ker(\eta^k|_{U_k})$ mapping isomorphically to a maximal positive subspace of $E_k$ in $B_k$. For $\beta\in D_k$, write
\begin{equation}\label{gp:eq:finite-correction}
 S_{k,\varepsilon}\beta=\sum_{j=0}^{k-1}(-\varepsilon)^j
  \eta^{k-1-j}\omega^j\beta,\qquad v_{k,\varepsilon}(\beta)=
  \lambda S_{k,\varepsilon}\beta .
\end{equation}
Multiplying the finite sum by $\eta+\varepsilon\omega$ and using $\Theta\eta^k\beta=0$, we obtain
\begin{equation}\label{gp:eq:exact-cancellation}
 \Omega_\varepsilon S_{k,\varepsilon}\beta
   =(-1)^{k-1}\varepsilon^k\Theta\Psi\omega^k\beta .
\end{equation}
Decreasing each of the two components of the bidegree by $k$ changes the Hodge sign by $(-1)^k$.
Using $S_{k,\varepsilon}\gamma=\eta^{k-1}\gamma+O(\varepsilon)$,
we get
\begin{equation}\label{gp:eq:all-order-sign}
 Q_{\Omega_\varepsilon C}
 (v_{k,\varepsilon}(\beta),v_{k,\varepsilon}(\gamma))
 =-\varepsilon^k E_k(\beta,\gamma)+O(\varepsilon^{k+1}).
\end{equation}
Using \eqref{gp:eq:exact-cancellation} and Hermitian symmetry, the pairing of $v_{k,\varepsilon}(\beta)$ with $v_{\ell,\varepsilon}(\gamma)$ under $Q_{\Omega_\varepsilon C}$ is $O(\varepsilon^{\max(k,\ell)})$. Rescale $v_{k,\varepsilon}(\beta)$ by $\varepsilon^{-k/2}$. Pairings with $k\ne\ell$ tend to zero as do pairings with elements of any given maximal negative subspace for $Q_{\Omega C}$. By \eqref{gp:eq:all-order-sign}, on the subspaces chosen above, the diagonal blocks approach $-E_k$ which is negative definite. Thus, by \eqref{gp:eq:positive-stage-count}, for sufficiently small $\varepsilon>0$ we have
\[
 n_-(Q_{\Omega_\varepsilon C}^{p,q})\ge n_-(Q_{\Omega C}^{p,q})+\Delta.
\]
By induction on $b$, the radical of $Q^{p,q}_{\Omega_\varepsilon C}$ is invariant on deformation of elements of $C$ to a single element of $\Gamma_{\Omega_\varepsilon}$. Thus, the index is $N_{M_{\Omega_\varepsilon}}(p,q)$. From this we subtract \eqref{gp:eq:delta-index} to obtain $n_-(Q_{\Omega C}^{p,q})\le N_{M_\Omega}(p,q)$. To establish the other inequality we use the induction hypothesis on $m$ and deduce that $Q_{\Omega C\lambda^2}^{p-1,q-1}$ is non-degenerate on $M_\Omega^{p-1,q-1}$. Using \eqref{gp:eq:shift-sign}, positive directions for $Q_{\Omega C\lambda^2}^{p-1,q-1}$ go to negative directions for $Q_{\Omega C}^{p,q}$ on multiplication by $\lambda$. The number of such directions is given by \eqref{gp:eq:expected-index} as $N_{M_\Omega}(p,q)$. Thus, we have
\begin{equation}\label{gp:eq:boundary-sharp-index}
 n_-(Q_{\Omega C}^{p,q})=N_{M_\Omega}(p,q).
\end{equation}

If $n=m$, the middle trace pairing has radical $\Ann(\Omega)$ by
exterior duality. For $n > m$, let $C = C' \mu$. If $u \neq 0$ and $|t|$ is small enough, then $\mu + t \rho_u \in \Gamma_{\Omega}$. By \eqref{gp:eq:boundary-sharp-index}, the negative index is constant for small values of $t$. Applying Lemma~\ref{gp:lem:radical-variation}, we get that $Q_{\Omega C' \rho_u}(\alpha,\alpha) = 0$ for each $\alpha \in \mathrm{rad}\, Q_{\Omega C}^{p,q}$. On a general hyperplane $W = \ker u$, the vanishing of $Q_{\Omega C'\rho_u}(\alpha,\alpha)$ implies that $\alpha|_W$ has zero norm with respect to $Q_{(\Omega C')|_W}$. Since $\Omega C\alpha=0$, $\alpha|_W$ is primitive with respect to $((\Omega C')|_W,\mu|_W)$. By Lemma~\ref{gp:lem:hyperbolic-facts}, the prefix restricted to $W$ has at most $b$ boundary steps. Using the induction hypothesis on $(m,b,d)$ on dimension $d-1$, we get $(\Omega \alpha)|_W = 0$. Since both the indices of $\Omega \alpha$ are less than $d$, Lemma~\ref{gp:lem:detection} gives $\Omega \alpha = 0$, thereby proving \eqref{gp:eq:boundary-radical}.

The subspace $\lambda M_{\Omega}^{p-1,q-1} \subseteq M_{\Omega}^{p,q}$ is non-degenerate and its negative index is $N_{M_{\Omega}}(p,q)$. Its orthogonal complement is the primitive space and is positive definite by \eqref{gp:eq:boundary-sharp-index}. Since perfect duality holds, the map induced by multiplication by $C$ is a Lefschetz isomorphism.
\end{proof}

\begin{corollary}\label{gp:cor:last-reference}
In Theorem~\ref{gp:thm:boundary}, the primitive reference may be any
$\lambda\in\Gamma_{\Omega C}$. More generally, the theorem applies
after any admissible strict continuation of an almost-good prefix.
\end{corollary}
\begin{proof}
We can apply the middle degree assertion to $\Omega C$. From \eqref{gp:eq:boundary-radical} we have $\Ann(\Omega C)^{p,q} = \Ann(\Omega)^{p,q}$. Since multiplication by strict factors is injective, we can use this to extend this result to the lower bidegrees.
\end{proof}

\section{Proof of higher Hodge--Riemann relations for weak \texorpdfstring{$(m,l)$}{(m,l)}-positive forms}\label{sec:lorentz}\label{sec:higher}

\subsection{Reduction to two eigenvalues}\label{codim:sec-reduction}
Theorem~\ref{codim:theorem} holds for all $m,p,q$ in \eqref{codim:eq-range}. We show this using lexicographic induction on $(l,d)$. For $l=0$, if $m\ge1$, Theorem~\ref{thm:ordinary-mixed} can be applied with primitive space defined by $\omega$ and the factors following $\omega^{d-m}$ are $\eta$. If $m=0$, then the pairing is $\omega^d$.

By selecting appropriate coordinates we can write $\omega=\sum_i\rho_i$ and $\eta=\sum_i\lambda_i\rho_i$ where $\rho_i=\sqrt{-1}\,dz_i\wedge d\bar z_i$. Let $\mathcal B_d=\mathbb R[\rho_1,\ldots,\rho_d]/(\rho_1^2,\ldots,\rho_d^2)$. For degree $(k,k)$ elements, the weight 0 subspace under the action of the diagonal torus is $\mathcal B_d^k\otimes_{\mathbb R}\mathbb C$. In the monomial basis the pairing is given up to positive multiple by $(-1)^kH_k^m(\lambda)$.
\begin{equation}\label{codim:eq-central}
 (H_k^m)_{I,J}=\begin{cases}\sigma_{m-2k}(\lambda_{(I\cup J)^c}),&I\cap J=\varnothing,\\0,&I\cap J\ne\varnothing,\end{cases}
 \qquad |I|=|J|=k .
\end{equation}
At $\lambda = \bf1$, summing the primitive dimensions $\binom dj-\binom d{j-1}$ separately for the two signs in the Lefschetz decomposition yields the classical inertia $\bigl(\binom{d-1}{k},\binom{d-1}{k-1}\bigr)$. Whenever the induction hypothesis applies, orthogonality of the weight spaces implies that the zero-weight pairing is nondegenerate along the segment to $\bf1$ and therefore has classical inertia.

\begin{lemma}[Inertia of a block matrix]\label{codim:lem-saddle}
Let $A$ be a symmetric form on a finite dimensional real vector space $E$ and $B\colon F\to E^*$ be a linear map with zero kernel. Then the inertia of $\left(\begin{smallmatrix}A&B\\B^*&0\end{smallmatrix}\right)$ is inertia of $A|_{\ker B^*}$ plus $(\dim F,\dim F)$.
\end{lemma}
\begin{proof}
Let $E'$ be a subspace of $E$ which is complementary to $\ker B^*$. There is a perfect pairing between $E'$ and $F$, so $E'\oplus F$ is isomorphic to the direct sum of $\dim F$ hyperbolic planes. Adding vectors in $F$ to vectors in $\ker B^*$, we can identify $(\ker B^*,A|_{\ker B^*})$ with $(E'\oplus F)^\perp$.
\end{proof}

\begin{proposition}\label{codim:prop-two-level-reduction}
Suppose the induction hypotheses hold with $k\ge l\ge1$ and $s=m-2k>0$. If $(-1)^kH_k^m(\lambda)$ does not have classical inertia for some $\lambda\in C_{m,l}^d$, then neither does it for some $\lambda=(1^P,(-a)^Q)\in C_{m,l}^d$ with $a>0$ and $P\ge m$.
\end{proposition}
\begin{proof}
Let us refer to such $\lambda$ as a counterexample. When $d=m$, all eigenvalues of $\eta$ are positive by Lemma~\ref{codim:cones}(b), so classical inertia follows from Theorem~\ref{thm:ordinary-mixed}. Assume $d>m$.
Let $\lambda=(\beta,t)$ and list the elements of the monomial basis of $\mathcal B_d^k$ by first listing elements not containing $\rho_d$ followed by the remaining elements.
\begin{equation}\label{codim:eq-one-saddle}
 H_k^m(\beta,t)=\begin{pmatrix}A_0+tA_1&B\\ B^*&0\end{pmatrix},
 \qquad A_1=H_k^{m-1}(\beta).
\end{equation}
Let $L_\eta$ be multiplication on $\mathcal B_{d-1}$ by $\sum_i \beta_i \rho_i$. Using Euler’s identity for elementary symmetric polynomials we have
\begin{equation}\label{codim:eq-one-identities}
 A_1L_\eta=sB,\qquad L_\eta^*A_1L_\eta=s(s+1)H_{k-1}^{m-1}(\beta).
\end{equation}
According to Lemma~\ref{codim:cones}(c), $\beta$ is in $C_{m-1, l-1}^{d-1}$. This allows us to apply the induction hypothesis on $A_1$ and $H_{k-1}^{m-1}(\beta)$. Thus from \eqref{codim:eq-one-identities} we can conclude that $L_\eta$ and $B$ are injective. The pullback of $(-1)^kA_1$ under $L_\eta$ is a negative scalar multiple of $(-1)^{k-1}H_{k-1}^{m-1}(\beta)$. Thus the pullback is nondegenerate and accounts for all $\binom{d-2}{k-1}$ negative directions of $(-1)^kA_1$. Hence $(-1)^kA_1$ is positive definite on $(\operatorname{im}L_\eta)^\perp=\ker B^*$. Classical inertia as in \eqref{codim:eq-one-saddle} is equivalent to $(-1)^k(A_0 + t A_1) > 0$ on $\ker B^*$ by Lemma~\ref{codim:lem-saddle}. Decreasing $t$ along the cone maintains the failure of classical inertia.

Let $J$ be the set of indices associated with the $l$ largest coordinates in the tuple. Set these coordinates to $\min_{j \in J} \lambda_j$. From Lemma~\ref{codim:cones} (a), (b), the tuple still lies in $C_{m,l}^d$, since the coordinates outside of $J$ remain unmodified, and the minimum value is positive. Rescale the tuple so that this minimum is equal to 1.

Every coordinate of a counterexample is non zero. If a coordinate is zero, then permute the coordinates so that $t=0$. From Lemma~\ref{codim:cones} (d), $\beta \in C_{m,l}^{d-1}$ and $A_0 = H_k^m(\beta)$. Thus, from \eqref{codim:eq-central}, $A_0L_\omega = (d-m)B$ on $\mathcal B_{d-1}^{k-1}$. By the hard Lefschetz theorem, $L_\omega$ is injective. By induction on $d$, $A_0$ is nonsingular. Thus $B$ is injective as well. Now $\ker B^*$ is the $\omega$-primitive subspace of $\mathcal B_{d-1}^k$ on which, by induction, $(-1)^k A_0$ is positive definite.

Vary two coordinates $x$ and $y$ not in $J$, while leaving the other coordinates $\beta$ fixed. Group the monomials of $\mathcal B_d^k$ into 4 classes, based on whether they contain neither coordinate generator indexed by $x,y$, only the first, only the second, or both. Using \eqref{codim:eq-central} we may then write
\begin{equation}\label{codim:eq-pair-block}
 \begin{pmatrix}
 C_0+(x+y)C_1+xyC_2 & Q_0+yQ_1 & Q_0+xQ_1 & J_0\\
 *&0&F&0\\ *&F&0&0\\ J_0^*&0&0&0
 \end{pmatrix},\qquad F=H_{k-1}^{m-2}(\beta).
\end{equation}
All coefficient matrices are functions of $\beta$. Since $\beta\in C_{m-2,\max(l-2,0)}^{d-2}$ by Lemma~\ref{codim:cones}(c), by induction $F$ is nonsingular. For $k\ge2$, the expansion of $L_\eta^2$ in \eqref{codim:eq-central} is $J_0^*L_\eta^2=(s+1)(s+2)H_{k-2}^{m-2}(\beta)$. By induction $H_{k-2}^{m-2}(\beta)$ is nonsingular, hence $J_0$ is injective. For $k=1$ remove the last block. Taking the Schur complement of the middle two blocks in \eqref{codim:eq-pair-block} subtracts $(Q_0+yQ_1)F^{-1}(Q_0+xQ_1)^*$ and its transpose from the first block. By Lemma~\ref{codim:lem-saddle}, after splitting off the hyperbolic planes, the signed pairing reduces to
\begin{equation}\label{codim:eq-residual}
 R(x,y)=V+(x+y)U+xyW\quad\text{on }\ker J_0^* .
\end{equation}
Hyperbolic terms contribute $\binom{d-2}{k-1}+\binom{d-2}{k-2}=\binom{d-1}{k-1}$ negative directions. Thus $(-1)^kH_k^m(\lambda)$ has classical inertia if and only if $R(x,y)>0$.

Let $r=m-l\ge2$. Let $\gamma$ be the coordinates of $\lambda_{\widehat J}$ different from $x$ and $y$, and let $A=\sigma_r(\gamma)$, $B=\sigma_{r-1}(\gamma)$ and $C=\sigma_{r-2}(\gamma)$. Let us fix the level
\begin{equation}\label{codim:eq-level}
 \delta=\sigma_r(\lambda_{\widehat J})=A+B(x+y)+Cxy>0.
\end{equation}
If $C>0$, then on this level set $R(x,y)$ is an affine function of $x+y$. Let $\lambda_{\widehat J}$ vary along a path on this level set starting at the initial tuple. Since $\sigma_r$ vanishes on the boundary of $\Gamma_r$, $\lambda_{\widehat J}$ stays in $\Gamma_r$. If coordinates in $J$ continue to be the largest, then by Lemma~\ref{codim:cones}(a), $\lambda$ continues to be in $C_{m,l}^d$.

Let $0<x<1$ and $y<0$. Eliminating $x$ and $y$ from $\lambda_{\widehat J}\in\Gamma_r$ yields $C>0$ and eliminating $x$ gives $B+Cy>0$, hence $B>0$. Therefore the level set is of the form $(t,y(t))$ where
\[
 y(t)=\frac{\delta-A-Bt}{B+Ct},\qquad
 \max\{0,(\delta-A)/B\}\le t\le1.
\]
This path passes through $(x,y)$. On this path, $y(t)\le0$ and on the interior of the parameter interval $(t+y(t))'=C(t-y(t))/(B+Ct)>0$. Since one of the coordinates is zero at the bottom endpoint of the path, $R$ is positive definite there. If $R$ were positive definite also at $t=1$, then since $R$ is affine in $t+y(t)$ we would have $R(x,y)>0$. Therefore the endpoint $t=1$ is a counterexample. Any counterexample must contain a negative coordinate, since any nonnegative tuple with at least $m$ positive coordinates is $m$-positive. Applying this to each of the coordinates in $(0,1)$ we get a counterexample of the form $(1^P,-a_1,\ldots,-a_Q)$ where $a_i>0$ and $P\ge m$ by Lemma~\ref{codim:cones}(b).

Fix $P, Q, \delta$. On the level set $\sigma_r = \delta > 0$, the set of admissible tuples $(a_1,\ldots,a_Q)$ with $a_i \ge 0$ is compact. Membership in $\Gamma_r$ is closed on this set as $\sigma_r=0$ on $\partial\Gamma_r$. Positive trace condition implies $\sum_i a_i<P-l$. The condition of having classical inertia is open, so counterexamples are in a closed set. Let us find a counterexample in this set for which $\sum_i a_i$ is maximized. None of its coordinates is zero.

Assume $-t$ and $-u$ are coordinates, where $t<u$. Removing these coordinates preserves $\Gamma_r$ by Lemma~\ref{codim:cones}(d), so $A, B, C > 0$. Let $S = B/C$. Positivity of $\sigma_{r-1}$ after deleting either $-t$ or $-u$ implies $t, u < S$. By \eqref{codim:eq-level}, we have $(S-t)(S-u) = D$, where $D > 0$. On this level set we obtain $(0, S-D/S)$ and $(v,v)$ in the nonnegative quadrant with $v = S - \sqrt{D}$. The sum $t+u$ is between these sums, as
\[
 t+u-(S-D/S)=tu/S>0,\qquad
 2v-t-u=(\sqrt{S-t}-\sqrt{S-u})^2>0 .
\]
By the zero-coordinate case, $R(0,-(S-D/S))>0$. Since $R(-t,-u)$ is affine in $t+u$ on this level set and fails to be positive definite at the chosen pair, $R(-v,-v)$ cannot be positive definite. Replacing $(-t,-u)$ by $(-v,-v)$ then increases $\sum_i a_i$, contrary to maximality. Therefore $a_1=\dots=a_Q$.
\end{proof}

\subsection{The binary determinant and completion of the induction}\label{codim:sec-quadrature}\label{codim:sec-completion}
Let $\lambda = (1^P, (-a)^Q)$. Write $\mathcal B_d=\mathcal B_P\otimes\mathcal B_Q$ as a direct sum of tensor products of Lefschetz strings. Let $j_+$ and $j_-$ be primitive degrees and let $J = j_+ + j_-$. On each summand of the decomposition the pairing is given by a tensor product of a nondegenerate primitive pairing and the pairing in degree $k-J$ of
\[
 \mathbb R[X,Y]/(X^{P-2j_++1},Y^{Q-2j_-+1}),
 \qquad \omega=X+Y,\quad\eta=X-aY .
\]
The strings of orthogonal primitive classes are pairwise orthogonal, and on each string, the top pairing is non-zero. Let $J>0$. From $\lambda$ delete $2j_+$ positive and $2j_-$ negative entries. According to Lemma~\ref{codim:cones}(c), the new tuple is in $C_{m-2J,l'}^{d-2J}$ for $l'=\max\{l-2J,0\}<l$. A nonzero summand in the decomposition has $k-J\ge l'$. So, by the induction hypothesis, the pairing on $\mathcal B_{d-2J}^{k-J}$ is nondegenerate. Averaging over permutations within each of the two groups of coordinates defines a self-adjoint projection. Thus the pairing on the image of this projection is nondegenerate, and this pairing is the one mentioned above.

Now let $J=0$. Set $b=d-m$ and $s=m-2k$. Note $P\ge m>2k$. Let $X^{k-i}Y^i$, with $0\le i\le\min(k,Q)$, be a basis in degree $k$. The pairing has entry in the $(i,j)$ position, given by the coefficient of $X^PY^Q$ in $X^{k-i}Y^i X^{k-j}Y^j\omega^b\eta^s$. This is
\begin{equation}\label{codim:eq-binary}
 H_{ij}=c_{Q-i-j},\qquad c_j=[z^j](1+z)^b(1-az)^s.
\end{equation}
By reversing the columns of $H$, we obtain a triangular matrix with 1’s on the diagonal for $Q\le k$. For $Q>k$, let $R=Q-k$. Reversing the columns of $H$ and using the dual Jacobi--Trudi identity, we get
\begin{equation}\label{codim:eq-rectangle}
 (-1)^{k(k+1)/2}\det H=s_{((k+1)^R)}(1^b,(-a)^s).
\end{equation}

\begin{lemma}\label{codim:lem-grassmannian}
Let $U$ be an orthonormal $R$-frame chosen with Haar probability on the space of orthonormal $R$-frames in $\mathbb C^n$, $0\le R\le n$. For any Hermitian endomorphism $T$ on $\mathbb C^n$ and any integer $w\ge0$,
\begin{equation}\label{codim:eq-grassmannian}
 \mathbb E\det(U^*TU)^w=s_{(w^R)}(T)/s_{(w^R)}(1^n).
\end{equation}
Let $n=A+S$ where $A\ge R$. Let $P_S$ be the orthogonal projection onto the last $S$ coordinates and let $\Delta$ be the Vandermonde determinant. With probability 1, the matrix $U^*P_SU$ has $h=\min(R,S)$ non zero eigenvalues. The joint probability density of these eigenvalues is proportional to
\begin{equation}\label{codim:eq-jacobi-density}
\Delta(t)^2\prod_{i=1}^h t_i^{|R-S|}(1-t_i)^{A-R},\qquad 0<t_i<1.
\end{equation}
\end{lemma}
\begin{proof}
Let $u=u_1\wedge\cdots\wedge u_R$ be the unit vector corresponding to $U$. Then $u^{\otimes w}$ belongs to the irreducible representation of highest weight $(w^R)$. From Schur’s lemma it follows that the average of the orthogonal projections onto the lines $\mathbb C u^{\otimes w}$ acts on this representation as multiplication by a scalar. The trace of this average is 1, so the scalar is $1/s_{(w^R)}(1^n)$. Taking the trace of its product with the induced action of $T$ proves \eqref{codim:eq-grassmannian}. The squares of the non zero singular values of the last $S$ rows of $U$ are the non zero eigenvalues of $U^*P_SU$. Applying Forrester’s formula \cite[p. 6868, (3.15)]{forrester-jacobi} with $\beta=2$, we see that the density of these eigenvalues is proportional to \eqref{codim:eq-jacobi-density}. In case $S<R$, the adjoint block formula can be used.
\end{proof}

Apply Lemma~\ref{codim:lem-grassmannian} with $T=\diag(1^b,(-a)^s)$ and $w=k+1$. Put
\begin{equation}\label{codim:eq-parameters}
 r=m-k,\quad N=\min(Q,r),\quad h=N-k=\min(R,s),\quad
 \alpha=|Q-r|,\quad B=P-m,\quad x=(1+a)^{-1}.
\end{equation}
In \eqref{codim:eq-jacobi-density}, we have $|R-s|=\alpha$ and $b-R=B+k\ge0$. Since $\det(U^*TU)=(1+a)^h\prod_i(x-t_i)$, the rectangular Schur polynomial in \eqref{codim:eq-rectangle} is of the same sign as
\begin{equation}\label{codim:eq-integral}
 \mathcal I(x)=\int_{[0,1]^h}\Delta(t)^2
 \prod_{i=1}^h(1-t_i)^k(x-t_i)^{k+1}w_0(t_i)\,dt_i,
 \qquad w_0(t)=t^\alpha(1-t)^B .
\end{equation}

\begin{lemma}\label{codim:lem-cone-polynomial}
Let $p_N$ be the degree-$N$ monic orthogonal polynomial for $w_0$ on $[0,1]$. Then
\begin{equation}\label{codim:eq-cone-polynomial}
 x^N\sigma_r(1^{P-k},(-a)^Q)=\binom{P+Q-k}{r}\,p_N(x).
\end{equation}
Let $\xi_N$ be the largest zero of $p_N$. If $(1^P,(-a)^Q)\in C_{m,l}^d$ and $k\ge l$, then $x>\xi_N$.
\end{lemma}
\begin{proof}
Using the basis $x^{N-j}(x-1)^j$ the polynomial on the left of \eqref{codim:eq-cone-polynomial} has coefficients $\binom{P-k}{r-j}\binom Qj$. By Jacobi’s binomial formula \cite[(18.5.8)]{dlmf}, the coefficients of $P_N^{(B,\alpha)}(2x-1)$ are $\binom{N+B}{N-j}\binom{N+\alpha}{j}$. If $N=r$, the coefficients of the two polynomials agree term by term. If $N=Q$, then the ratio of the coefficients of the polynomial on the left to those of the Jacobi polynomial is $(P-k)!Q!/((Q+B)!r!)$ which does not depend on $j$. Consequently, the left polynomial is a multiple of $p_N$. Its leading coefficient is $\binom{P+Q-k}{r}$ by the Vandermonde identity.

By Lemma~\ref{codim:cones}(c), $(1^{P-k},(-a)^Q)\in\Gamma_r$, since $k\ge l$. As $a$ approaches 0, the tuple remains in $\Gamma_r$ by coordinatewise monotonicity of the cone. Additionally, $x$ approaches 1. Therefore, \eqref{codim:eq-cone-polynomial} implies that $p_N$ has no zeros in the interval $[x,1]$. Since all zeros of $p_N$ lie in $(0,1)$ by orthogonality, this proves $x>\xi_N$.
\end{proof}

\begin{proposition}\label{codim:prop-quadrature}
Given the parameters in \eqref{codim:eq-parameters}, the condition $x>\xi_N$ implies $\mathcal I(x)>0$. This in turn implies that the binary matrix in \eqref{codim:eq-binary} is nonsingular for all admissible tuples.
\end{proposition}
\begin{proof}
The positive weights $w_j$ in Gaussian quadrature for $w_0$ correspond to the zeros $0<\xi_1<\dots<\xi_N<1$ of the polynomial, and the approximation is exact for any polynomial of degree up to $2N-1$ \cite[\S3.5(v)]{dlmf}. The polynomial in the integrand of \eqref{codim:eq-integral} has degree $2(h-1)+k+(k+1)=2N-1$ in each of the variables. Using Gaussian quadrature with respect to each of the variables, we have, since the Vandermonde factor is 0 for repeated nodes and for $x>\xi_N$ the rest of the factors in the integrand are positive,
\begin{equation}\label{codim:eq-positive-certificate}
 \mathcal I(x)=h!\sum_{\substack{I\subseteq[N]\\|I|=h}}
 \Delta(\xi_I)^2\prod_{i\in I}w_i(1-\xi_i)^k(x-\xi_i)^{k+1}>0 .
\end{equation}
For $Q>k$, by Lemma~\ref{codim:lem-cone-polynomial}, $x>\xi_N$. Thus $\det H\ne0$ by \eqref{codim:eq-rectangle}, \eqref{codim:eq-integral}, and \eqref{codim:eq-positive-certificate}. For $Q\le k$, the matrix $H$ has 1 on the diagonal with entries below the diagonal equal to zero after the column reversal.
\end{proof}

All blocks remain nonsingular as $a$ decreases to zero in $\lambda=(1^P,(-a)^Q)$. Since $P\ge m$, the endpoint $(1^P,0^Q)$ is ordinarily $m$-positive. Therefore, $(-1)^kH_k^m(\lambda)$ has classical inertia. Proposition~\ref{codim:prop-two-level-reduction} then shows classical inertia for all $\lambda\in C_{m,l}^d$ with $k\ge l$ and $s>0$. If $s=0$, then \eqref{codim:eq-central} is independent of $\lambda$.

We now fix a bidegree $(p,q)$ as in \eqref{codim:eq-range}. A torus weight is determined by two disjoint sets $A$ and $B$ of indices occurring only holomorphically or only antiholomorphically, respectively. Let $t=|A|+|B|$ and $c=p-|A|=q-|B|$. The weight space has basis given by $dz_A\wedge d\bar z_B\wedge\rho_I$, where $I\subseteq (A\cup B)^c$ has length $c$. Except for a nonzero scalar, the pairing has matrix $H_c^{m-t}(\lambda_{\widehat{A\cup B}})$. When $t=0$, we have already shown that the pairing on the zero-weight space has classical inertia. When $t>0$, by Lemma~\ref{codim:cones}(c), the other eigenvalues are in $C_{m-t,l'}^{d-t}$, where $l'=\max(l-t,0)<l$. Since $\max(p,q)\ge l$, we also have $c\ge l'$. Hence, the pairing is nondegenerate by induction, including the case $c=l'=0$. Since weights are orthogonal, $Q_{\Omega_{p,q}}$ is nondegenerate. Exterior duality shows multiplication by $\Omega_{p,q}$ is an isomorphism. Deforming $\lambda$ to $\mathbf1$ in $C_{m,l}^d$, the pairing remains nondegenerate and hence has classical inertia.

When $p=0$ or $q=0$, the whole space is primitive and the pairing is positive definite. Suppose $p,q\ge1$. The pullback of $Q_{\Omega_{p,q}}$ by the map $\omega\wedge-\colon\Lambda^{p-1,q-1}V^*\to\Lambda^{p,q}V^*$ is $-Q_{\Omega_{p,q}\wedge\omega^2}$ where
\[
 \Omega_{p,q}\wedge\omega^2
 =\omega^{d-(m-2)}\wedge\eta^{(m-2)-(p-1)-(q-1)}.
\]
From Lemma~\ref{codim:cones}(e), we have $\lambda\in C_{m-2,\max(l-2,0)}^d$. Since $\max(p-1,q-1)\ge l-1$, we may apply the induction hypothesis to $Q_{\Omega_{p,q}\wedge\omega^2}$ which then implies it has classical inertia. Nondegeneracy of $Q_{\Omega_{p,q}}$ restricted to the image of $\omega\wedge-\colon\Lambda^{p-1,q-1}V^*\to\Lambda^{p,q}V^*$ follows. Its orthogonal complement is $P^{p,q}(\Omega_{p,q},\omega)$ by exterior duality. Subtracting the inertia on the image of $\omega\wedge-$ from that on the full space shows that $Q_{\Omega_{p,q}}$ is positive definite on this primitive subspace.

\section{From pointwise Hodge--Riemann relations to cohomology}\label{glob:section}

We adapt Dinh and Nguy\^en's construction of primitive representatives
\cite[Propositions~2.2--2.4]{dn}, \cite{dn-abstract} to a real $(1,1)$-form $\lambda$.

\begin{theorem}\label{glob:theorem}
Let $X$ be a compact K\"ahler $d$-fold, and let $p,q\ge0$ with $p+q\le d$.
Let $\Omega$ and $\lambda$ be smooth closed real forms of bidegrees
$(d-p-q,d-p-q)$ and $(1,1)$, respectively. Suppose that at every point of $X$ and for every $0\le j\le\min(p,q)$, $\Omega\wedge\lambda^{2j}$ satisfies hard Lefschetz and the pair
$(\Omega\wedge\lambda^{2j},\lambda)$ satisfies
the Hodge--Riemann relations in bidegree $(p-j,q-j)$. Then cup product with $[\Omega]$ defines an isomorphism
\[
 [\Omega]\cdot-\colon H^{p,q}(X)\xrightarrow{\ \sim\ } H^{d-q,d-p}(X).
\]
Let $P^{p,q}(X;\Omega,\lambda)=\ker([\Omega]\cdot[\lambda]\cdot-)$ in $H^{p,q}(X)$. The pairing $\mathcal Q_\Omega$ obtained from \eqref{eq:Q} by integration is positive definite on $P^{p,q}(X;\Omega,\lambda)$, and we have the orthogonal decomposition
\[
 H^{p,q}(X)=P^{p,q}(X;\Omega,\lambda)\oplus[\lambda]\cdot H^{p-1,q-1}(X).
\]
If a Hodge group has a negative index in the bidegree, it is zero. For $p=q=0$, the positive volume form $\Omega$ is the only requirement.
\end{theorem}

\begin{proof}
Using a fixed Hermitian metric, let $dV$ be the volume form. Let $q_\Omega$ be the pointwise pairing with respect to $dV$. Let $p,q\ge 1$. Express each $(p,q)$-covector as $a = a_0 + \lambda\wedge b$, where $a_0$ is primitive with respect to $(\Omega,\lambda)$. Since $\Omega\wedge\lambda^2\wedge b = \Omega\wedge\lambda\wedge a$, by the inverse Lefschetz map we have $|b| \le c|\Omega\wedge\lambda\wedge a|$ for some constant $c$. Since the primitive subspace is positive, and by orthogonality, we obtain
\begin{equation}\label{glob:coercive}
 |a|^2\le C_1|\Omega\wedge\lambda\wedge a|^2+C_2q_\Omega(a,a).
\end{equation}
The constants in the inequality are uniform by compactness. By conjugation, \eqref{glob:coercive} is also valid in bidegree $(q,p)$. Let $\mathcal D = \sqrt{-1}\partial\bar\partial$. Let $\psi$ be a smooth $(q-1, p-1)$-form. Since $d\Omega = 0$, by Stokes’ theorem $\int_X q_\Omega(\mathcal D\psi, \mathcal D\psi) dV = 0$. By taking $a = \mathcal D\psi$ in \eqref{glob:coercive}, we get
\begin{equation}\label{glob:L2}
 \|\mathcal D\psi\|_{L^2}\le C\|\Omega\wedge\lambda\wedge\mathcal D\psi\|_{L^2}.
\end{equation}
Let $f$ be a smooth closed form, representing a class which is primitive. Since $[\Omega]\cdot[\lambda]\cdot[f]=0$, the $\partial\bar\partial$-lemma allows us to express $\Omega\wedge\lambda\wedge f$ in the form $\mathcal Dg$ for some smooth $(d-q,d-p)$-form $g$. Integration by parts shows that the functional
\[
 \Omega\wedge\lambda\wedge\mathcal D\psi\longmapsto
 \int_X f\wedge\Omega\wedge\lambda\wedge\psi
 =\int_X g\wedge\mathcal D\psi
\]
is well defined and bounded on the subspace of forms $\Omega\wedge\lambda\wedge\mathcal D\psi$ in $L^2_{d-p+1,d-q+1}(X)$ by \eqref{glob:L2}. By the Hahn--Banach theorem and $L^2$-duality there exists a $(p-1,q-1)$-form $u \in L^2$ such that
\begin{equation}\label{glob:correction}
 \Omega\wedge\lambda\wedge\mathcal Du=\Omega\wedge\lambda\wedge f
\end{equation}
holds in the sense of distributions. Thus, $h = f - \mathcal D u$ is primitive and closed.

We define the Hermitian metric on $\Lambda^{p,q}T^*X$ as the orthogonal direct sum of $q_\Omega$ on the primitive subbundle and any smooth Hermitian metric on $\lambda\wedge\Lambda^{p-1,q-1}T^*X$. Smooth Hermitian metrics are arbitrarily chosen on other bidegree bundles as well. For a smooth $(p,q-1)$-form $v$, the primitivity of $h$ implies $(\bar\partial v,h)=\mathcal Q_\Omega(\bar\partial v,h)$. Since $d\Omega=0$ and $\partial h=0$, Stokes’ theorem implies that the latter vanishes as well. Therefore, $\bar\partial h=\bar\partial^*h=0$ in the sense of distributions. Since the Dolbeault Laplacian is elliptic, $h$ must be smooth. Since smooth forms and currents compute the same de Rham cohomology, $[h]=[f]$. If $[f] \neq 0$, pointwise positivity of primitive forms gives $\mathcal Q_\Omega([f],[f])=\int_Xq_\Omega(h,h)\,dV>0$. If either $p=0$ or $q=0$, then any $(p,q)$-form is pointwise primitive, hence the result follows on integration.

If $[\Omega]\cdot\alpha=0$, then $\alpha$ is primitive and $\mathcal Q_\Omega(\alpha,\alpha)=0$. By positivity on the primitive subspace we have $\alpha=0$. From Poincar\'e duality and Hodge symmetry it follows that $\dim H^{p,q}(X)=\dim H^{d-q,d-p}(X)$, hence hard Lefschetz holds. Repeat the argument above for $(\Omega\wedge\lambda^{2j},\lambda)$ in bidegree $(p-j,q-j)$. For $p,q\ge1$, $[\Omega]\cdot[\lambda]^2$ yields an isomorphism $H^{p-1,q-1}(X)\simeq H^{d-q+1,d-p+1}(X)$. Hence, for any $\alpha\in H^{p,q}(X)$, there exists a unique $\beta\in H^{p-1,q-1}(X)$ such that $[\Omega]\cdot[\lambda]^2\cdot\beta=[\Omega]\cdot[\lambda]\cdot\alpha$. Thus $\alpha-[\lambda]\cdot\beta$ is primitive. Orthogonality is verified by the definition of $\mathcal Q_\Omega$.
\end{proof}

The following corollaries provide some classes $([\Omega],[\lambda])$ for which hard Lefschetz and Hodge--Riemann hold on $H^{p,q}(X)$, including the orthogonal primitive decomposition of Theorem~\ref{glob:theorem}. Here $p,q\ge0$ and $X$ is a compact K\"ahler $d$-fold, and all $(1,1)$-forms are real, smooth and closed.

\begin{corollary}\label{glob:ordinary}
Let $p+q\le m\le d$. If $\eta_1,\ldots,\eta_{m-p-q},\lambda$ are $m$-positive with respect to a K\"ahler form $\omega$, then Theorem~\ref{glob:theorem} holds with
\[
 \Omega=\omega^{d-m}\wedge\eta_1\wedge\cdots\wedge\eta_{m-p-q}.
\]
If $\eta$ is $m$-positive and $\Omega=\omega^{d-m}\wedge\eta^{m-p-q}$, one may take $\lambda=\omega$ or $\lambda=\eta$.
\end{corollary}
\begin{proof}
Apply Theorem~\ref{thm:ordinary-mixed} with $2j$ additional copies of $\lambda$ when $(p-j)+(q-j)>0$. In total degree zero the mixed top product is positive by G\aa rding's inequality if $m\ge1$, and equals $\omega^d$ if $m=0$.
\end{proof}

\begin{corollary}\label{glob:weak}
Let $\eta$ be weakly $(m,l)$-positive with respect to a K\"ahler form $\omega$, with $0\le l\le m\le d$. If $p+q\le m$ and $\max\{p,q\}\ge l$, then Theorem~\ref{glob:theorem} may be applied with $\Omega=\omega^{d-m}\wedge\eta^{m-p-q}$ and $\lambda=\omega$. If $l=1$ we may consider any bidegree with $1\le p+q\le m$.
\end{corollary}
\begin{proof}
Let $m_j=m-2j$ and $l_j=\max\{l-2j,0\}$ for $0\le j\le\min(p,q)$. Ordinary cones being nested, Lemma~\ref{codim:cones}(e) shows $C_{m,l}^d\subseteq C_{m_j,l_j}^d$. Since $p+q-2j\le m_j$ and $\max\{p-j,q-j\}\ge l_j$, Theorem~\ref{codim:theorem} can be applied to $\Omega\wedge\omega^{2j}=\omega^{d-m_j}\wedge\eta^{m_j-(p-j)-(q-j)}$ in bidegree $(p-j,q-j)$.
\end{proof}

\begin{corollary}[Sequential Hodge--Riemann relations]\label{glob:sequential}
Let $(\eta_1,\ldots,\eta_{d-p-q},\lambda)$ be a globally good prefix, with $1\le p+q\le d$. Then Theorem~\ref{glob:theorem} holds with $\Omega = \eta_1\wedge\cdots\wedge\eta_{d-p-q}$. In degree zero, any full globally good prefix has positive top intersection.
\end{corollary}
\begin{proof}
Since differentiation with respect to $\lambda$ increases the hyperbolicity cone, we can append $2j$ copies of $\lambda$ to $\Omega$ for each $j$ with $0\le j\le\min(p,q)$. Apply Theorem~\ref{seq:cone-comparison} to $(\Omega\wedge\lambda^{2j},\lambda)$ for each positive total degree. For total degree 0, integrate the positive volume form.
\end{proof}
The corollaries use local Lefschetz isomorphisms on the exterior algebra. For a boundary prefix the dimensions of the graded components of $\Lambda^{*,*}T_x^*X/\Ann(\Omega_x)$ can vary with $x\in X$.

\subsection{The cone of ordinary $m$-positive classes}\label{glob:ordinary-cone}
Let $\chi$ be a K\"ahler form on $X$, and let $\kappa = [\chi]$. Let $\mathcal K_m(X;\chi)$ denote the cone of real $(1,1)$-classes with an $m$-positive representative with respect to $\chi$.

\begin{theorem}[Fu--Li--Zhang--Zhang {\cite[Theorem~10.1]{flzz}}]\label{num:ordinary}
Let $1\le m\le d$. For any irreducible analytic subvariety $Z\subseteq X$ of dimension $p=\dim Z\ge d-m+1$, including $Z=X$, let $I_Z(\alpha)=\int_Z\alpha^{m-d+p}\cdot\kappa^{d-m}$. Let
$\mathcal P_m(X;\kappa)=\{\alpha\in H^{1,1}(X,\R):I_Z(\alpha)>0\text{ for every such }Z\}$.
Then $\mathcal K_m(X;\chi)$ is the connected component containing $\kappa$ of $\mathcal P_m(X;\kappa)$. Equivalently, for any fixed K\"ahler class $\beta$,
\[
 \alpha\in\mathcal K_m(X;\chi)
 \quad\Longleftrightarrow\quad
 I_Z(\alpha+t\beta)>0\quad\text{for all such }Z\text{ and every }t\ge0.
\]
\end{theorem}
When $m=d$, this is the K\"ahler cone theorem of Demailly and P\u aun \cite[Theorem 0.1]{demailly-paun}. Let $X$ be projective, with $2\le m\le d$. If $\int_X\alpha^m\cdot\kappa^{d-m}>0$, then by the result of Chen, Nie and Xu \cite[Corollary 1.5]{cnx}, $\alpha\in\mathcal K_m(X;\chi)$ is equivalent to there existing some $\varepsilon_0>0$ such that
\[
 \int_Z\alpha^{m-d+p}\cdot\kappa^{d-m}\ge\varepsilon_0\int_Z\kappa^p
 \qquad (d-m\le p=\dim Z<d)
\]
holds for all irreducible analytic subvarieties $Z$ in the indicated range of dimensions, with $\varepsilon_0$ independent of $Z$. For $p=d-m$, the inequality is satisfied for $\varepsilon_0\le1$. On the 2-Hessian equation in dimension 3 see also Fu, Zhang and Zhang \cite{fzz}.

\begin{corollary}[Independence of the reference form]\label{num:reference-class}
If $\chi_0,\chi_1$ are cohomologous K\"ahler forms, then
$\mathcal K_m(X;\chi_0)=\mathcal K_m(X;\chi_1)$ for $1\le m\le d$.
Suppose there exists an $m$-positive smooth closed real $(1,1)$-form $\eta$ on $X$ relative to $\chi_0$. Then there is $u\in C^\infty(X,\R)$ such that $(\eta+\sqrt{-1}\partial\bar\partial u)^j\wedge\chi_1^{d-j}>0$ for each $j=1,\ldots,m$.
\end{corollary}
\begin{proof}
From Theorem~\ref{num:ordinary}, both cones are the connected component containing $\kappa=[\chi_0]=[\chi_1]$ in $\mathcal P_m(X;\kappa)$. The existence of $u$ follows from the $\partial\bar\partial$-lemma.
\end{proof}
Corollary~\ref{num:reference-class} answers Xiao’s question \cite[Question 2.11]{xiao-index}. If $X$ is projective and $m\ge 2$, it also follows from the theorem of Chen, Nie and Xu. We will denote this cone by $\mathcal K_m(X;\kappa)$.

\begin{example}\label{num:reference}
On $X=(\mathbb P^1)^d$, let $h_i$ be the hyperplane class of the $i$-th factor, and let $\alpha = \sum_{i<d} h_i - h_d$, $\kappa_s = \sum_{i<d} h_i + s h_d$ for $s > 0$. For $1 \le m < d$, $\alpha \in \mathcal K_m(X;\kappa_s)$ if and only if $s > m/(d-m)$. Let each $h_i$ be represented by the pullback of the Fubini--Study form. Then the representative of $\alpha$ relative to that of $\kappa_s$ has eigenvalues $\lambda=(1,\ldots,1,-1/s)$, and $\sigma_j(\lambda) = \binom{d-1}{j} - s^{-1}\binom{d-1}{j-1}$. Thus for $s > m/(d-m)$ we have $\sigma_j(\lambda)>0$ for $1\le j\le m$. On the other hand, if $s \le m/(d-m)$, then $\sigma_m(\lambda) \le 0$ and $\int_X \alpha^m \cdot \kappa_s^{d-m} \le 0$.
\end{example}
Since $\mathcal K_m(X;\kappa)$ is an open convex cone containing $\kappa$, $\alpha$ is in its closure if and only if $\alpha + \varepsilon \kappa \in \mathcal K_m(X;\kappa)$ for any $\varepsilon > 0$.

\subsection{A complete-intersection criterion}

\begin{proposition}\label{glob:complete-intersection}
Let $X$ be a smooth projective $d$-fold and $H$ be a very ample line bundle.  Set $\kappa = c_1(H)$ and select $1\le m\le d$. Let $Y=H_1\cap\cdots\cap H_{d-m}$ be a smooth complete intersection of generic members of $|H|$. If a smooth closed real $(1,1)$-form $\eta$ restricts to a K\"ahler form on $Y$, or more generally if $[\eta]|_Y$ is a K\"ahler class, then $[\eta]\in\mathcal K_m(X;\kappa)$. Consequently, for each K\"ahler form $\chi\in\kappa$, there exists a function $u\in C^\infty(X,\R)$ such that
\[
 (\eta+\sqrt{-1}\partial\bar\partial u)^j\wedge\chi^{d-j}>0
 \qquad(1\le j\le m).
\]
\end{proposition}
\begin{proof}
Since $[\eta]|_Y$ is a K\"ahler class, the $\partial\bar\partial$-lemma on $Y$ yields a $v\in C^\infty(Y,\R)$ with $\eta|_Y+\sqrt{-1}\partial\bar\partial v>0$. Let us extend $v$ smoothly to $X$ and replace $\eta$ with $\eta+\sqrt{-1}\partial\bar\partial v$. After this change, we have $\eta|_Y>0$. Thus, for $m=d$, the result is established. Let us assume $k=d-m>0$, and choose a K\"ahler form $\chi\in\kappa$. Since $Y$ is compact, and complete intersections vary smoothly, there is $\delta>0$ such that $\eta|_{Y'}>\delta\chi|_{Y'}$ for all $Y'$ in some neighborhood of $Y.$

Let $Z$ be an irreducible subvariety of $X$ of dimension $p\geq k+1$ (possibly $Z=X$). Then a general complete intersection of $k$ members of $|H|$ meets $Z$ properly. Choose such a complete intersection $Y'$ near $Y.$ Since $[Y']=\kappa^k$ and $Z\cdot Y'$ is an effective cycle of dimension $p-k$, for all $t\geq 0$ we have
\[
 \int_Z([\eta]+t\kappa)^{p-k}\cdot\kappa^k
 =\int_{Z\cdot Y'}(\eta+t\chi)^{p-k}
 \ge(\delta+t)^{p-k}\int_Z\kappa^p>0 .
\]
By Theorem~\ref{num:ordinary}, $[\eta]\in\mathcal K_m(X;\kappa)$. Apply Corollary~\ref{num:reference-class} for every K\"ahler form representing $\kappa.$
\end{proof}

\section{Numerical criteria, generalized Hessian equations, and intersection inequalities}\label{num:section}\label{rkt:section}

On a compact K\"ahler $d$-fold $X$, fix a globally good prefix of length $r=d-m$ with product $\Theta$, where $m\ge1$. The empty product is $1$. Call a smooth closed real $(1,1)$-form $\theta$ \emph{admissible} if $\theta_x\in\Gamma_{\Theta_x}$ for every $x\in X$. Call a class admissible if it has an admissible representative. By Proposition~\ref{seq:hyperbolicity}, the cone $\mathcal A_\Theta(X)$ of admissible classes is convex and contains the K\"ahler cone. It is open because sufficiently small harmonic perturbations of an admissible form remain admissible on compact $X$.

\subsection{Numerical criteria for sequential positivity}
Let $\mathcal P_{[\Theta]}$ consist of the real $(1,1)$-classes $\alpha$ such that
\[
 \int_Y[\Theta]\cdot\alpha^{p-r}>0
\]
for any irreducible analytic subvariety $Y\subseteq X$ of dimension $r+1\le p\le d$ (including $Y=X$).  If $\dim Y=r$, then $\int_Y[\Theta]>0$, as in Proposition~\ref{seq:restriction}.

\begin{conjecture}[Numerical characterization of sequential $m$-positivity]\label{num:component}
The cone $\mathcal A_\Theta(X)$ is the connected component of
$\mathcal P_{[\Theta]}$ containing the K\"ahler cone.
\end{conjecture}

If Conjecture~\ref{num:component} is true, then $\mathcal A_\Theta(X)$ is determined by $[\Theta]$. The projective criterion of Chen, Nie and Xu \cite[Corollary~1.5]{cnx} suggests the following uniform criterion.

\begin{conjecture}[Uniform sequential criterion]\label{num:uniform}
With respect to a fixed K\"ahler class $\kappa$, a class $\alpha\in H^{1,1}(X,\R)$ is in $\mathcal A_\Theta(X)$ if and only if there exists $\varepsilon>0$ such that
\begin{equation}\label{num:uniformineq}
 \int_Y[\Theta]\cdot(\alpha+t\kappa)^{p-r}
 \ge\varepsilon(1+t)^{p-r}\int_Y\kappa^p
\end{equation}
for all $t\ge0$ and for any irreducible analytic subvariety $Y\subseteq X$ of dimension $r+1\le p\le d$.
\end{conjecture}

\begin{proposition}\label{num:necessary}
Each admissible class $\alpha$ satisfies \eqref{num:uniformineq} for some $\varepsilon > 0$. Any class satisfying this bound is in the connected component of $\mathcal P_{[\Theta]}$ containing the K\"ahler cone.
\end{proposition}
\begin{proof}
Let $\theta \in \alpha$ be an admissible representative and let $\widehat \kappa \in \kappa$ be a K\"ahler form. By convexity, $\theta_s = (1-s)\theta + s \widehat\kappa$ is also admissible for any $s \in [0,1]$. By Propositions~\ref{seq:hyperbolicity} and~\ref{seq:restriction}, we can append $p-r$ many $\theta_s$’s and restrict to a complex $p$-plane $E$ to obtain
\[
 \frac{(\Theta\wedge\theta_s^{p-r})|_E}{\widehat\kappa^p|_E}>0.
\]
As $[0,1]\times\operatorname{Gr}_p(TX)$ is compact, and there are finitely many values of $p$, there exists some $\varepsilon>0$ smaller than any of the ratios. For $s = t/(1+t)$, integration of $(1+t)^{p-r}\Theta\wedge\theta_s^{p-r}$ on the regular part of $Y$ gives \eqref{num:uniformineq}.

Now suppose $\alpha$ satisfies \eqref{num:uniformineq}. For $t \ge 0$, $\alpha + t\kappa \in \mathcal P_{[\Theta]}$. For sufficiently large $t$, this ray is in the K\"ahler cone.
\end{proof}

By Proposition~\ref{num:necessary}, Conjecture~\ref{num:component} implies Conjecture~\ref{num:uniform}. Whether the converse is also true is open. The component criterion holds for $\Theta=\chi^{d-m}$, where $\chi$ is a K\"ahler form, as proven in Theorem~\ref{num:ordinary}. For $m=1$, it is proven in Proposition~\ref{num:linear}. For $r=0$, the component criterion is precisely the K\"ahler cone theorem of Demailly and P\u aun \cite[Theorem~0.1]{demailly-paun}.

\subsection{Generalized Hessian equations}\label{hess:section}
Let $\omega$ be a K\"ahler form and let $\theta$ be a smooth closed real $(1,1)$-form. For $F\in C^\infty(X,\R)$, we seek a smooth real-valued function $u$ on $X$ satisfying the following conditions.
\begin{equation}\label{num:equation}
 \Theta\wedge\theta_u^m=e^F\omega^d,\qquad
 \theta_u=\theta+\sqrt{-1}\partial\bar\partial u,
 \qquad\theta_u(x)\in\Gamma_{\Theta_x}\quad(x\in X).
\end{equation}
By Stokes' theorem, any solution of \eqref{num:equation} satisfies
\begin{equation}\label{num:normalization}
 \int_Xe^F\omega^d=\int_X\Theta\wedge\theta^m>0.
\end{equation}

\begin{proposition}\label{num:structure}
The operator in \eqref{num:equation} is elliptic on $\Gamma_{\Theta_x}$ and the $m$-th root of this operator is concave there. Any two smooth, admissible solutions differ by a constant.
\end{proposition}
\begin{proof}
For $A\in\Gamma_{\Theta_x}$ the linearization is given by the map $H\mapsto m\Theta_x\wedge A^{m-1}\wedge H$. If $H$ is semipositive of rank 1, then by Propositions~\ref{seq:hyperbolicity} and~\ref{seq:restriction}, $(\Theta_x\wedge A^{m-1})|_{\ker H}$ is positive, hence the linearization is elliptic. According to G\aa rding’s inequality, the map $A\mapsto (*(\Theta_x\wedge A^m))^{1/m}$ is concave on $\Gamma_{\Theta_x}$.

For smooth admissible solutions $u,v$, convexity of the admissible cone implies the line segment between $\theta_u$ and $\theta_v$ lies in the cone. Let $L$ be the average of the linearizations on this segment. Then $L(u-v)=0$. By compactness, $L$ is uniformly elliptic. Therefore, by the maximum principle, $u-v$ is a constant.
\end{proof}

\begin{proposition}[The case $m=1$]\label{num:linear}
If $m=1$, then
\[
 \mathcal A_\Theta(X)=\left\{\alpha\in H^{1,1}(X,\R):\int_X[\Theta]\cdot\alpha>0\right\}.
\]
Let $\nu$ be any positive smooth volume form with $\int_X\nu=\int_X\Theta\wedge\theta$. Then the equation $\Theta\wedge\theta_u=\nu$ has a smooth solution, unique up to an additive constant.
\end{proposition}
\begin{proof}
According to Proposition~\ref{seq:restriction}, $\Theta$ is strictly positive. Therefore, $Lu=*(\Theta\wedge\sqrt{-1}\partial\bar\partial u)$ is elliptic. Since $d\Theta = 0$, integration by parts shows that $L$ is self-adjoint on $L^2(X,\vol_\omega)$. Consequently, the maximum principle implies the kernel of $L$ is the set of constant functions. By the Fredholm alternative, the equation $Lu=*(\nu-\Theta\wedge\theta)$ can be solved if and only if $\int_X\nu=\int_X\Theta\wedge\theta$. When $\int_X\Theta\wedge\theta>0$, we can choose $\nu$ to be a positive multiple of $\omega^d$ with the same integral as $\Theta\wedge\theta$. The resulting solution satisfies $\Theta\wedge\theta_u>0$, so $[\theta]$ is admissible. Conversely, integrating $\Theta\wedge\theta_u>0$ for an admissible representative gives $\int_X\Theta\wedge\theta>0$.
\end{proof}

For $\Theta=\omega^{d-m}$, solvability for any $m$-positive class was established in \cite[Corollary 10.2]{flzz} with the assumption \eqref{num:normalization}. Refer also to \cite{hmw,dk-estimates,dk,szekelyhidi}. For general $\Theta$, the spectral operator estimates in \cite{fang-ma,liu-xu} do not apply directly, since $*(\Theta\wedge\theta_u^m)$ can depend on both the eigenvectors and eigenvalues of $\theta_u$.

\begin{conjecture}[Solvability of the generalized $m$-Hessian equation]\label{num:existence}
Let $[\theta] \in \mathcal A_{\Theta}(X)$. For any $F \in C^{\infty}(X,\mathbb R)$ satisfying \eqref{num:normalization}, equation \eqref{num:equation} has a smooth admissible solution, unique up to addition of a constant.
\end{conjecture}

Together, Conjectures~\ref{num:existence} and~\ref{num:component} would characterize solvability of \eqref{num:equation} for normalized $F$ by the connected-component condition on $[\theta]$.

\subsection{Consequences of the Hodge index theorem}
Assume $m\ge2$. For $a_1,\ldots,a_m\in H^{1,1}(X,\R)$, set
\[
 I_\Theta(a_1,\ldots,a_m)=\int_X[\Theta]\cdot a_1\cdot\ldots\cdot a_m,
 \qquad V_\Theta(a)=I_\Theta(a^m).
\]
In any expression of $I_\Theta$, products represent lists of terms.

\begin{proposition}[Intersection inequalities from the sequential Hodge index theorem]\label{rkt:package}
Let all the classes $a,b,a_i,c_i$ below belong to $\mathcal A_\Theta(X)$.
\begin{enumerate}[label=(\alph*)]
\item Every mixed product $I_\Theta(a_1,\ldots,a_m)$ is positive. The pairing
\[
 (u,v)\longmapsto I_\Theta(u,v,a_3,\ldots,a_m)
\]
on $H^{1,1}(X,\R)$ has signature $(1,h^{1,1}(X)-1)$.
\item For every finite collection $a_1,\ldots,a_N$, the polynomial
\[
 f(t_1,\ldots,t_N)=V_\Theta(t_1a_1+\cdots+t_Na_N)
\]
is Lorentzian and has positive coefficients.
\item For $C=c_3\cdot\ldots\cdot c_m$, one has
\begin{equation}\label{rkt:forwardKT}
 I_\Theta(a,b,C)^2\ge I_\Theta(a,a,C)I_\Theta(b,b,C),
\end{equation}
where equality holds if and only if $a$ and $b$ are proportional. The sequence $(I_\Theta(a^{m-j}, b^j))$ is log-concave for $0 \le j \le m$, and $V_\Theta^{1/m}$ is concave on $\mathcal A_\Theta(X)$. Homogeneity of $V_\Theta^{1/m}$ means this concavity is equivalent to the Brunn--Minkowski inequality.
\[
 V_\Theta(a+b)^{1/m}\ge V_\Theta(a)^{1/m}+V_\Theta(b)^{1/m}.
\]
\end{enumerate}
The inequalities extend by continuity to $\overline{\mathcal A_\Theta(X)}$.
\end{proposition}
\begin{proof}
By Proposition~\ref{seq:hyperbolicity}, we can append admissible representatives to the prefix in any order. G\aa rding’s inequality implies the mixed products are positive. The pairing in part (a) has the signature given in Corollary~\ref{glob:sequential} in bidegree $(1,1)$.

All coefficients of monomials of degree $m$ are positive in part (b). Up to a positive factor, the Hessian of any coordinate derivative of order $m-2$ of $f$ is given by the pullback of the pairing in part (a) by the map $t \mapsto \sum t_i a_i$. It has at most one positive eigenvalue. Therefore, according to Br\"and\'en and Huh, $f$ is Lorentzian \cite{branden-huh,hu-xiao-rkt}.

To prove part (c), define $Q(u,v) = I_\Theta(u,v,C)$ and write $b = sa + b_0$ with $Q(a,b_0) = 0$. By (a), $Q$ has signature $(1,h^{1,1}(X)-1)$ and $Q(a,a)>0$, so $Q$ is negative definite on $a^\perp$. Hence
\[
 Q(a,b)^2-Q(a,a)Q(b,b)=-Q(a,a)Q(b_0,b_0)\ge0,
\]
with equality if and only if $b_0=0$. Taking $C=a^{m-j-1}b^{j-1}$ for $1\le j<m$ yields log-concavity. Let $Q_a(u,v)=I_\Theta(a^{m-2},u,v)$. Then the second derivative of $V_\Theta^{1/m}$ at $a$ in the direction of $v$ is a positive multiple of $Q_a(a,a)Q_a(v,v)-Q_a(a,v)^2\le0$.
\end{proof}

In the terms used by Hu and Xiao in \cite[Definition~3.4 and Proposition~3.5]{hu-xiao-shadow}, $V_\Theta$ is strictly $\mathcal A_\Theta(X)$-Lorentzian.  The pointwise Hodge index theorem means $P_\Theta$ is also strictly $\Gamma_\Theta$-Lorentzian.

\subsection{Sharp reverse Khovanskii--Teissier inequalities}
Let $P$ be a homogeneous polynomial of degree $m$. Let $\mathsf P$ be the symmetric polarization of $P$ such that $\mathsf P(a^m)=P(a)$.

\begin{lemma}[Hyperbolic reverse inequality]\label{rkt:hyperbolic}
Let $P$ be hyperbolic of degree $m$ with hyperbolicity cone $\Gamma$. Fix $1\le k\le m-1$. For $a,b_1,\ldots,b_k,c_1,\ldots,c_{m-k}\in\Gamma$,
\begin{align}\label{rkt:polarized}
 &P(a)\,\mathsf P(b_1,\ldots,b_k,c_1,\ldots,c_{m-k})\notag\\
 &\qquad\le\binom{m}{k}
 \mathsf P(a^{m-k},b_1,\ldots,b_k)
 \mathsf P(a^k,c_1,\ldots,c_{m-k}).
\end{align}
The inequality remains true for $b_i,c_j\in\overline\Gamma$ and $a\in\Gamma$. The constant $\binom{m}{k}$ is the best possible constant among all hyperbolic polynomials of degree $m$.
\end{lemma}
\begin{proof}
We use the hyperbolic Rayleigh inequality at $a\in\Gamma$,
\begin{equation}\label{rkt:rayleigh}
 (D_bP)(D_cP)-P D_bD_cP\ge0\qquad(b,c\in\Gamma).
\end{equation}
See \cite[Theorem~3.1]{kpv} for the case where $P$ has no multiple factors. In the general case, we may write $P=\prod_iP_i$ where the $P_i$ are irreducible, positive on $\Gamma$, and counted with multiplicity. Applying the square free case to each factor yields $-D_bD_c\log P=-\sum_iD_bD_c\log P_i\ge0$.

Let $D_B=\prod_{i=1}^kD_{b_i}$ and $D_C=\prod_{j=1}^{m-k}D_{c_j}$. Derivatives in these directions of degree larger than 0 are hyperbolic and positive on $\Gamma$. For $Q=P$ or any such derivative of degree greater than $k$, successive applications of \eqref{rkt:rayleigh} in the directions $b_1,\ldots,b_k$ give
\[
 \frac{D_cD_BQ}{D_BQ}\le\frac{D_cQ}{Q}
 \qquad(c\in\Gamma).
\]
For $1\le j\le m-k$, let $c=c_j$ and $Q=D_{c_1}\cdots D_{c_{j-1}}P$. By multiplying out and cancelling intermediate terms we get
\begin{equation}\label{rkt:iterated}
 P(a)D_BD_CP(a)\le D_BP(a)D_CP(a),
\end{equation}
as in \cite[Section~2.2, (15)]{hu-xiao-rkt}. A derivative of order $j$ at $a$ is equal to $m!/(m-j)!$ times the polarization with $m-j$ copies of $a$. Substituting this into \eqref{rkt:iterated} gives \eqref{rkt:polarized}. By continuity this holds for $b_i,c_j\in\overline\Gamma$ as well.

Take $P = x_1\dots x_m$ and $a = (1,\dots,1)$ to prove the optimality of the result. Let $b_1,\dots,b_k,c_1,\dots,c_{m-k}$ be the standard basis of $\mathbb{R}^m$. Both sides of \eqref{rkt:iterated} equal 1 in this case. Approximating these basis vectors with vectors in the positive orthant proves optimality on $\Gamma$.
\end{proof}

\begin{corollary}[Pointwise sequential reverse inequality]\label{rkt:pointwise}
Let $1\le k\le m-1$ and $a,b_1,\ldots,b_k,c_1,\ldots,c_{m-k}\in\Gamma_\Theta$. Put
$B=b_1\wedge\cdots\wedge b_k$ and $C=c_1\wedge\cdots\wedge c_{m-k}$. Then
\begin{equation}\label{rkt:pointineq}
\begin{aligned}
 &\bigl(*(\Theta\wedge a^m)\bigr)\bigl(*(\Theta\wedge B\wedge C)\bigr)\\
 &\qquad\le\binom{m}{k}\bigl(*(\Theta\wedge a^{m-k}\wedge B)\bigr)
                 \bigl(*(\Theta\wedge a^k\wedge C)\bigr).
\end{aligned}
\end{equation}
For $\Theta=\omega^{d-m}$, this applies to ordinary $m$-positive forms.
\end{corollary}
\begin{proof}
Now we apply Lemma \ref{rkt:hyperbolic} to $P_\Theta(A)=*(\Theta\wedge A^m)$, hyperbolic as in Proposition \ref{seq:hyperbolicity}. Polarization is defined by $(A_1,\ldots,A_m)\mapsto*(\Theta\wedge A_1\wedge\cdots\wedge A_m)$.
\end{proof}

To move to cohomology we use the normalization of Lehmann and Xiao \cite[Theorem 5.7]{lehmann-xiao}.

\begin{proposition}\label{rkt:conditional}
For each $1\le k\le m-1$, let $a, b_1, \ldots, b_k, c_1, \ldots, c_{m-k}$ be smooth, closed, admissible $(1,1)$ forms. Let $B=b_1\wedge\ldots\wedge b_k$, $C=c_1\wedge\ldots\wedge c_{m-k}$, $A = \int_X \Theta\wedge a^m$ and $D = \int_X \Theta\wedge B\wedge C$. Assume there is a real valued smooth function $u$ on $X$ with the property
\begin{equation}\label{rkt:normalizing}
 \Theta\wedge a_u^m=\frac{A}{D}\Theta\wedge B\wedge C,
 \qquad
 a_u=a+\sqrt{-1}\partial\bar\partial u\in\Gamma_\Theta.
\end{equation}
Then
\begin{equation}\label{rkt:sharpglobal}
 \left(\int_X\Theta\wedge a^m\right)\left(\int_X\Theta\wedge B\wedge C\right)
 \le\binom{m}{k}
 \left(\int_X\Theta\wedge a^{m-k}\wedge B\right)
 \left(\int_X\Theta\wedge a^k\wedge C\right).
\end{equation}
\end{proposition}
\begin{proof}
Set $\mu=\Theta\wedge B\wedge C$ and write
\[
 \Theta\wedge a_u^{m-k}\wedge B=f\mu,\qquad
 \Theta\wedge a_u^k\wedge C=g\mu.
\]
Since the volume forms are positive according to G\aa rding’s inequality, substituting \eqref{rkt:normalizing} into \eqref{rkt:pointineq} yields $fg\ge A/(D\binom{m}{k})$. Applying the Cauchy--Schwarz inequality and using $\int_X\mu=D$ gives
\[
 \left(\int_X f\mu\right)\left(\int_X g\mu\right)
 \ge\left(\int_X\sqrt{fg}\,\mu\right)^2
 \ge\frac{AD}{\binom{m}{k}}.
\]
Since $[a_u]=[a]$ and $\Theta,B,C$ are closed, we can replace $a_u$ by $a$ in both integrals.
\end{proof}

\begin{proof}[Proof of Proposition~\ref{intro:rkt-sharp}]
Represent each of the given classes by admissible forms. The form $(A/D)\Theta\wedge B\wedge C$ is a positive volume form of integral $A$. Conjecture~\ref{intro:solvability} therefore yields an admissible solution of \eqref{rkt:normalizing}, and Proposition~\ref{rkt:conditional} applies.
\end{proof}

\begin{corollary}[Sharp ordinary $m$-positive reverse inequality]\label{rkt:ordinary}
For $\Theta=\omega^{d-m}$, the inequality~\eqref{rkt:sharpglobal} is satisfied by smooth closed $m$-positive forms $a,b_i,c_j$ with respect to $\omega$ and classes in $\overline{\mathcal K_m(X;[\omega])}$.
\end{corollary}
\begin{proof}
An admissible solution to equation \eqref{rkt:normalizing} exists by \cite[Corollary 10.2]{flzz}. Applying Proposition \ref{rkt:conditional} to $m$-positive forms, we may extend the inequality to the closure of $\mathcal K_m(X;[\omega])$.
\end{proof}

\begin{remark}[An unconditional cohomological inequality]\label{rkt:global}
Let $1\le k\le m-1$. For $a, b_i, c_j \in \mathcal A_\Theta(X)$, using the Hodge index argument of Hu and Xiao \cite[Theorem 2.5 and Corollary 3.10]{hu-xiao-rkt} we have
\begin{align}\label{rkt:globalineq}
 &I_\Theta(a^m)I_\Theta(b_1,\ldots,b_k,c_1,\ldots,c_{m-k})\notag\\
 &\qquad\le2^{k(m-k)}
 I_\Theta(a^{m-k},b_1,\ldots,b_k)
 I_\Theta(a^k,c_1,\ldots,c_{m-k}).
\end{align}
Indeed, consider any mixed pairing $Q$ as in Proposition \ref{rkt:package}(a). Let $b=sa+b_0$ and $c=ta+c_0$ where $b_0, c_0 \in a^\perp$. Since $Q(b,b)>0$ and $Q(c,c)>0$, we have $-Q(b_0,b_0)\le s^2Q(a,a)$ and $-Q(c_0,c_0)\le t^2Q(a,a)$. Mixed positivity implies $s, t > 0$. From the Cauchy--Schwarz inequality for $-Q$ restricted to $a^\perp$ we have
\[
 Q(a,a)Q(b,c)\le 2Q(a,b)Q(a,c).
\]
Let $I_{i,j}=I_\Theta(a^{m-i-j},b_1,\ldots,b_i,c_1,\ldots,c_j)$. Replacing $b$ and $c$ in the inequality with $b=b_{i+1}$ and $c=c_{j+1}$ with respect to the pairing $Q$ defined by $a^{m-i-j-2},b_1,\ldots,b_i,c_1,\ldots,c_j$ we get
$I_{i,j}I_{i+1,j+1}\le2I_{i+1,j}I_{i,j+1}$.
Multiplying these for $0\le i<k$ and $0\le j<m-k$ and cancelling terms gives the global inequality~\eqref{rkt:globalineq}. By continuity this also holds on the closure of $\mathcal A_\Theta(X)$.
\end{remark}

\begin{remark}[Correlation inequalities for volume polynomials]\label{rkt:correlation}
For projective and K\"ahler volume polynomials we use the multiset cover inequalities in \cite[Theorem 1.2 and Corollary 1.3]{huh-michalek-wang-xiao} on blocks of size $k$ and $m-k$ to obtain a reverse Khovanskii--Teissier inequality with coefficient $\binom{m}{k}$. In the projective case, we can additionally use a common product of auxiliary nef divisor classes \cite[Corollary 1.4]{huh-michalek-wang-xiao}. To use their result for any other classes we would first have to write the intersection polynomial as a volume polynomial.
\end{remark}

\appendix
\section{Volume hyperbolicity and higher Lefschetz maps}\label{app:volume}

We construct a smooth projective variety $X$ with hyperbolic volume polynomial $F_X$ and a class $\alpha$ in its hyperbolicity cone such that
$\alpha\cdot-\colon H^{3,3}(X)\to H^{4,4}(X)$ has a nonzero kernel.
Since divisor classes generate $H^*(X,\mathbb{R})$ in this example by \eqref{vol:ring}, the apolar algebra of $F_X$ is isomorphic to $H^*(X,\mathbb{R})$ via Poincaré duality. Consequently, the comparison of cones stated in Theorem \ref{main:blocks} may not hold for the apolar algebra of a volume polynomial.

Let $h=c_1(\mathcal O_{\mathbb P^4}(1))$. The rank two Horrocks--Mumford bundle $H$ on $\mathbb P^4$ has $c_1(H)=5h$ and $c_2(H)=10h^2$ (see \cite{horrocks-mumford}). Let $\varphi\colon\mathbb P^4 \rightarrow \mathbb P^4$ and let $G$ be the bundle given by
\[
 \varphi([z_0:\cdots:z_4])=[z_0^4:\cdots:z_4^4],
 \qquad G=\varphi^*(H(-3))\otimes\mathcal O(2).
\]
Using the twisting formula for rank two bundles and the fact that $\varphi^*h=4h$ we have
\[
 (c_1,c_2)(H(-3))=(-h,4h^2),\quad
 (c_1,c_2)(\varphi^*(H(-3)))=(-4h,64h^2),\quad
 (c_1,c_2)(G)=(0,60h^2).
\]
Let $X$ be the projective bundle of one-dimensional quotients
\[
 X=\mathbb P_{\mathbb P^4}\bigl(G\oplus\mathcal O(5)\oplus\mathcal O(-5)\bigr),
 \qquad \xi=c_1(\mathcal O_X(1)).
\]
We will write $h$ for its pullback to $X$ as well. Using the projective bundle formula, we get
\begin{equation}\label{vol:ring}
 H^*(X,\R)=
 \frac{\R[h,\xi]}{\bigl(h^5,(\xi^2+60h^2)\cdot(\xi^2-25h^2)\bigr)},
 \qquad\int_Xh^4\cdot\xi^3=1.
\end{equation}
In particular, $H^{p,q}(X)=0$ for $p\ne q$. By \eqref{vol:ring},
\[
 F_X(a,b)=\int_X(ah+b\xi)^7
       =b^3(35a^4-735a^2b^2+2725b^4)
       =35b^3(a^2-r_-b^2)(a^2-r_+b^2),
\]
where $r_\pm=(147\pm\sqrt{6349})/14$. Since $0<r_-<r_+<25$, all linear factors of $F_X$ are real. They span the dual of $H^{1,1}(X,\mathbb{R})$ and so $F_X$ is complete. The hyperbolicity cone of $F_X$ containing $\xi + Nh$ for $N \gg 0$ is
\[
 \mathcal C_X=\{(a,b):b>0,\ a>\sqrt{r_+}\,b\}.
\]
The ample cone is contained in $\mathcal C_X$ since it is connected, contains $\xi+Nh$, and $F_X$ is positive on it.

Set
\[
 \alpha=\xi+5h,\qquad
 \beta=(\xi-5h)\cdot(\xi^2+60h^2).
\]
Since $5>\sqrt{r_+}$, we have $\alpha\in\mathcal C_X$. From \eqref{vol:ring}, we have $\alpha\cdot\beta=0$ and $\int_X\beta\cdot h^4=1$, so $\beta$ is a nonzero element in the kernel of the map
\[
 \alpha\cdot-\colon H^{3,3}(X)\longrightarrow H^{4,4}(X).
\]
In general, using the basis $(h^3,h^2\cdot\xi,h\cdot\xi^2,\xi^3)$ of $H^{3,3}(X)$, the matrix of the pairing $(u,v)\mapsto\int_X(ah+b\xi)\cdot u\cdot v$ is
\[
 M(a,b)=
 \begin{pmatrix}
 0&0&b&a\\
 0&b&a&-35b\\
 b&a&-35b&-35a\\
 a&-35b&-35a&2725b
 \end{pmatrix},
 \qquad
 \det M(a,b)=(a^2-25b^2)(a^2+60b^2).
\]
At $(a,b)=(5,1)$, the vector $(-300,60,-5,1)^T$ represents $\beta$.

To get a similar example with first factor ample we set $Y=X\times\mathbb P^1$. Let $e$ be an ample class on $X$ and let $u$ be the hyperplane class on $\mathbb P^1$. For simplicity we denote $\eta_1 = e + u$ and $\eta_2 = \alpha$ on $Y$ by suppressing the pullbacks. Then $F_Y(v,s) = 8s F_X(v)$ for any $v \in H^{1,1}(X,\mathbb{R})$ and $s \in \mathbb{R}$. Since $e,\alpha \in \mathcal C_X$, the polynomial
\[
 g(t)=F_Y(\eta_2+t\eta_1)=8tF_X(\alpha+te)
\]
has seven negative real roots counted with multiplicity and a simple root at 0. By Rolle's theorem, roots of $g'(t)$ are all negative as well. Note that $g'(t)=(D_{\eta_1}F_Y)(\eta_2+t\eta_1)$. Thus
\[
 \eta_2\in\Gamma(D_{\eta_1}F_Y,\eta_1).
\]
We can use $\eta_1$ for each subsequent derivative because each successive derivative cone contains the preceding cone. However, the pullback of $\beta$ is a non-zero element in the kernel of
\[
 \eta_1\cdot\eta_2\cdot-\colon H^{3,3}(Y)\longrightarrow H^{5,5}(Y).
\]

The classes $\eta_1,\eta_2$ do not have representatives that form a good prefix on $Y$. Let $S \subset X$ be the section of $X \to \mathbb P^4$ corresponding to the quotient $\mathcal O(-5)$. Then $S \simeq \mathbb P^4$ and $\xi |_S = -5h$. Let $Z= S \times \mathbb P^1$. Since $\eta_2$ is the pullback of $\xi + 5h$, we have $\eta_2 |_Z = 0$. For any K\"ahler class $\kappa$ on $Y$, we have
\[
 \int_Z\eta_1\cdot\eta_2\cdot\kappa^3=0.
\]
If representatives of $\eta_1, \eta_2$ formed a good prefix, then by Proposition~\ref{seq:restriction} their restrictions to $Z$ would form a good prefix as well. Their product would then be strictly weakly positive on $Z$, and the integral would be positive.

\end{document}